\documentclass[12pt]{amsart}
\usepackage{amssymb}
\usepackage{amsmath}
\usepackage{amsbsy}
\usepackage{amscd}
\usepackage{amsthm}
\usepackage{url}

\usepackage[dvipdfmx]{graphicx}
\usepackage{here}
\advance\textheight by \topskip
\makeatletter
\def\Bbb{\mathbb}

\newenvironment{pf*}[1]{\proof[#1]}{\endproof}
\renewcommand{\thesubsection}{\thesection(\@roman\c@subsection)}
\makeatother
\newtheorem{Theorem}[equation]{Theorem}
\newtheorem{Corollary}[equation]{Corollary}
\newtheorem{Lemma}[equation]{Lemma}
\newtheorem{Proposition}[equation]{Proposition}

\newtheorem{Furuta}{\thmref{thmMS}$'$}

\theoremstyle{definition}
\newtheorem{Definition}[equation]{Definition}

\newtheorem{example}{Example}

\makeatletter
\renewcommand\section{\@startsection{section}{1}%
  {\z@}{.7\linespacing\@plus\linespacing}{.5\linespacing}%
  {\reset@font\normalfont\bfseries\center}}
\makeatother

\theoremstyle{remark}
\newtheorem{Remark}[equation]{Remark}

\newtheorem*{Acknowledgements}{Acknowledgements}
\numberwithin{equation}{section}
\numberwithin{figure}{section}

\newcommand{\thmref}[1]{Theorem~\ref{#1}}
\newcommand{\secref}[1]{Section~\ref{#1}}
\newcommand{\lemref}[1]{Lemma~\ref{#1}}
\newcommand{\propref}[1]{Proposition~\ref{#1}}
\newcommand{\remref}[1]{Remark~\ref{#1}}

\newcommand{\corref}[1]{Corollary~\ref{#1}}

\newcommand{\figref}[1]{Figure~\ref{#1}}

\newcommand{\furutaref}[1]{\thmref{thmMS}$'$}

\newcommand{\Romnum}[1]{\expandafter\uppercase\expandafter{\romannumeral #1}} 
\newcommand{\C}{{\Bbb C}}
\newcommand{\Z}{{\Bbb Z}}

\newcommand{\CP}{\operatorname{\C P}}

\begin{document}
\title[On homotopy $K3$ surfaces]{On homotopy $K3$ surfaces constructed by two knots and their applications}
\author{Masatsuna~Tsuchiya}
\address{Department of mathematics, Gakushuin University, 5-1, Mejiro 1-chome, Toshima-ku, Tokyo, 171-8588, Japan}
\email{tsuchiya@math.gakushuin.ac.jp}

\begin{abstract}
Let $LHT$ be a left handed trefoil knot and $K$ be any knot. We define $M_n(K)$ to be the homology $3$-sphere which is represented by a simple link of $LHT$ and $LHT \sharp K$ with framings $0$ and $n$ respectively. Starting with this link, we construct homotopy $K3$ and spin rational homology $K3$ surfaces containing $M_n(K)$. Then we apply the adjunction inequality to show that if $n>2g^n_s(K)-2$, $M_n(K)$ does not bound any smooth spin rational $4$-ball, and that under the same assumption the negative $n$-twisted Whitehead double of $LHT \sharp K$ is not a slice knot, where $g^n_s(K)$ is the $n$-shake genus of $K$.
\end{abstract}
\maketitle
%
%
\section{Introduction}\label{sec:intro}
\begin{figure}[H]
 \centering
  \includegraphics[height=15mm]{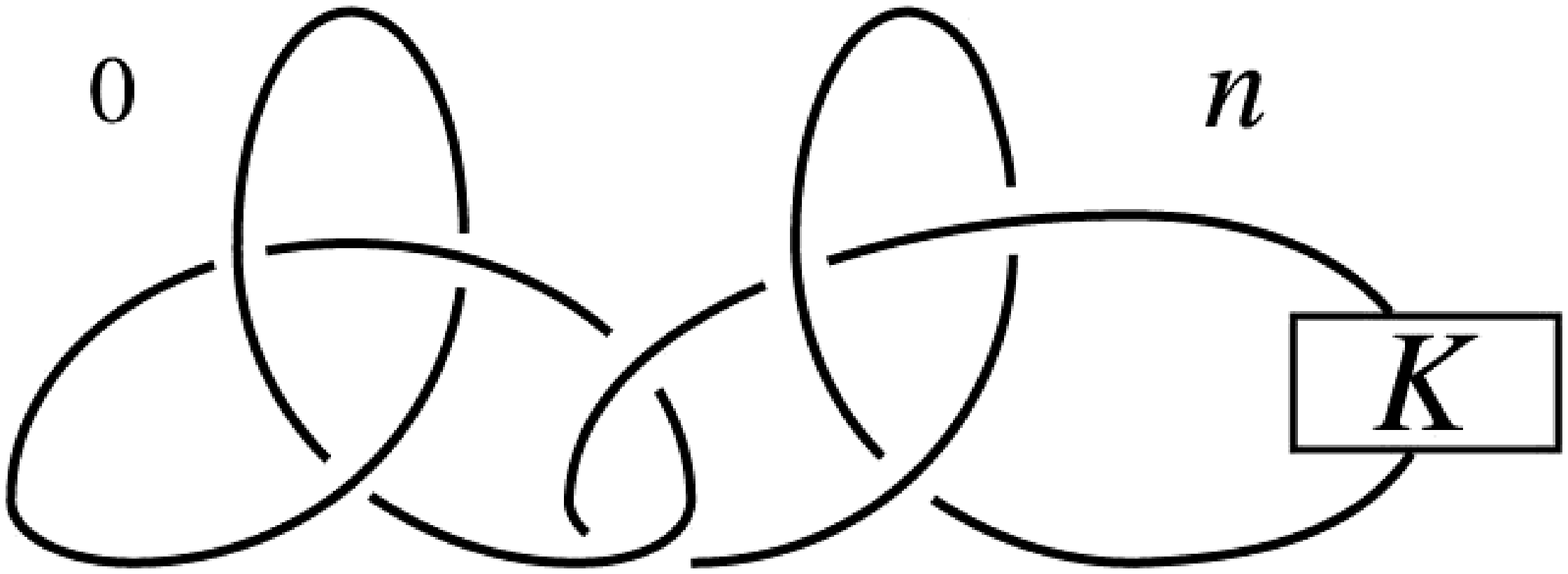}
 \caption{$X_n(K)$}
 \label{XnK}
\end{figure}

%
%
Let $K$ be a knot in $S^3$. We define $X_n (K)$ to be the $4$-dimensional handlebody which has a handle decomposition represented by \figref{XnK}, and define $M_n (K)$ to be $\partial (X_n (K))$. The boundary $M_n (K)$ is a homology $3$-sphere. Note that the right side knot of this link is the connected sum of a left handed trefoil knot $LHT$ and $K$.

\begin{Definition}[$r$-shake genus of $K$]
Let $N_{K, r}$ be a $4$-dimensional handlebody which is constructed by attaching a $2$-handle to $D^4$ along a knot $K$ with $r$-framing. We define the \textit{$r$-shake genus} of $K$ to be the minimal genus of smoothly embedded closed oriented surfaces in $N_{K, r}$ representing the generator of $H_2(N_{K, r})$. We denote the $r$-shake genus of $K$ by $g^r_s(K)$.
\end{Definition}

\begin{Remark}\label{shake}
Let $g_4(K)$  be the $4$-ball genus of $K$. Then we have $g^r_s(K) \leq g_4(K)$, for any $r \in \Z$.
\end{Remark}

In this paper, we show the following

\begin{Theorem}\label{thmt}
If $n>2g^n_s(K)-2$, $M_n(K)$ does not bound any smooth spin rational $4$-ball.
\end{Theorem}

If $K$ is an unknot, \thmref{thmt} is closely related to M.~Tange's result (see \remref{rem1}), who uses the Heegaard Floer homology $HF^+(M_n(U))$ and the correction term $d(M_n(U))$, while we apply the adjunction inequality.

\begin{Corollary}\label{cort}
If $n>2g^n_s(K)-2$, the negative $n$-twisted Whitehead double of $LHT \sharp K$ is not a slice knot, where $LHT \sharp K$ is the connected sum of $LHT$ and $K$.
\end{Corollary}

\begin{Remark}
By \corref{cort}, if $LHT \sharp K$ is a slice knot, $g^0_s(K)$ is not equal to $0$.
\end{Remark}

If $\tau(K)=g^n_s(K)$, \corref{cort} is a special case of Hedden's result (see \cite[Theorem~1.5]{H}), where $\tau(K)$ is the $\tau$-invariant of $K$.

We will prove \thmref{thmt} and \corref{cort} in \secref{sec:thm}.

\begin{Theorem}\label{thm:sec}
Let $(a_1, a_2, a_3, a_4)$ be any permutation of $(-1, -1, -2k-1,-2m-1)$, where $m$ and $k$ are integers and $k\geq0$. Then the knot represented by \figref{thmt-2} is not a slice knot, where $a_i$ is the number of full twists $(i=1, 2, 3, 4)$.
\begin{figure}[H]
 \centering
  \includegraphics[height=12mm]{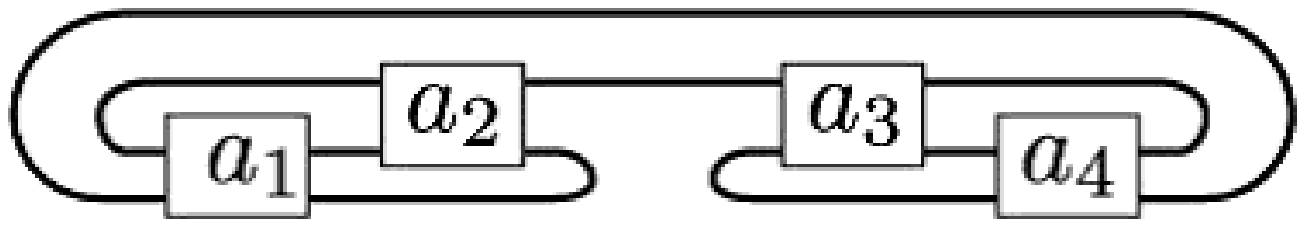}
\caption{}
 \label{thmt-2}
\end{figure}
\end{Theorem}

We will prove \thmref{thm:sec} in \secref{sec:3}.

\begin{Remark}\label{rem1}
Let $U$ be an unknot. Y.~Matsumoto asked in \cite[Problem~4.28]{K} whether $M_0(U)$ bounds a contractible $4$-manifold or not. By Gordon's result \cite{G}, if $n$ is odd, $M_n(U) $ does not bound any contractible $4$-manifold (cf.~\cite[\S3.1]{Y}). If $n$ is equal to $-6$, N.~Maruyama \cite{M} proved that $M_{-6}(U) $ bounds a contractible $4$-manifold. If $n$ is  equal to $0$, S.~Akbulut \cite{A} proved that $M_0(U) $ does not bound any contractible $4$-manifold. If $n>-2$, M.~Tange \cite{T} proved that $M_n(U)$ does not bound any negative definite $4$-manifold by computing the Heegaard Floer homology $HF^+(M_n(U))$ and the correction term $d(M_n(U))$.
\end{Remark}
\begin{Acknowledgements}
The author would like to thank Tetsuya Abe, Mikio Furuta, Yukio Matsumoto, Nobuhiro Nakamura and Motoo Tange for their useful comments and encouragement. 
\end{Acknowledgements}
%
%
\section{Proof of \thmref{thmt}}\label{sec:thm}
\noindent{\bf Notation.}
\begin{enumerate}
\item We represent a knot $K$ by \figref{knot}. If $K$ is a right handed trefoil knot $RHT$, $X_n(RHT)$ is represented by \figref{XnRHT}.
\begin{figure}[H]
\begin{minipage}{0.5\hsize}
  \centering
   \includegraphics[height=15mm]{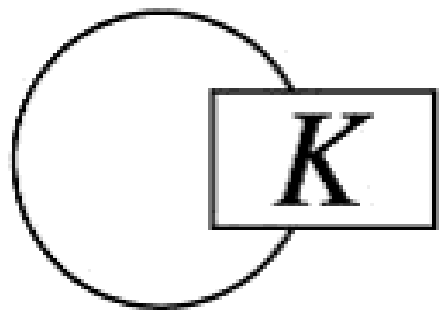}
   \caption{$K$}
   \label{knot}
 \end{minipage}%
 \begin{minipage}{0.5\hsize}
  \centering
   \includegraphics[height=15mm]{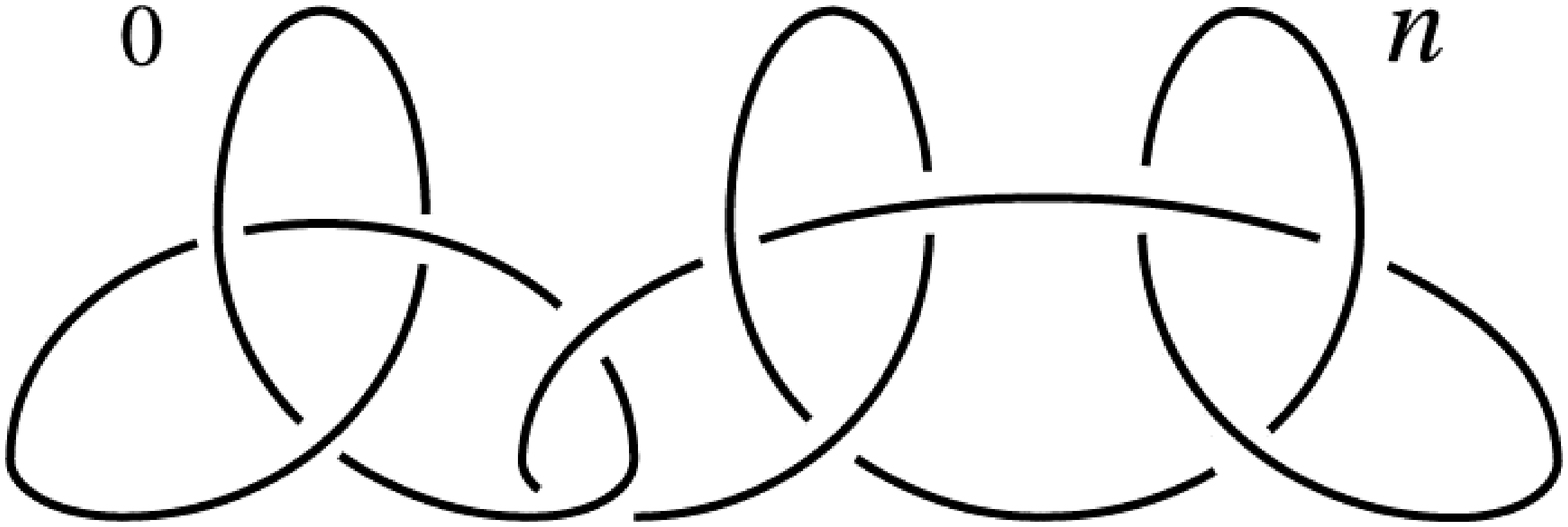}
  \caption{$X_n (RHT)$}
  \label{XnRHT}
 \end{minipage}
\end{figure}
\item Let $D_-(K, n)$ (resp.~$D_+(K,n)$) be the negative (resp.~positive) $n$-twisted Whitehead double of $K$, and we represent $D_-(K, n)$ by \figref{WD-1}. For example, $D_-(LHT, -6)$ is represented by \figref{WD-2}. To simplify the diagram, we usually use \figref{WD-3} instead of \figref{WD-2}.
\begin{figure}[H]
\begin{minipage}{0.32\hsize}
  \centering
   \includegraphics[height=19mm]{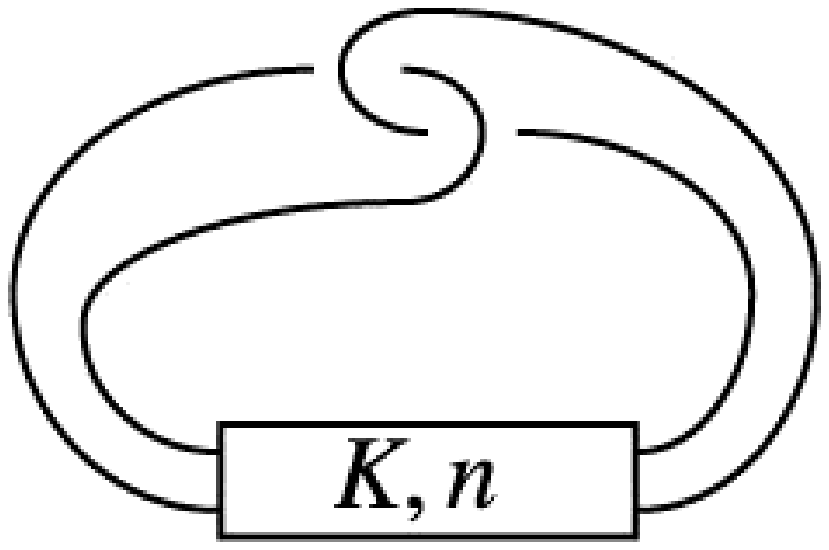}
   \caption{$D_-(K, n)$}
   \label{WD-1}
 \end{minipage}%
\begin{minipage}{0.34\hsize}
  \centering
   \includegraphics[height=22mm]{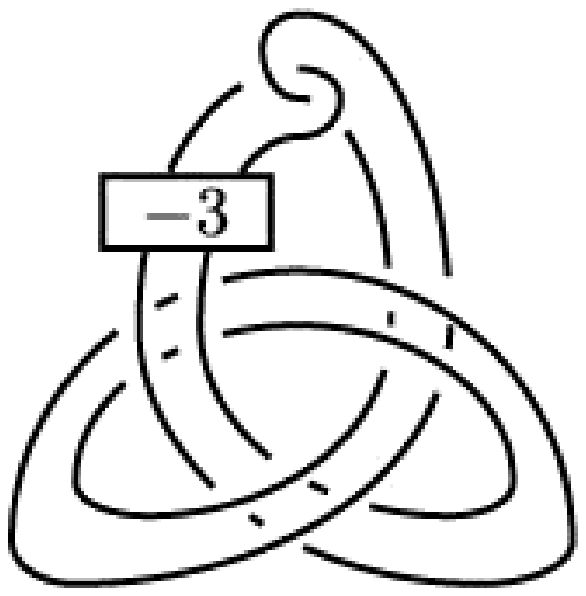}
   \caption{$D_-(LHT, -6)$}
   \label{WD-2}
 \end{minipage}%
 \begin{minipage}{0.34\hsize}
  \centering
   \includegraphics[height=19mm]{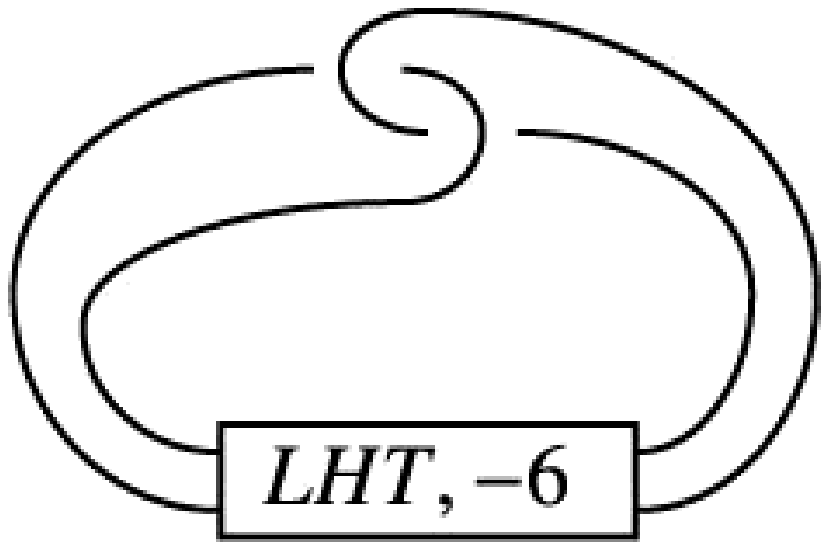}
  \caption{$D_-(LHT, -6)$}
  \label{WD-3}
 \end{minipage}
\end{figure}
\end{enumerate}

Let  $V_n(K)$ be the 4-dimensional handlebody represented by \figref{prop1-2}. 
\begin{Proposition}\label{prop1}
The $4$-dimensional handlebodies $X_n(K)$ and $V_n(K)$ have the same boundary.
\begin{figure}[H]
\begin{minipage}{0.4\hsize}
  \centering
   \includegraphics[height=16mm]{XnK.eps}
   \caption{$X_n (K)$}
   \label{prop1-1}
 \end{minipage}%
\begin{minipage}{0.2\hsize}
  \centering
   \includegraphics[height=10mm]{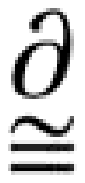}
 \end{minipage}%
 \begin{minipage}{0.4\hsize}
  \centering
   \includegraphics[height=20mm]{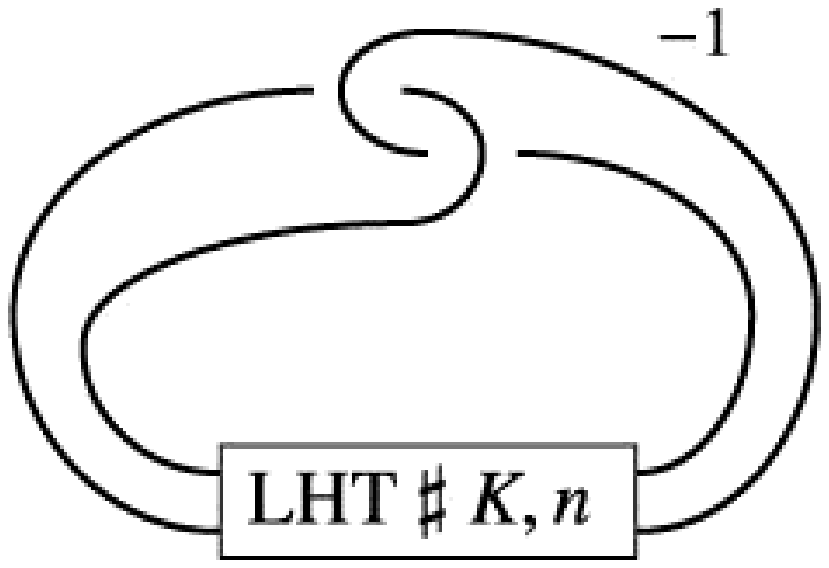}
  \caption{$V_n(K)$}
  \label{prop1-2}
 \end{minipage}
\end{figure}
\end{Proposition}
\proof
We show \propref{prop1} by the following handle calculus:
\begin{figure}[H]
 \begin{minipage}{0.35\hsize}
  \centering
   \includegraphics[height=17mm]{XnK.eps}
  \caption{$X_n$}
 \end{minipage}%
\begin{minipage}{0.3\hsize}
  \centering
   \includegraphics[height=12mm]{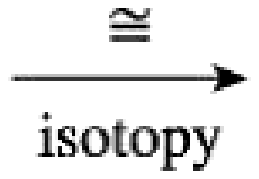}
 \end{minipage}%
 \begin{minipage}{0.35\hsize}
  \centering
   \includegraphics[height=19mm]{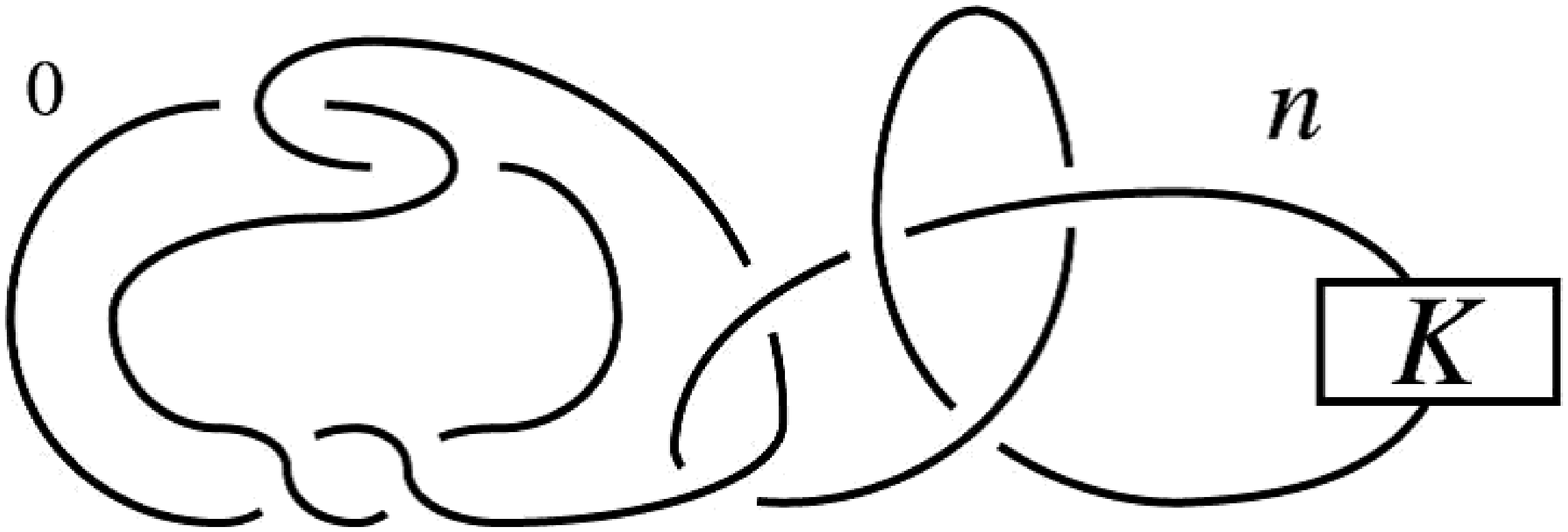}
  \caption{}
 \end{minipage}
\end{figure}

\begin{figure}[H]
\begin{minipage}{0.2\hsize}
  \centering
   \includegraphics[height=11mm]{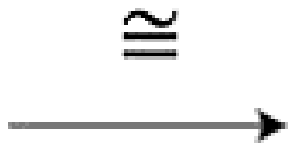}
 \end{minipage}%
 \begin{minipage}{0.3\hsize}
  \centering
   \includegraphics[height=30mm]{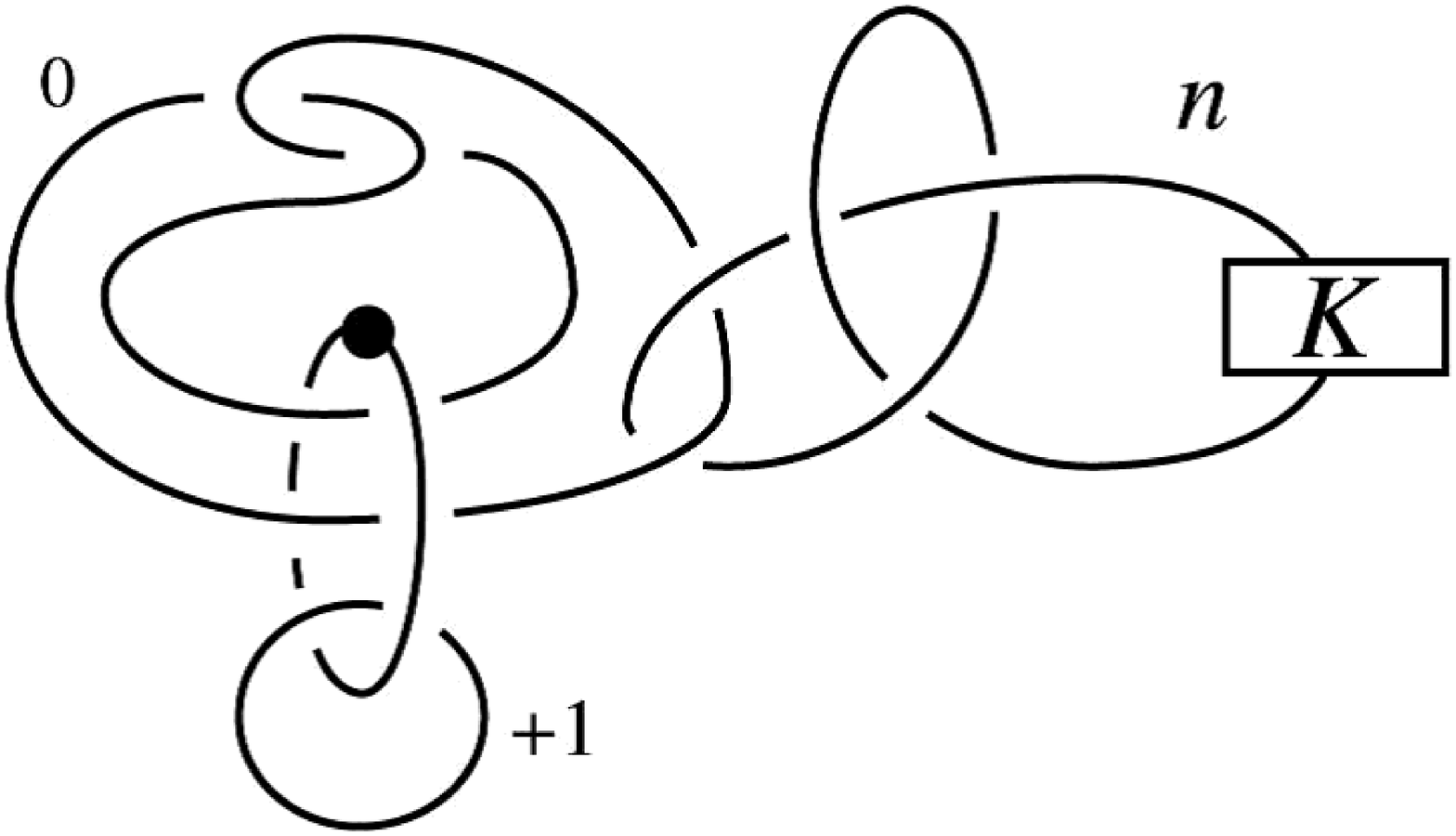}
   \caption{}
 \end{minipage}%
\begin{minipage}{0.2\hsize}
  \centering
   \includegraphics[height=12mm]{isotopy.eps}
 \end{minipage}%
 \begin{minipage}{0.3\hsize}
  \centering
   \includegraphics[height=30mm]{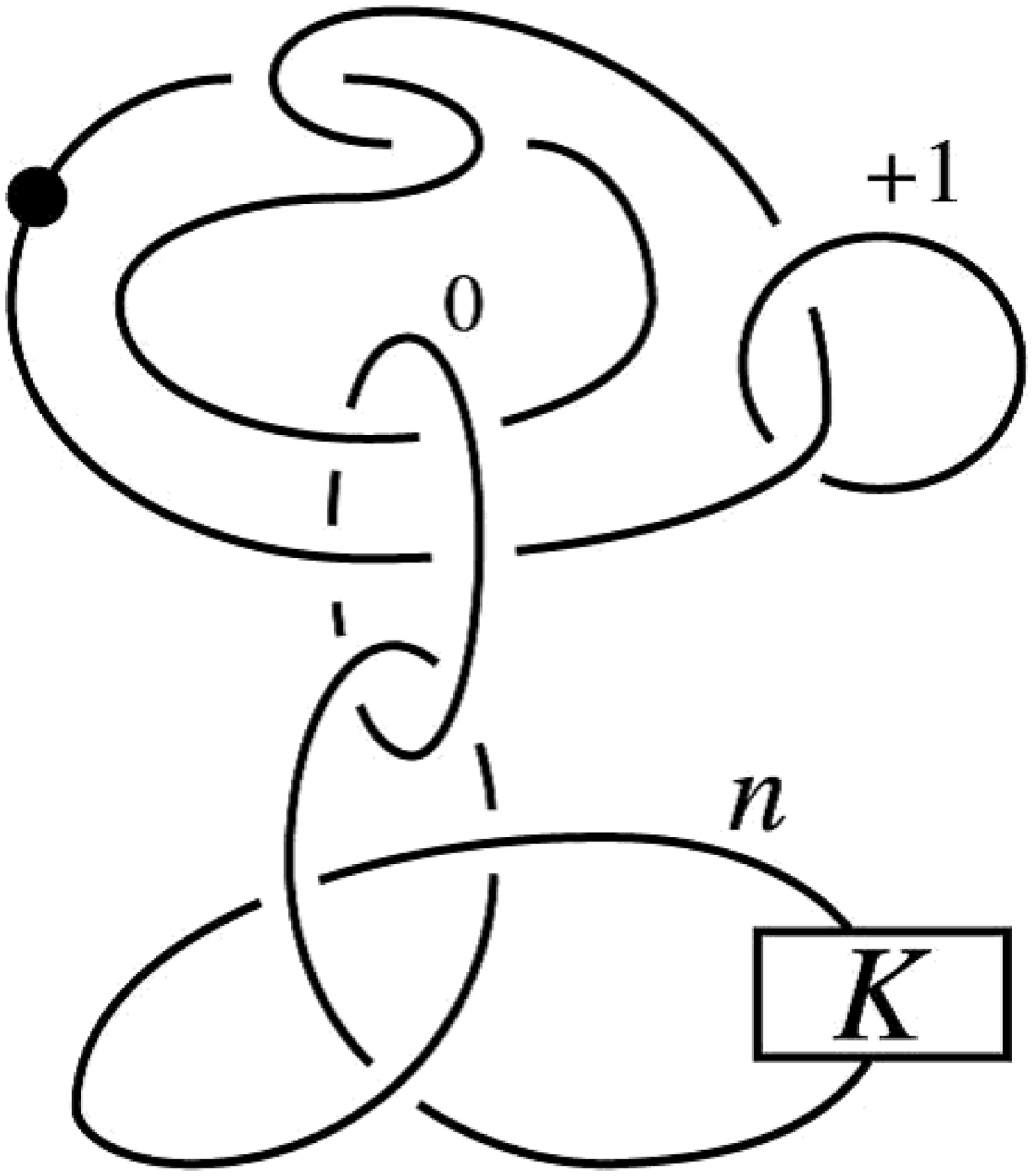}
  \caption{}
 \end{minipage}
\end{figure}

\begin{figure}[H]
\begin{minipage}{0.2\hsize}
  \centering
   \includegraphics[height=20mm]{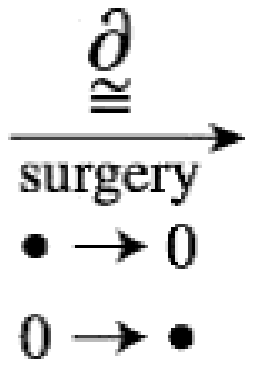}
 \end{minipage}%
 \begin{minipage}{0.3\hsize}
  \centering
   \includegraphics[height=32mm]{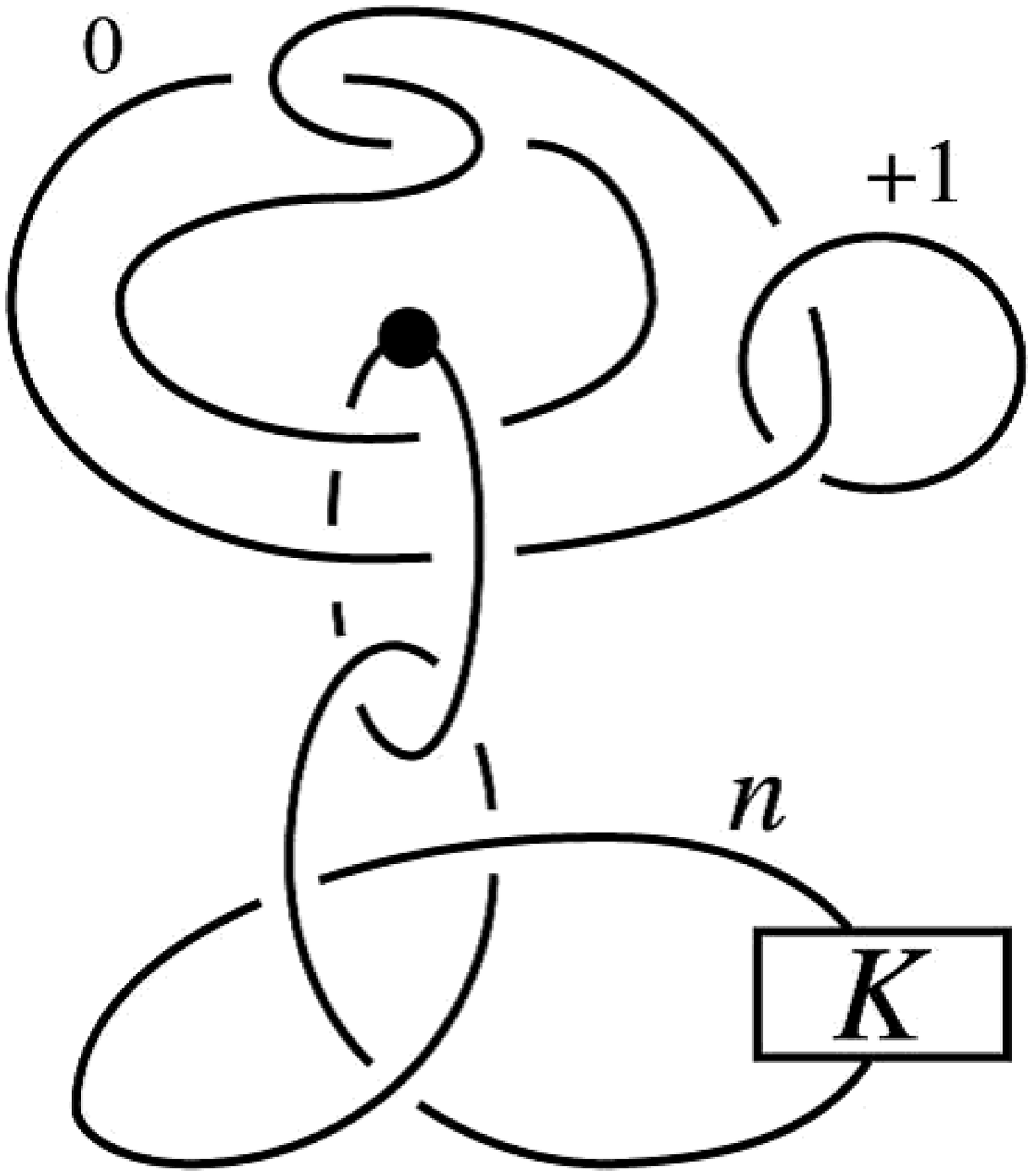}
   \caption{}
 \end{minipage}%
\begin{minipage}{0.2\hsize}
  \centering
   \includegraphics[height=10mm]{diff.eps}
 \end{minipage}%
 \begin{minipage}{0.3\hsize}
  \centering
   \includegraphics[height=20mm]{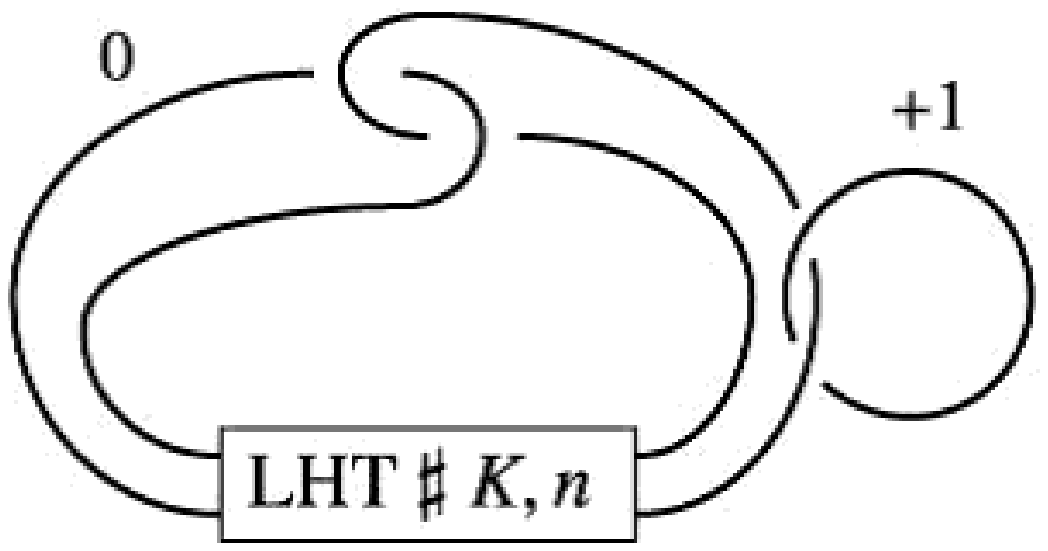}
  \caption{}
 \end{minipage}
\end{figure}

\begin{figure}[H]
\begin{minipage}{0.2\hsize}
  \centering
   \includegraphics[height=15mm]{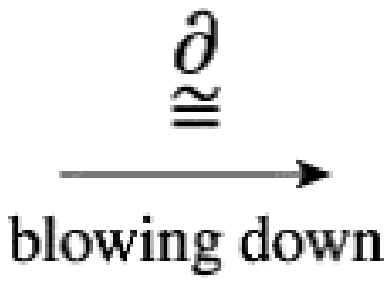}
 \end{minipage}%
 \begin{minipage}{0.3\hsize}
  \centering
   \includegraphics[height=21mm]{VnK.eps}
   \caption{$V_n(K)$}
 \end{minipage}
\end{figure}
\endproof

\begin{Remark}\label{rem2}
By \figref{prop1-2}, we can compute the Casson invariant $\lambda(M_n(K))$ as follows (cf.~\cite[\S2(vii)]{Tsu}):
\begin{center}
$\lambda(M_n(K))=-n.$
\end{center}
The Casson invariant, when reduced modulo 2, is the Rohlin invariant:
\begin{center}
$\lambda(M_n(K))\equiv\mu(M_n(K))\mod2.$
\end{center}
Therefore, if $n$ is odd, $M_n(K)$ does not bound any smooth spin rational $4$-ball. 
\end{Remark}

\begin{Remark}\label{rem3}
If $D_- (LHT \sharp K, n)$ is a slice knot, we have a smooth $S^2$ with self intersection $-1$ in $V_n(K)$ representing the generator of $H_2(V_n(K))$. We can blow down this $S^2$ to get a smooth contractible $4$-manifold $W_n(K)$. For example, Casson showed that $D_-(LHT, -6)$ is a slice knot. Therefore $V_{-6}(U)$ can be blown down to a contractible $4$-manifold $W_{-6}(U)$ represented by \figref{W-6U}.
\begin{figure}[H]
\begin{minipage}{0.38\hsize}
 \centering
   \includegraphics[height=27mm]{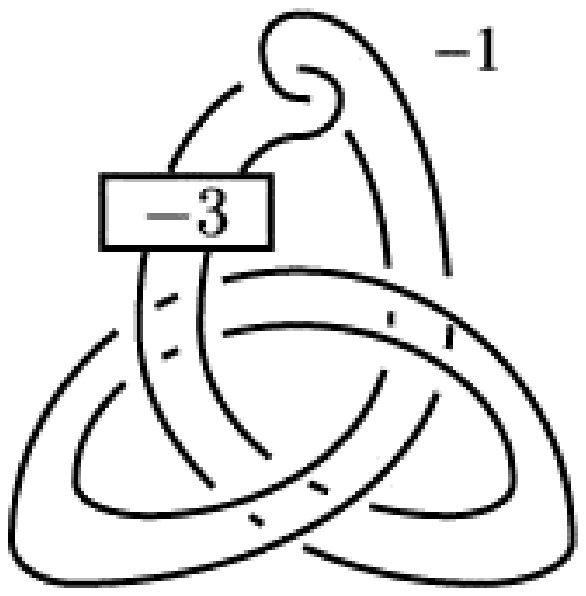}
   \caption{$V_{-6}(U)$}
 \end{minipage}%
\begin{minipage}{0.18\hsize}
  \centering
   \includegraphics[height=10mm]{bdydiff.eps}
 \end{minipage}%
 \begin{minipage}{0.44\hsize}
  \centering
  \includegraphics[height=27mm]{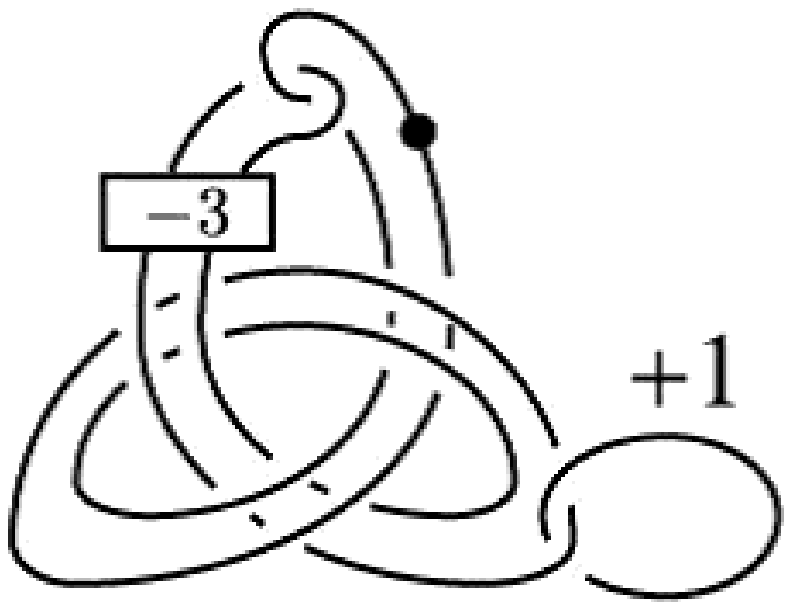}
 \caption{$W_{-6}(U)$}
 \label{W-6U}
 \end{minipage}
\end{figure}
\end{Remark}

\begin{figure}[H]
 \centering
  \includegraphics[height=81mm]{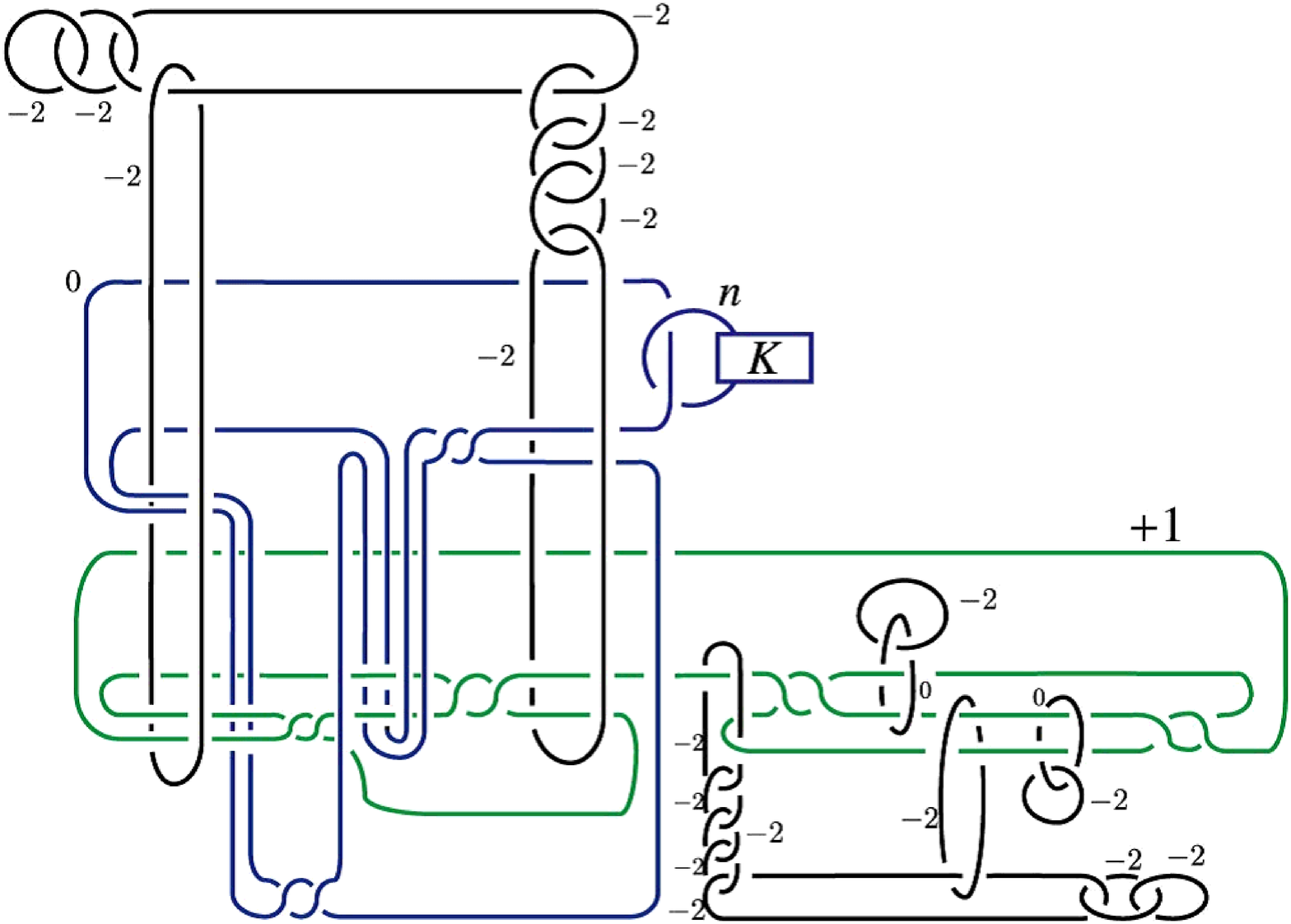}
 \caption{$Y_n (K)$}
 \label{YnK}
\end{figure}

Let $Y_n(K)$ be the $4$-dimensional handlebody represented by \figref{YnK} with intersection from $2E_8\oplus2\displaystyle\bigl(\begin{smallmatrix} 0 & 1 \\ 1 & 0 \end{smallmatrix} \bigr)\oplus\displaystyle\bigl(\begin{smallmatrix} 0 & 1 \\ 1 & n \end{smallmatrix} \bigr)\oplus1$.

\begin{figure}[H]
 \centering
  \includegraphics[height=81mm]{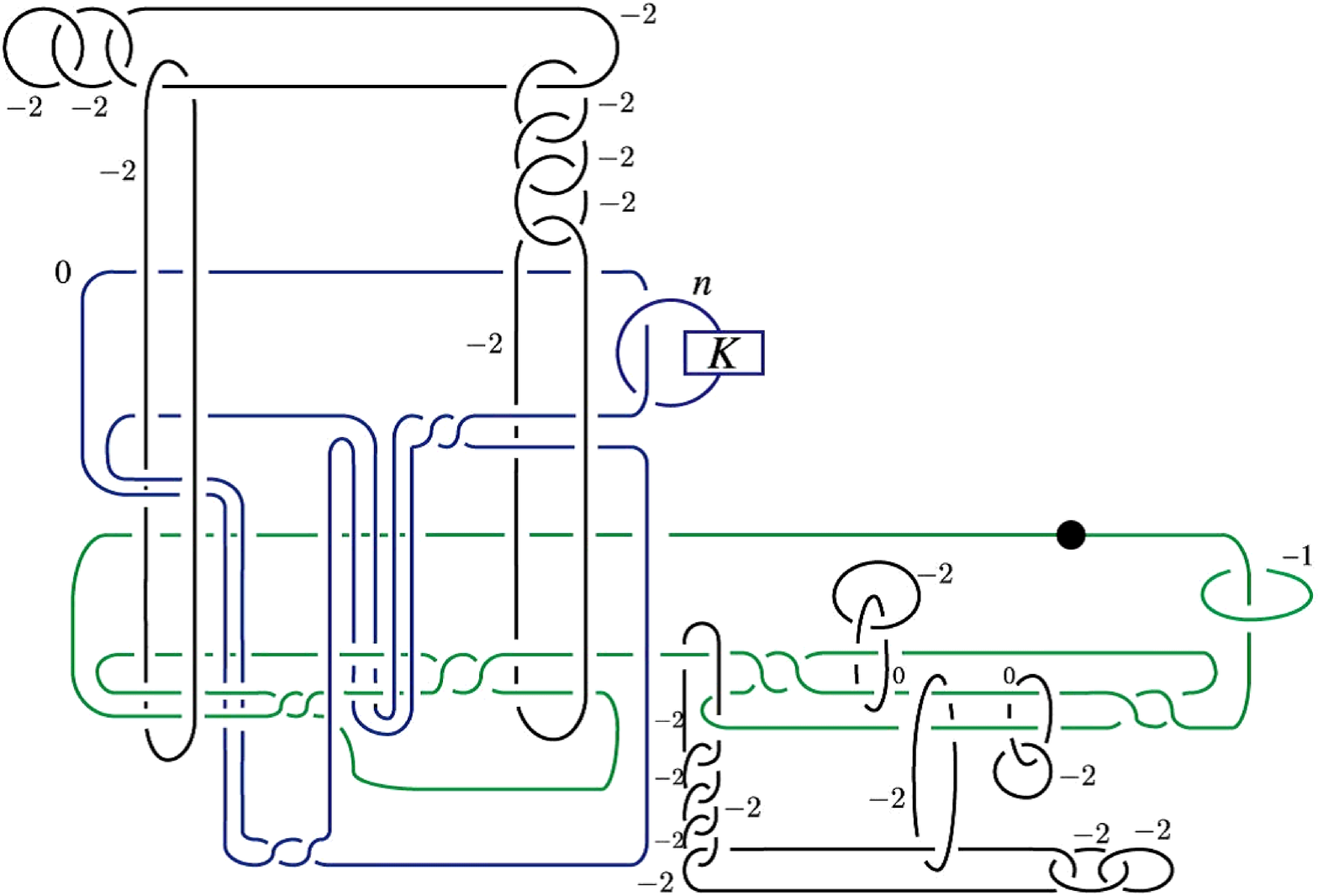}
 \caption{$Z_n (K)$}
 \label{ZnK}
\end{figure}

Because the $+1$-framed knot in \figref{YnK} is a slice knot ($RHT \sharp LHT$), we can blow it down. Then we have a smooth simply connected $4$-dimensional handlebody $Z_n(K)$ represented by \figref{ZnK} with intersection form  $2E_8\oplus2\displaystyle\bigl(\begin{smallmatrix} 0 & 1 \\ 1 & 0 \end{smallmatrix} \bigr)\oplus\displaystyle\bigl(\begin{smallmatrix} 0 & 1 \\ 1 & n \end{smallmatrix} \bigr)$. Note that $Y_n(K)$ and $Z_n(K)$ have the same boundary.

\begin{Proposition}\label{prop2}
The  $4$-dimensional handlebodies $X_n(K)$ and $Z_n(K)$ have the same boundary.
\end{Proposition}
\proof
Because $Y_n(K)$ and $Z_n(K)$ have the same boundary, we show that $X_n(K)$ and $Y_n(K)$ have the same boundary by \textquotedblleft S.~Akbulut's blowing up down process\textquotedblright (see \cite[Figures~1--4]{A} and \cite[Figures~9--19]{A1}) as follows:
\begin{figure}[H]
 \begin{minipage}{0.35\hsize}
  \centering
   \includegraphics[height=20mm]{XnK.eps}
  \caption{$X_n$}
  \label{prop2-1}
 \end{minipage}%
 \begin{minipage}{0.3\hsize}
  \centering
   \includegraphics[height=15mm]{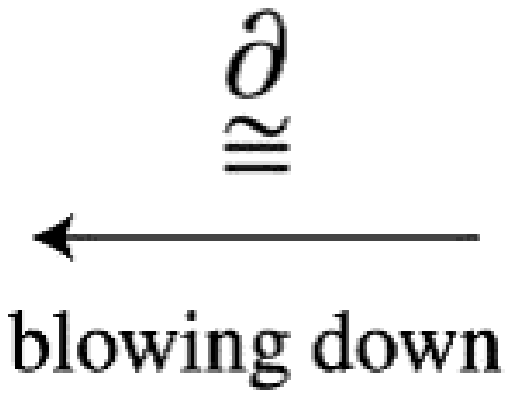}
 \end{minipage}%
 \begin{minipage}{0.3\hsize}
  \centering
   \includegraphics[height=25mm]{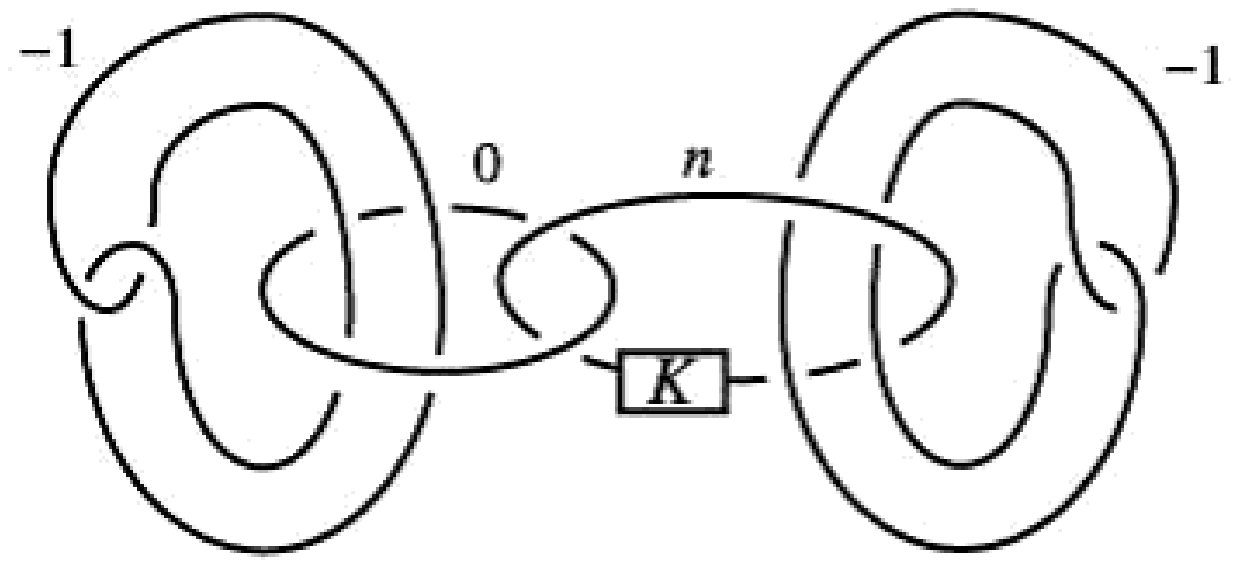}
  \caption{}
 \end{minipage}
\end{figure}  
\begin{figure}[H]
\begin{minipage}{0.15\hsize}
  \centering
   \includegraphics[height=17mm]{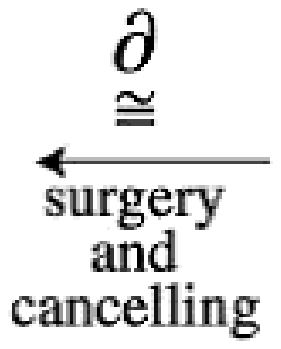}
 \end{minipage}%
 \begin{minipage}{0.39\hsize}
  \centering
   \includegraphics[height=30mm]{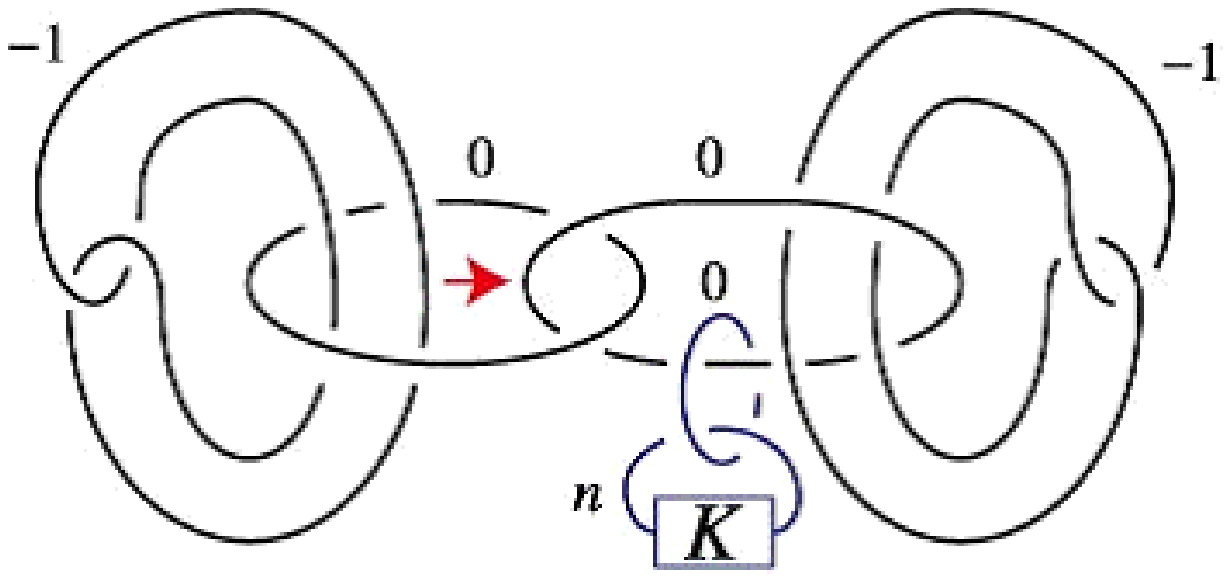}
   \caption{}
 \end{minipage}%
\begin{minipage}{0.16\hsize}
  \centering
   \includegraphics[height=17mm]{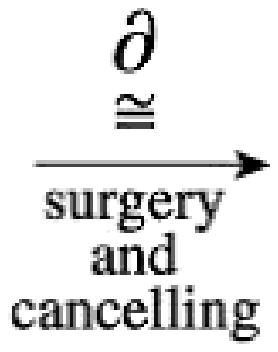}
 \end{minipage}%
 \begin{minipage}{0.3\hsize}
  \centering
   \includegraphics[height=35mm]{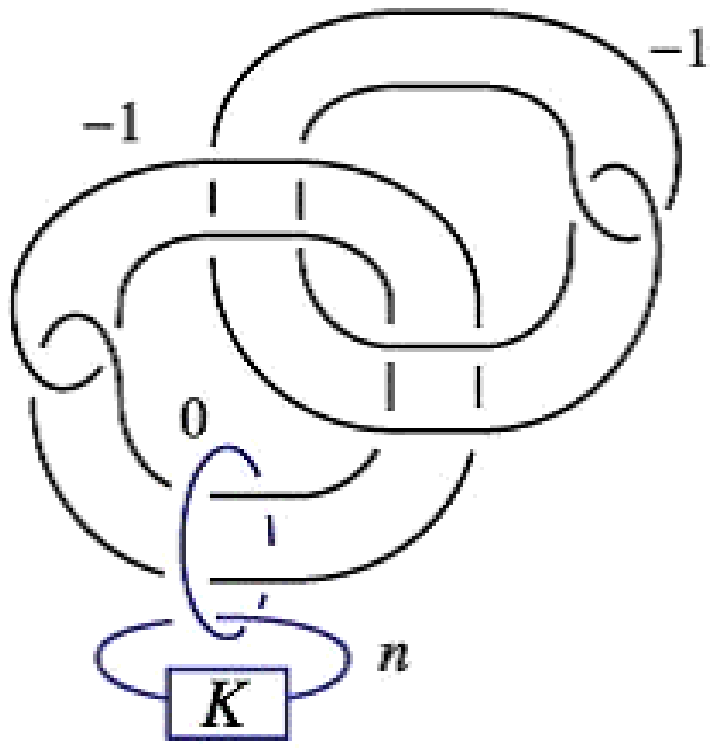}
  \caption{}
 \end{minipage}
\end{figure}  
\begin{figure}[H]
\begin{minipage}{0.2\hsize}
  \begin{center}
   \includegraphics[height=15mm]{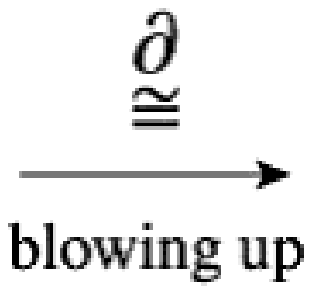}
  \end{center}
 \end{minipage}%
 \begin{minipage}{0.3\hsize}
  \begin{center}
   \includegraphics[height=30mm]{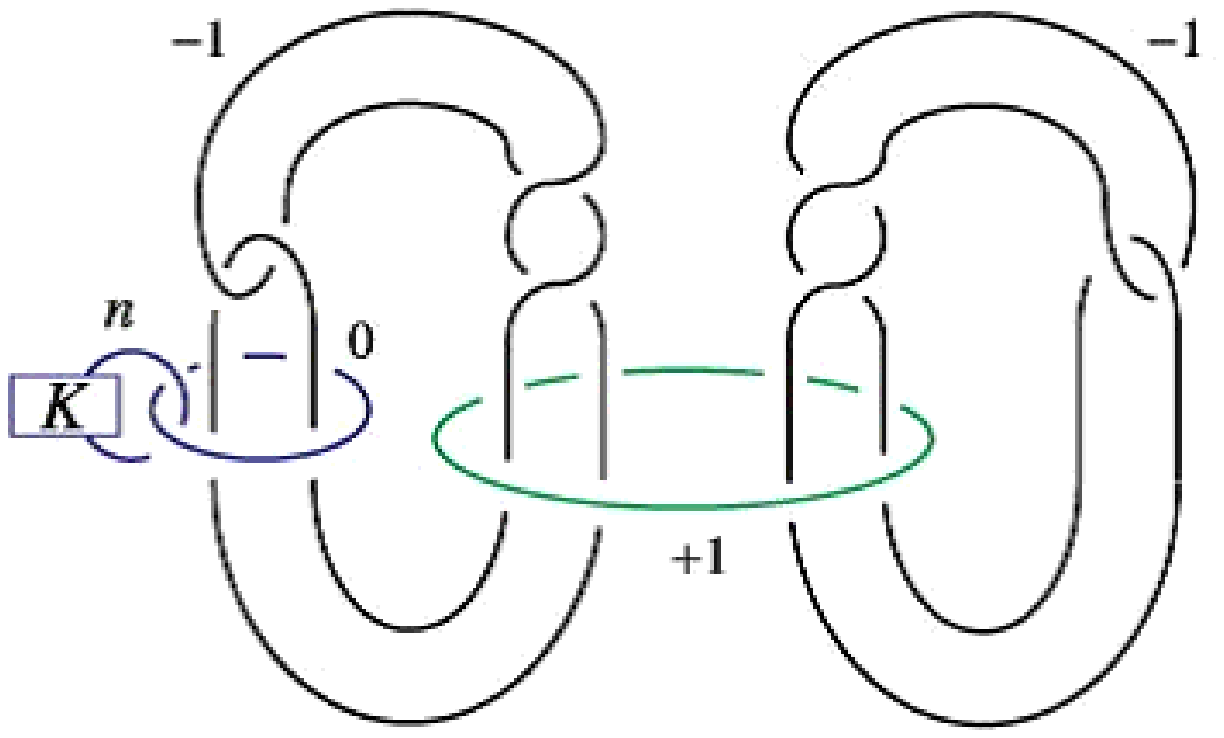}
  \end{center}
   \caption{}
 \end{minipage}%
\begin{minipage}{0.19\hsize}
  \begin{center}
   \includegraphics[height=11mm]{isotopy.eps}
  \end{center}
 \end{minipage}%
 \begin{minipage}{0.3\hsize}
  \begin{center}
   \includegraphics[height=35mm]{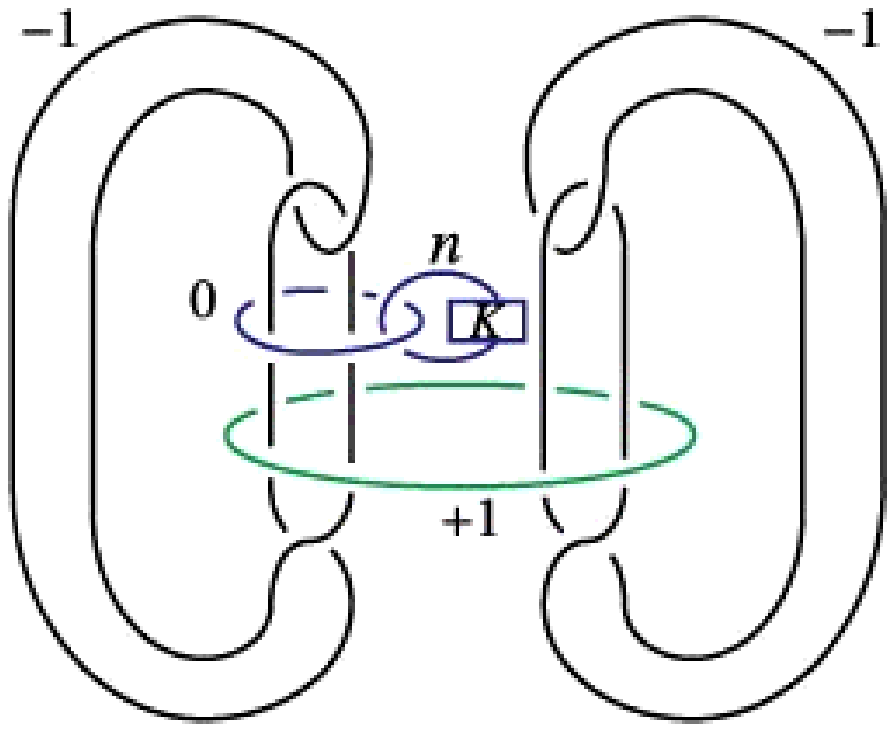}
  \end{center}
  \caption{}
 \end{minipage}
\end{figure}
\begin{figure}[H]
\begin{minipage}{0.15\hsize}
  \begin{center}
   \includegraphics[height=11mm]{isotopy.eps}
  \end{center}
 \end{minipage}%
 \begin{minipage}{0.35\hsize}
  \begin{center}
   \includegraphics[height=30mm]{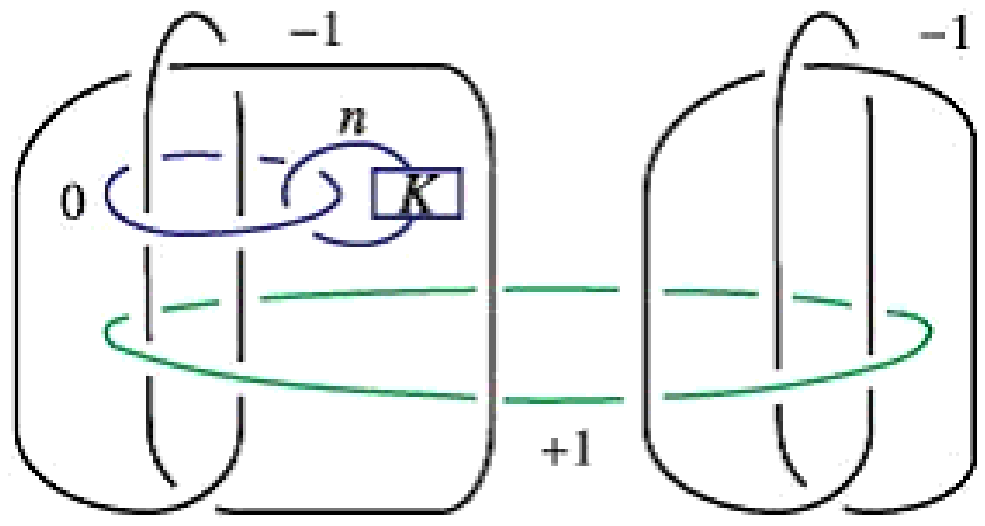}
  \end{center}
   \caption{}
 \end{minipage}%
\begin{minipage}{0.15\hsize}
  \begin{center}
   \includegraphics[height=15mm]{blowup.eps}
  \end{center}
 \end{minipage}%
 \begin{minipage}{0.35\hsize}
  \begin{center}
   \includegraphics[height=32mm]{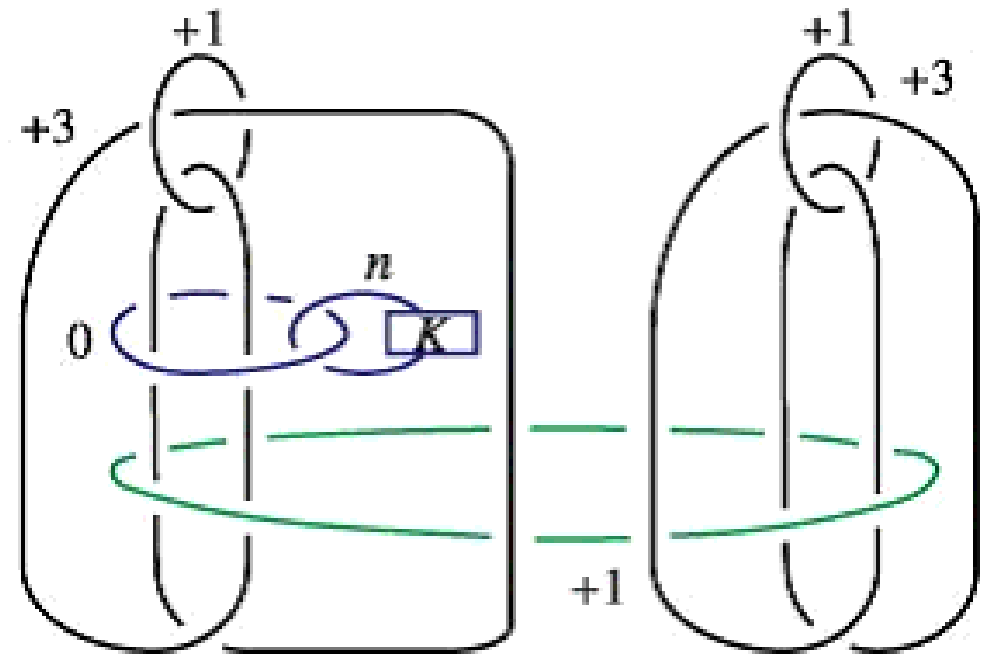}
  \end{center}
  \caption{}
 \end{minipage}
\end{figure}
\begin{figure}[H]
\begin{minipage}{0.15\hsize}
  \begin{center}
   \includegraphics[height=10mm]{isotopy.eps}
  \end{center}
 \end{minipage}%
 \begin{minipage}{0.35\hsize}
  \begin{center}
   \includegraphics[height=30mm]{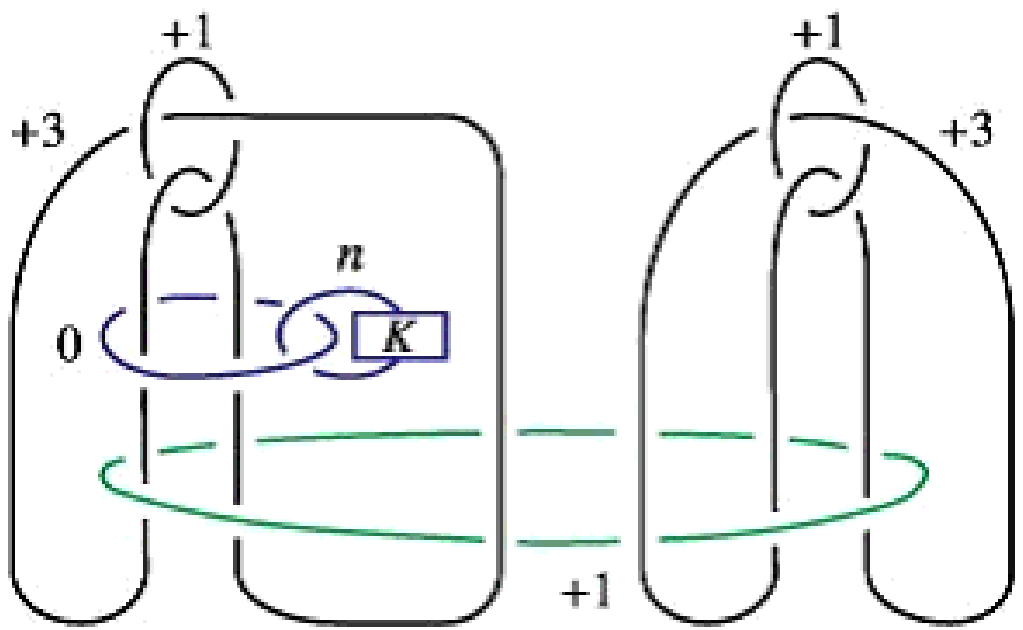}
  \end{center}
   \caption{}
 \end{minipage}%
\begin{minipage}{0.1\hsize}
  \begin{center}
   \includegraphics[height=10mm]{isotopy.eps}
  \end{center}
 \end{minipage}%
 \begin{minipage}{0.4\hsize}
  \begin{center}
   \includegraphics[height=30mm]{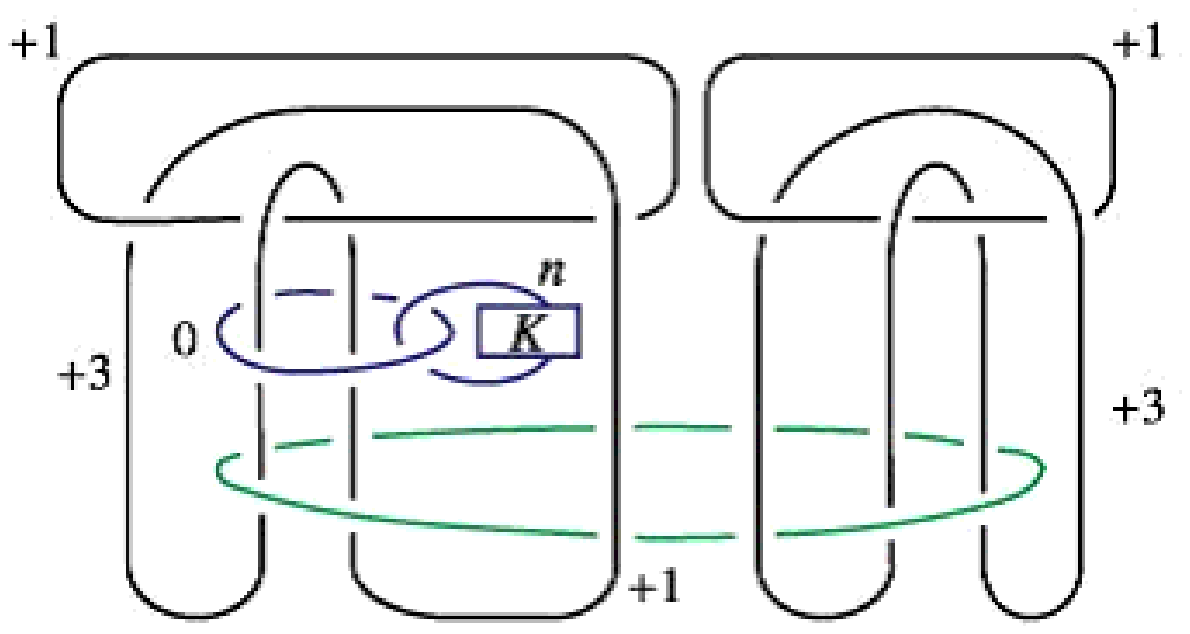}
  \end{center}
  \caption{}
 \end{minipage}
\end{figure}
\begin{figure}[H]
\begin{minipage}{0.1\hsize}
  \begin{center}
   \includegraphics[height=10mm]{isotopy.eps}
  \end{center}
 \end{minipage}%
 \begin{minipage}{0.35\hsize}
  \begin{center}
   \includegraphics[height=30mm]{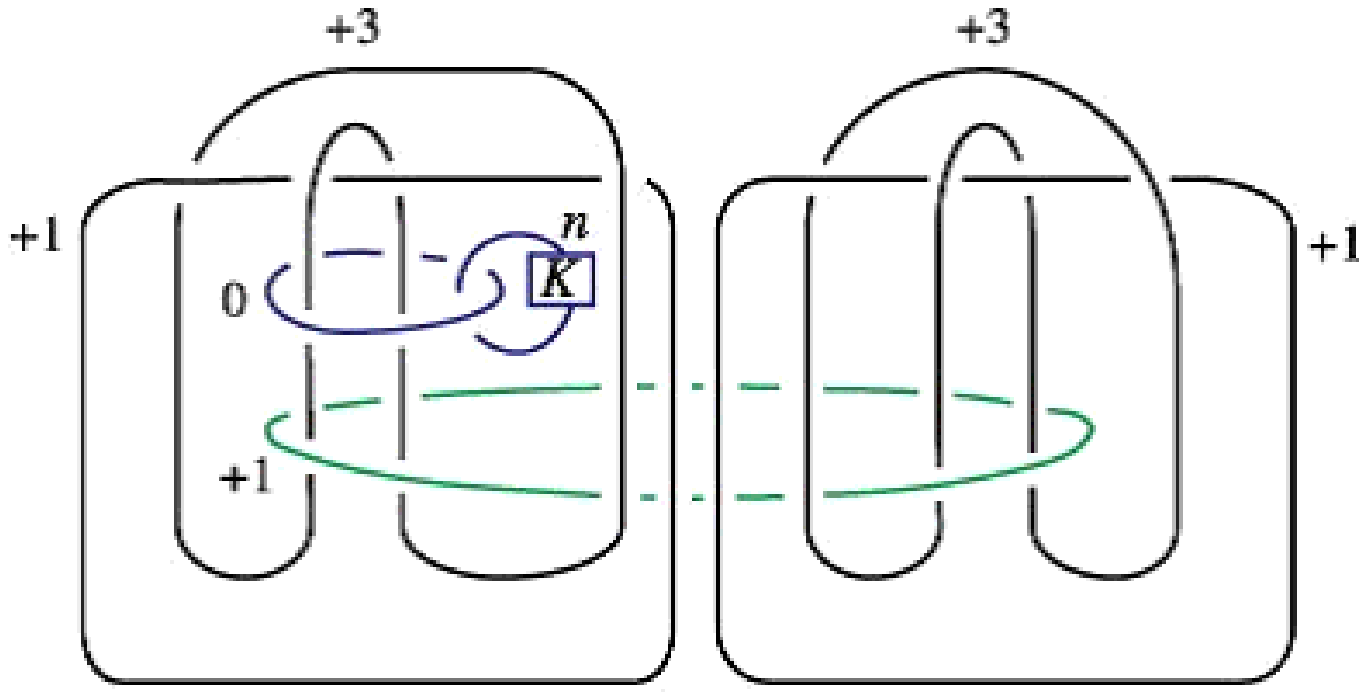}
  \end{center}
   \caption{}
 \end{minipage}%
\begin{minipage}{0.13\hsize}
  \begin{center}
   \includegraphics[height=10mm]{isotopy.eps}
  \end{center}
 \end{minipage}%
 \begin{minipage}{0.35\hsize}
  \begin{center}
   \includegraphics[height=30mm]{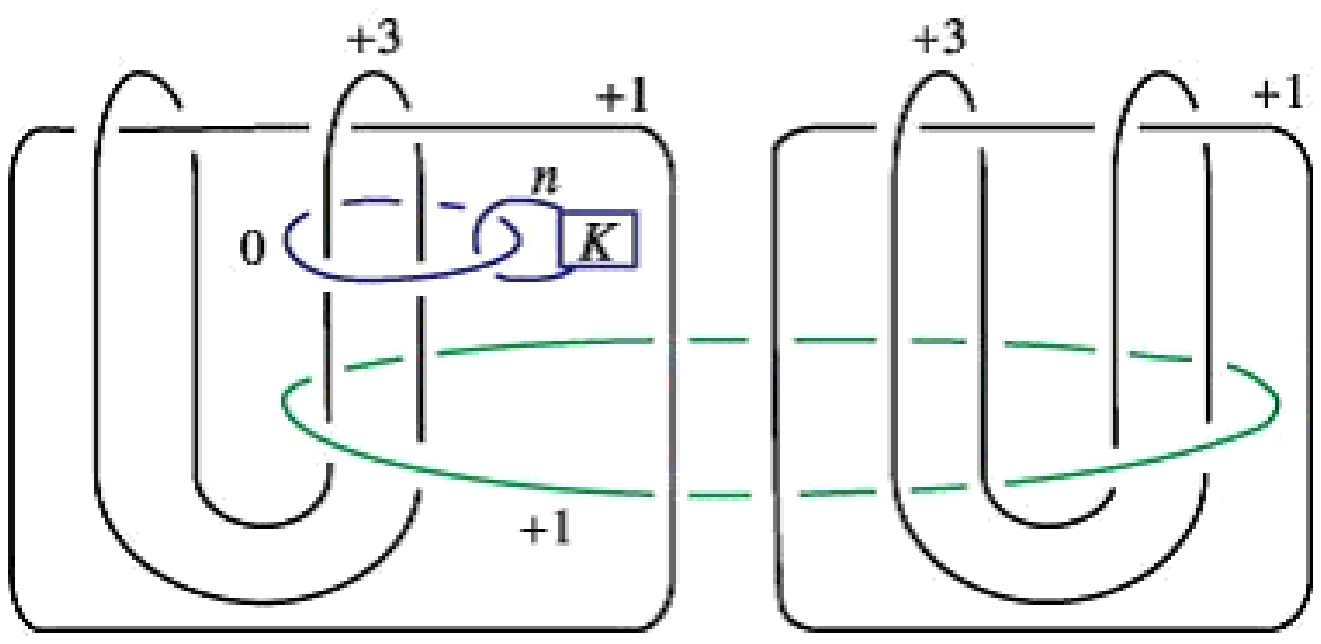}
  \end{center}
  \caption{}
 \end{minipage}
\end{figure}
\begin{figure}[H]
\begin{minipage}{0.1\hsize}
  \begin{center}
   \includegraphics[height=10mm]{isotopy.eps}
  \end{center}
 \end{minipage}%
 \begin{minipage}{0.35\hsize}
  \begin{center}
   \includegraphics[height=34mm]{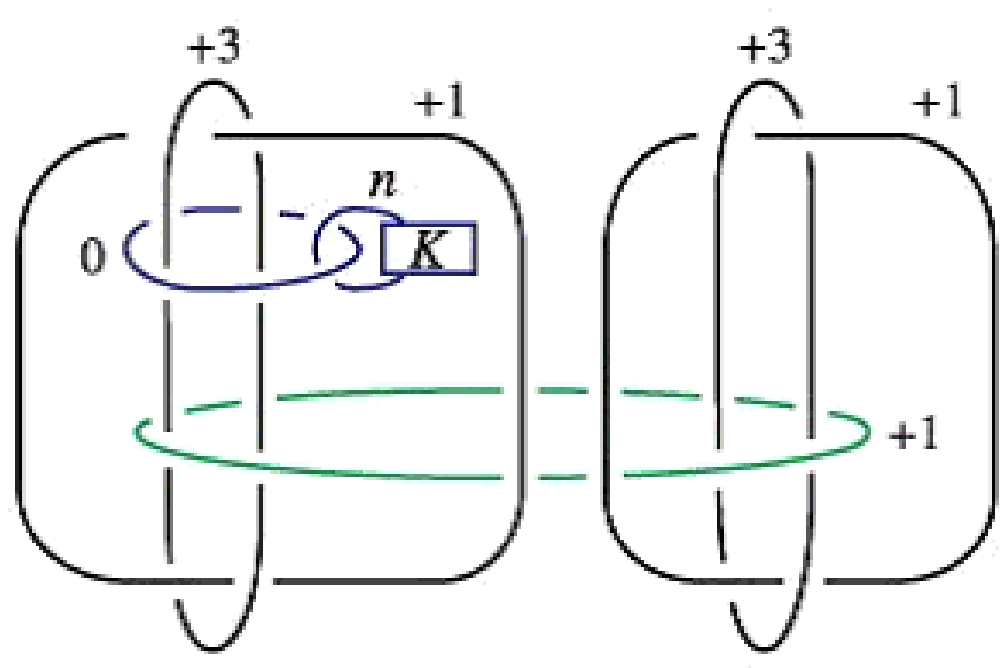}
  \end{center}
   \caption{}
 \end{minipage}%
\begin{minipage}{0.13\hsize}
  \begin{center}
   \includegraphics[height=13mm]{blowup.eps}
  \end{center}
 \end{minipage}%
 \begin{minipage}{0.34\hsize}
  \begin{center}
   \includegraphics[height=34mm]{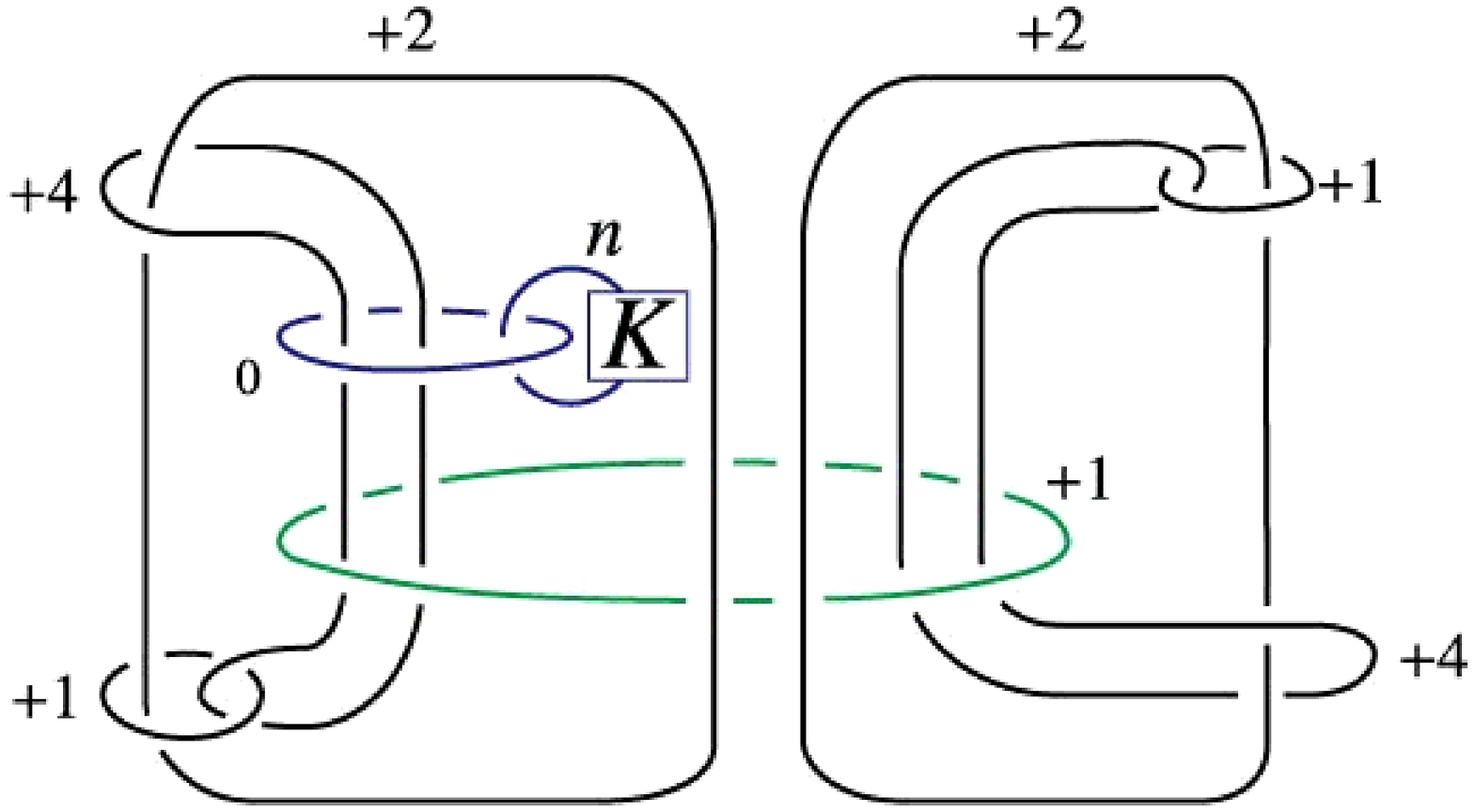}
  \end{center}
  \caption{}
 \end{minipage}
\end{figure}
\begin{figure}[H]
\begin{minipage}{0.1\hsize}
  \begin{center}
   \includegraphics[height=10mm]{isotopy.eps}
  \end{center}
 \end{minipage}%
 \begin{minipage}{0.38\hsize}
  \begin{center}
   \includegraphics[height=31mm]{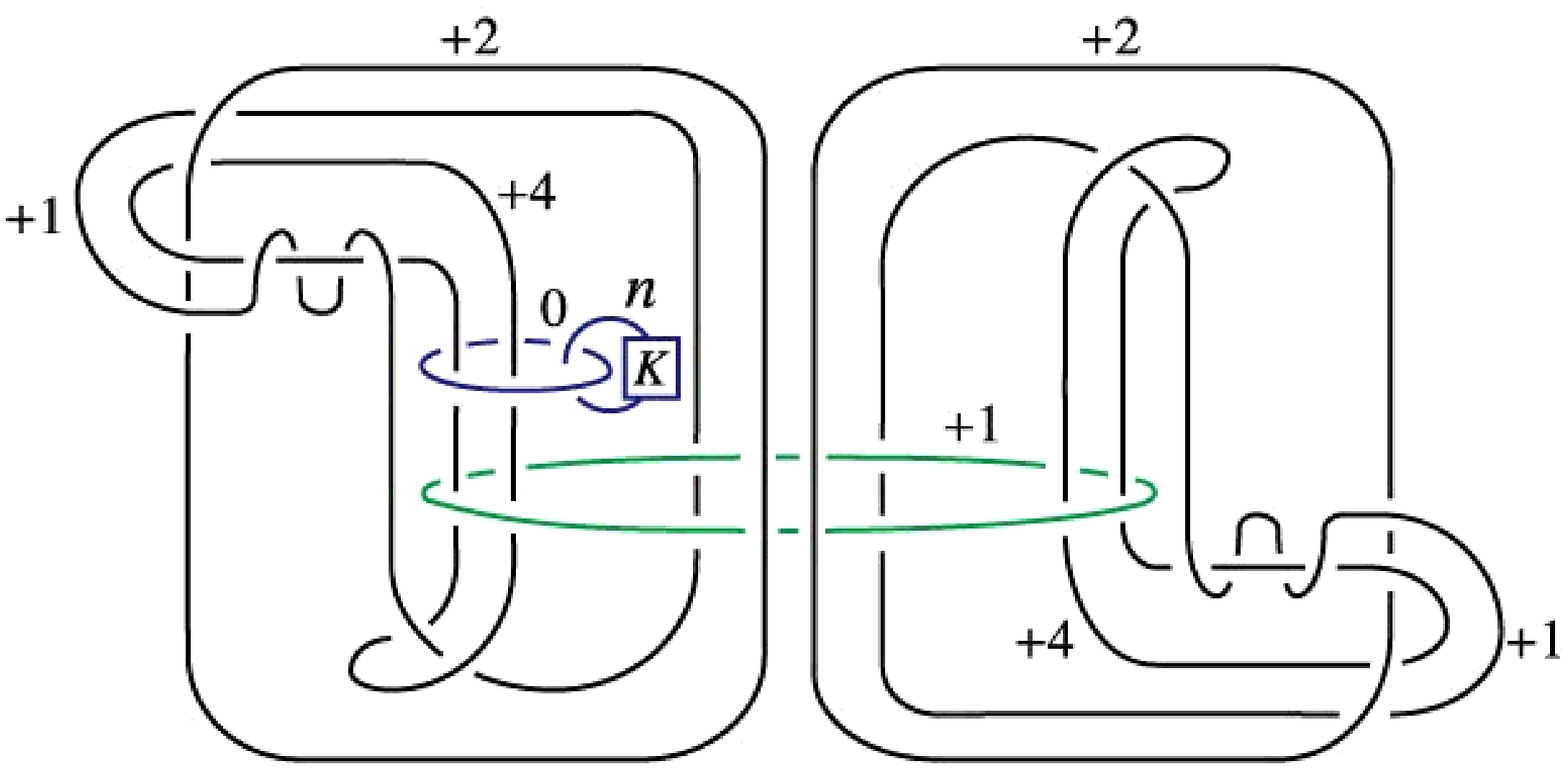}
  \end{center}
   \caption{}
 \end{minipage}%
\begin{minipage}{0.15\hsize}
  \begin{center}
   \includegraphics[height=13mm]{blowup.eps}
  \end{center}
 \end{minipage}%
 \begin{minipage}{0.35\hsize}
  \begin{center}
   \includegraphics[height=41mm]{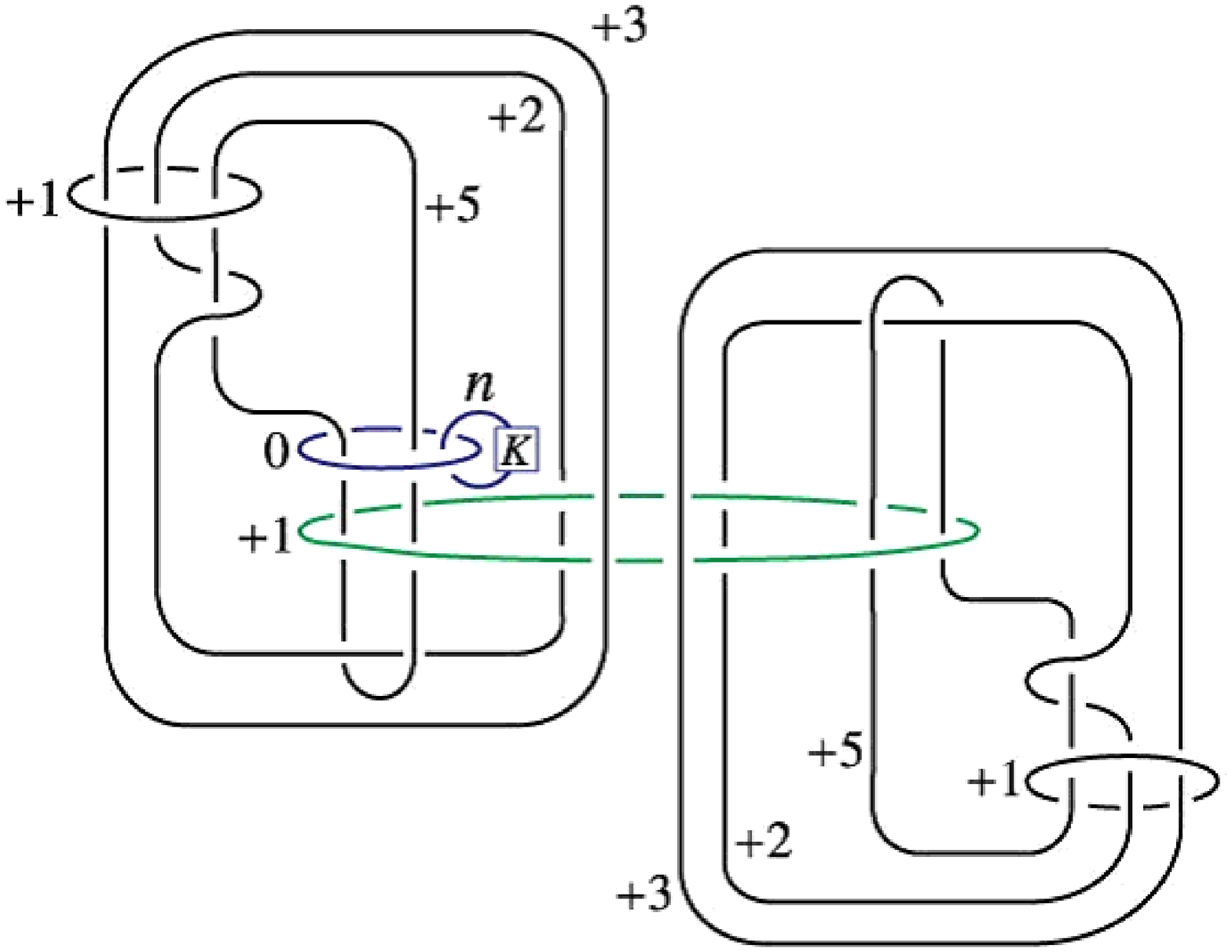}
  \end{center}
  \caption{}
 \end{minipage}
\end{figure}
\begin{figure}[H]
\begin{minipage}{0.08\hsize}
  \begin{center}
   \includegraphics[height=10mm]{isotopy.eps}
  \end{center}
 \end{minipage}%
 \begin{minipage}{0.38\hsize}
  \begin{center}
   \includegraphics[height=47mm]{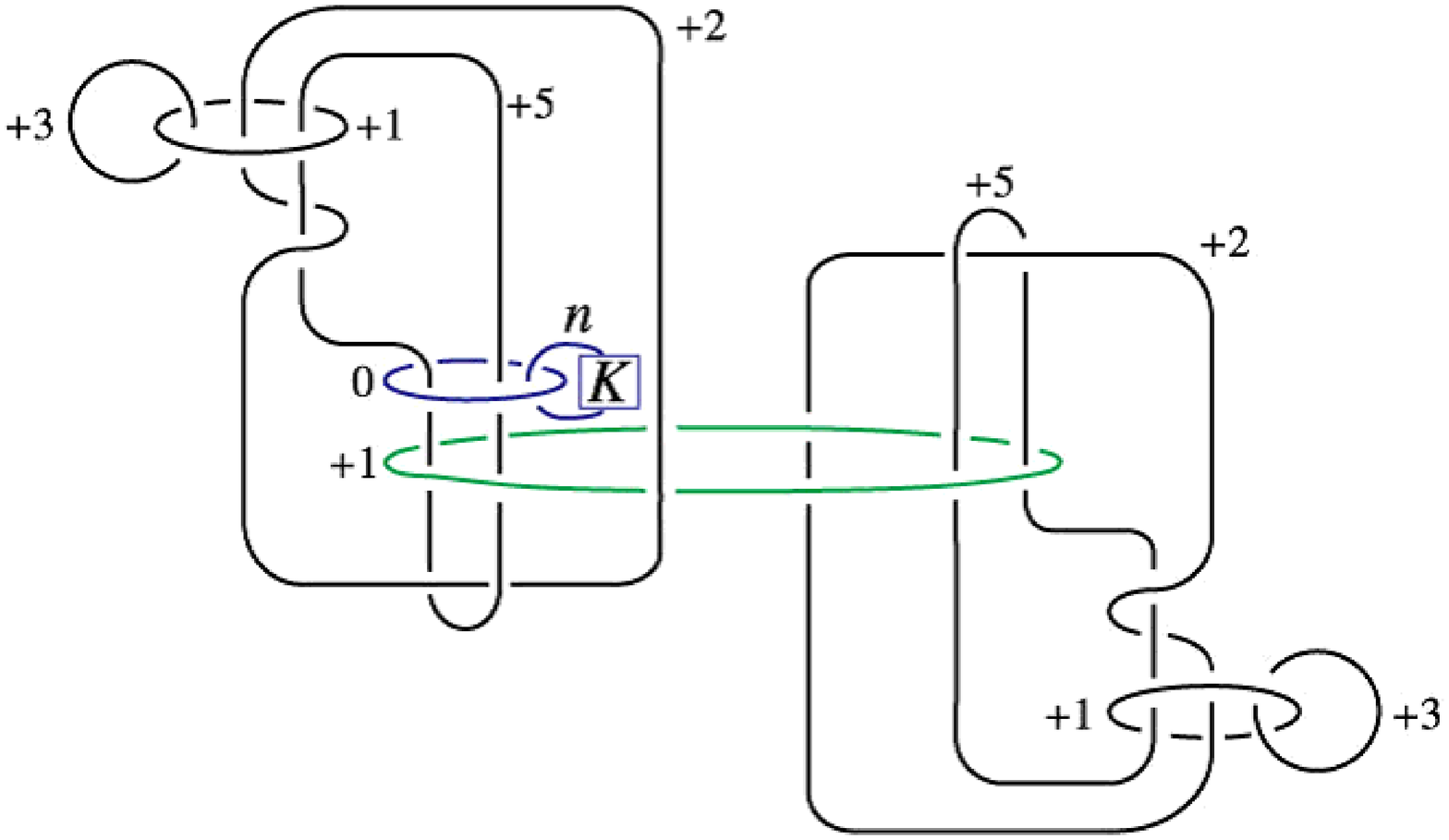}
  \end{center}
   \caption{}
 \end{minipage}%
\begin{minipage}{0.12\hsize}
  \begin{center}
   \includegraphics[height=10mm]{isotopy.eps}
  \end{center}
 \end{minipage}%
 \begin{minipage}{0.39\hsize}
  \begin{center}
   \includegraphics[height=43mm]{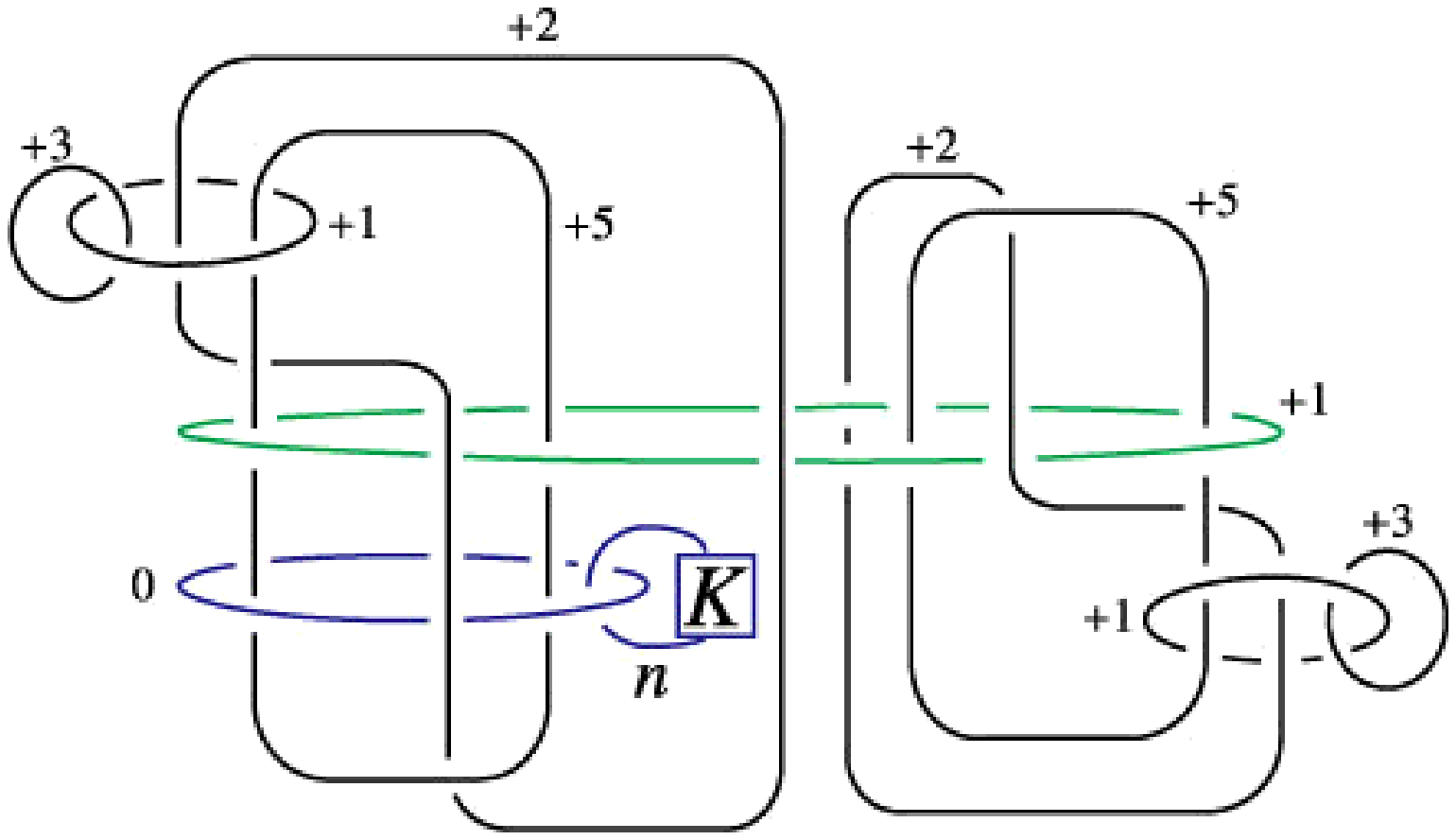}
  \end{center}
  \caption{}
 \end{minipage}
\end{figure}
\begin{figure}[H]
\begin{minipage}{0.1\hsize}
  \begin{center}
   \includegraphics[height=10mm]{isotopy.eps}
  \end{center}
 \end{minipage}%
 \begin{minipage}{0.36\hsize}
  \begin{center}
   \includegraphics[height=46mm]{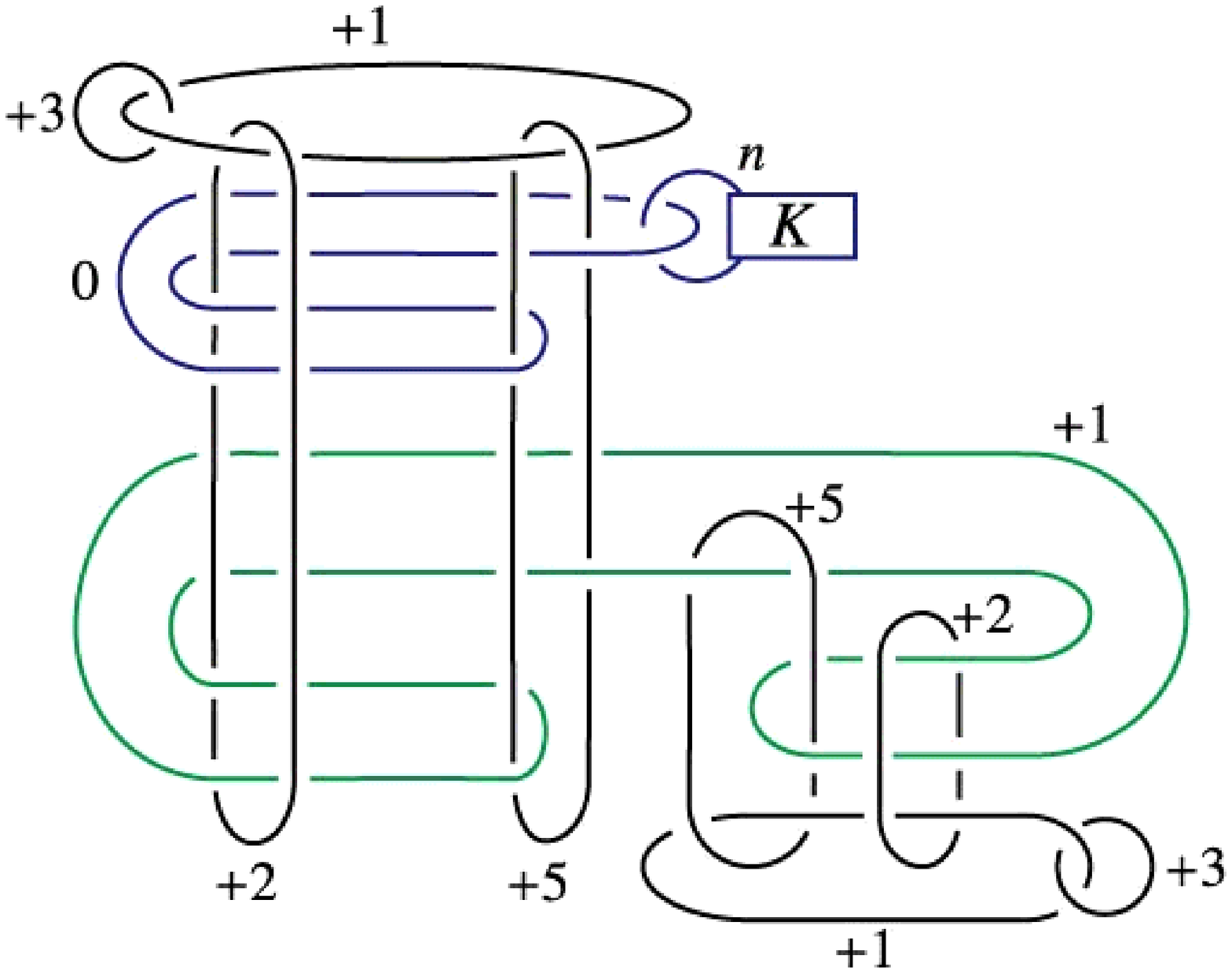}
  \end{center}
   \caption{}
 \end{minipage}%
\begin{minipage}{0.15\hsize}
  \begin{center}
   \includegraphics[height=18mm]{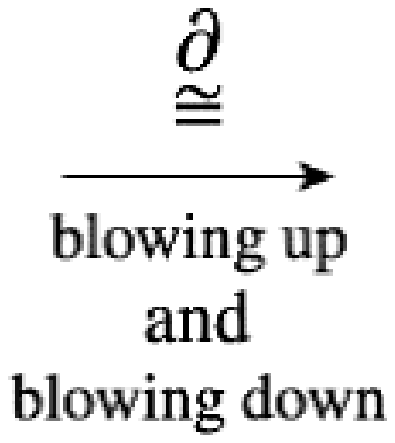}
  \end{center}
 \end{minipage}%
 \begin{minipage}{0.4\hsize}
  \begin{center}
   \includegraphics[height=65mm]{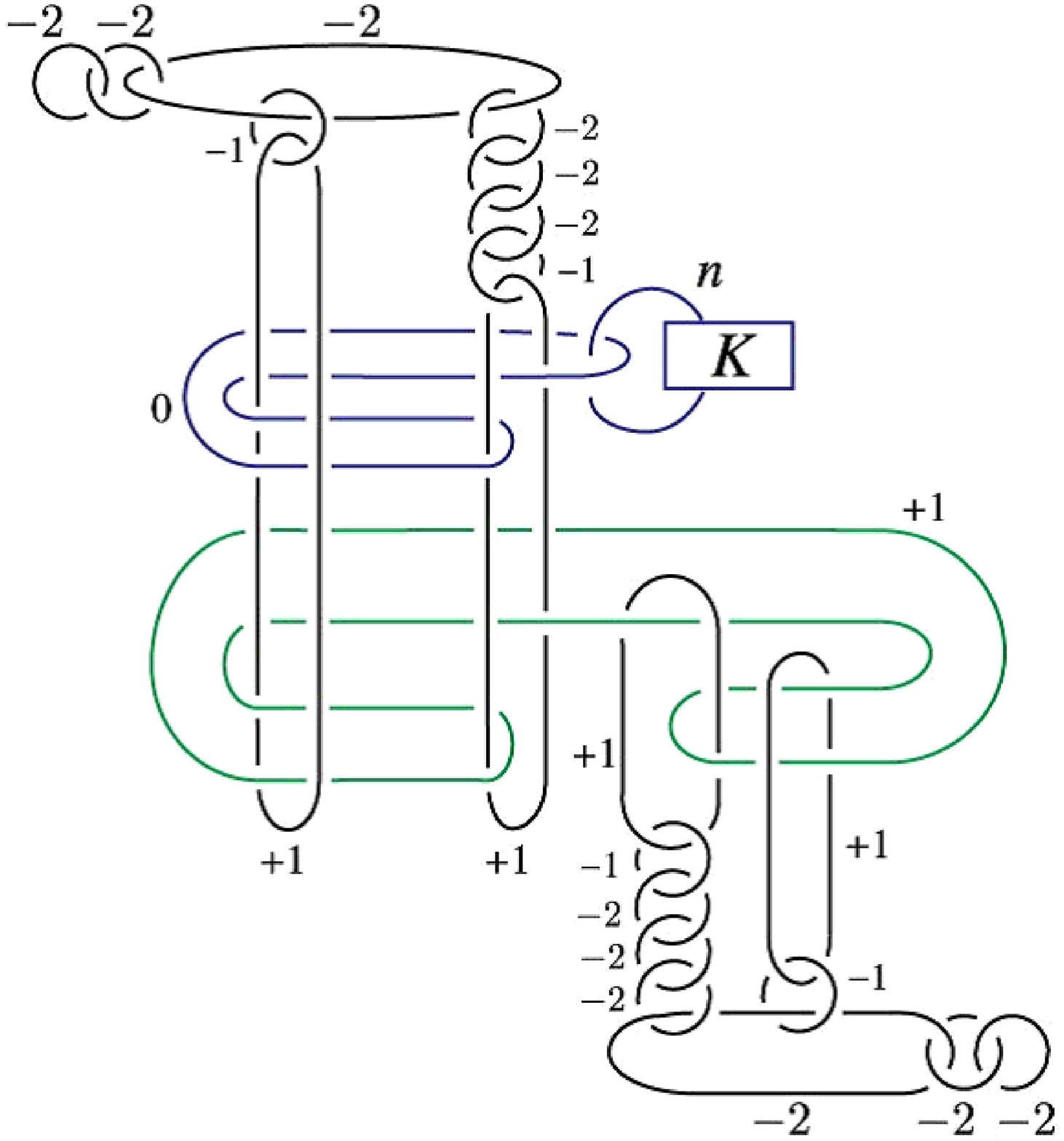}
  \end{center}
  \caption{}
 \end{minipage}
\end{figure}
\begin{figure}[H]
\begin{minipage}{0.07\hsize}
  \begin{center}
   \includegraphics[height=14mm]{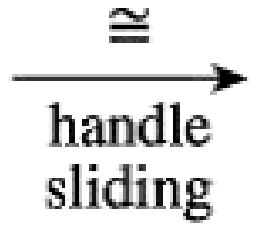}
  \end{center}
 \end{minipage}%
 \begin{minipage}{0.39\hsize}
  \begin{center}
   \includegraphics[height=57mm]{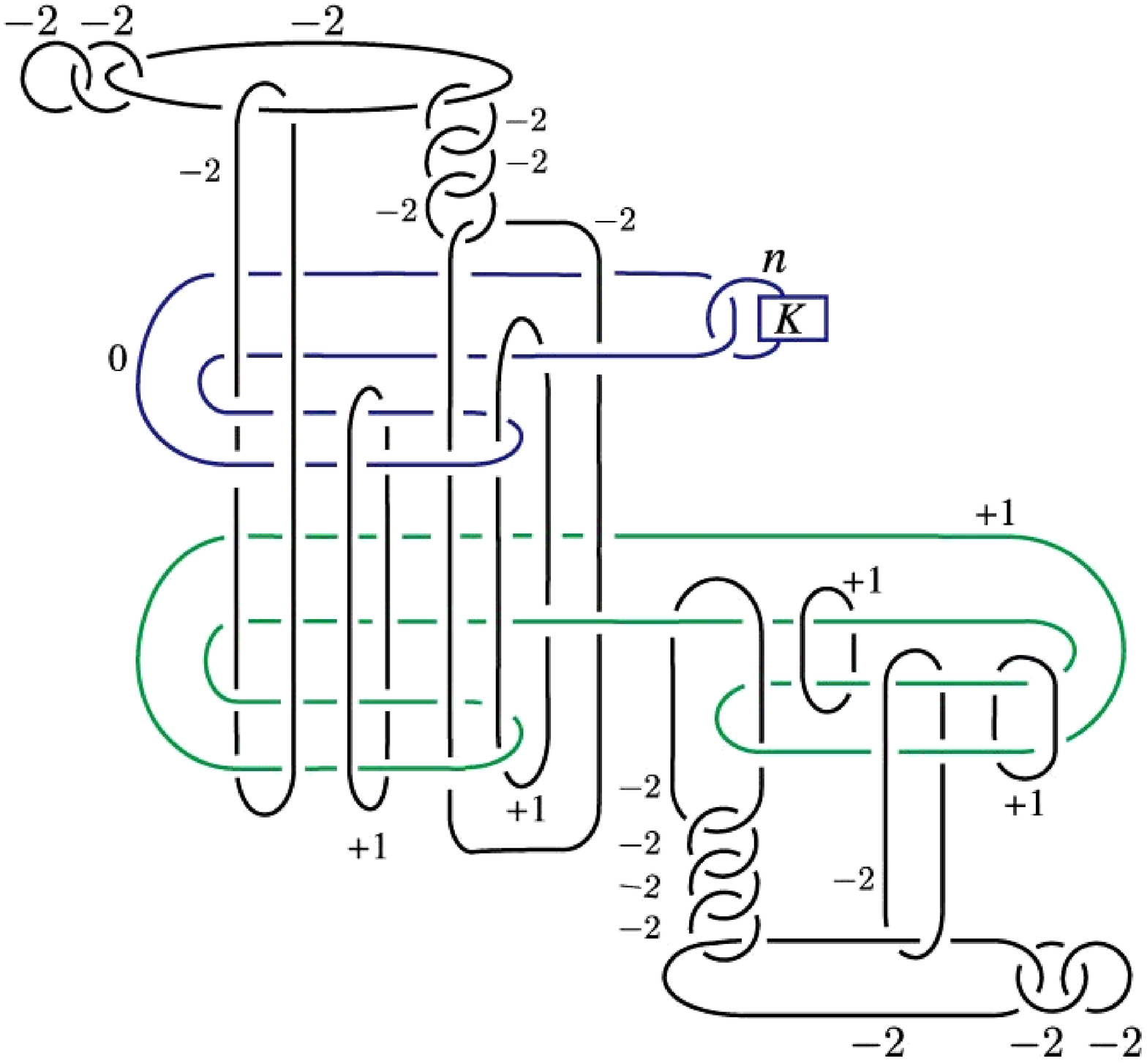}
  \end{center}
   \caption{}
 \end{minipage}%
\begin{minipage}{0.1\hsize}
  \begin{center}
   \includegraphics[height=13mm]{blowdown.eps}
  \end{center}
 \end{minipage}%
 \begin{minipage}{0.4\hsize}
  \begin{center}
   \includegraphics[height=57mm]{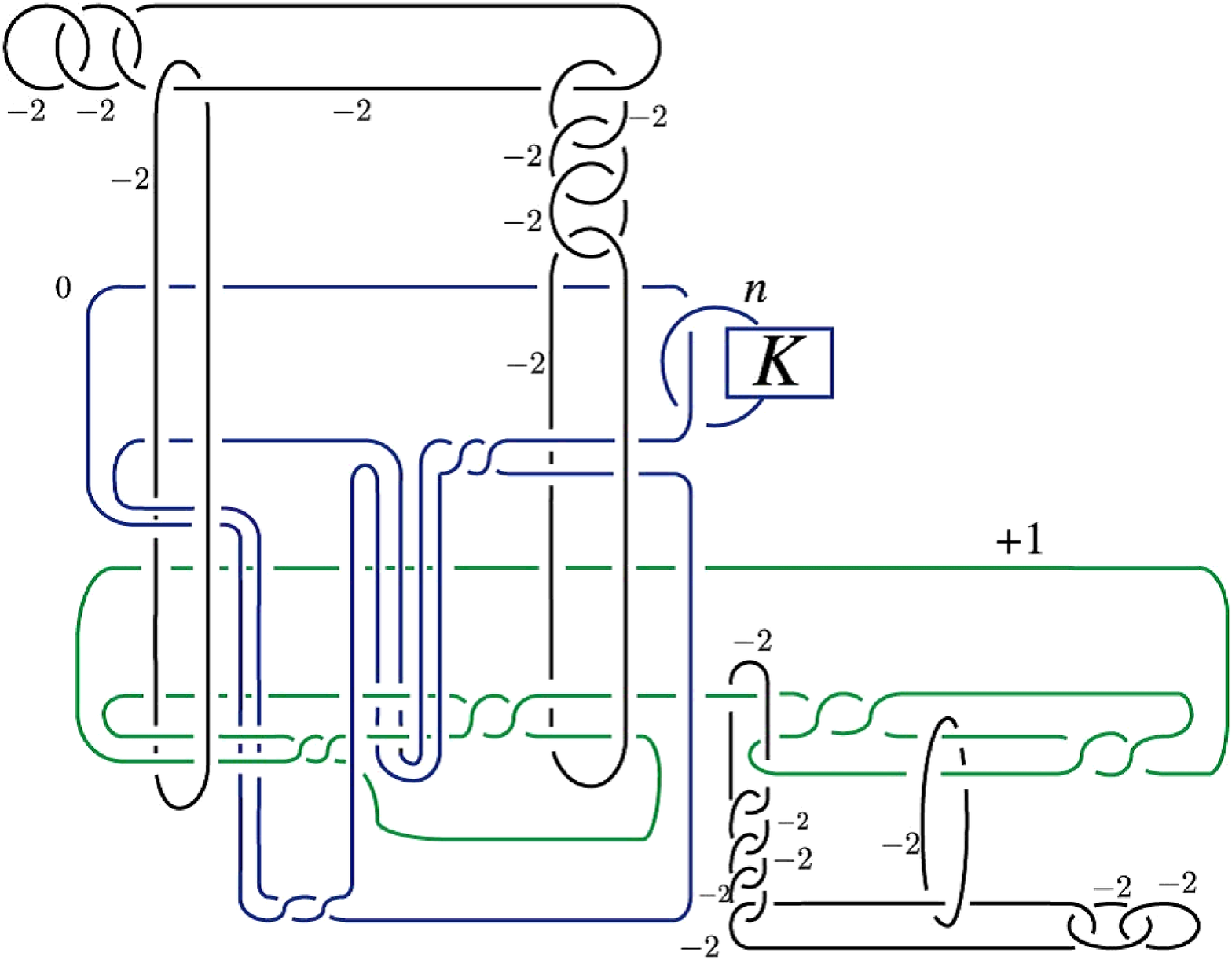}
  \end{center}
  \caption{}
 \end{minipage}
\end{figure}
\begin{figure}[H]
\begin{minipage}{0.3\hsize}
  \centering
   \includegraphics[height=20mm]{surgery3.eps}
 \end{minipage}%
 \begin{minipage}{0.7\hsize}
  \centering
   \includegraphics[height=75mm]{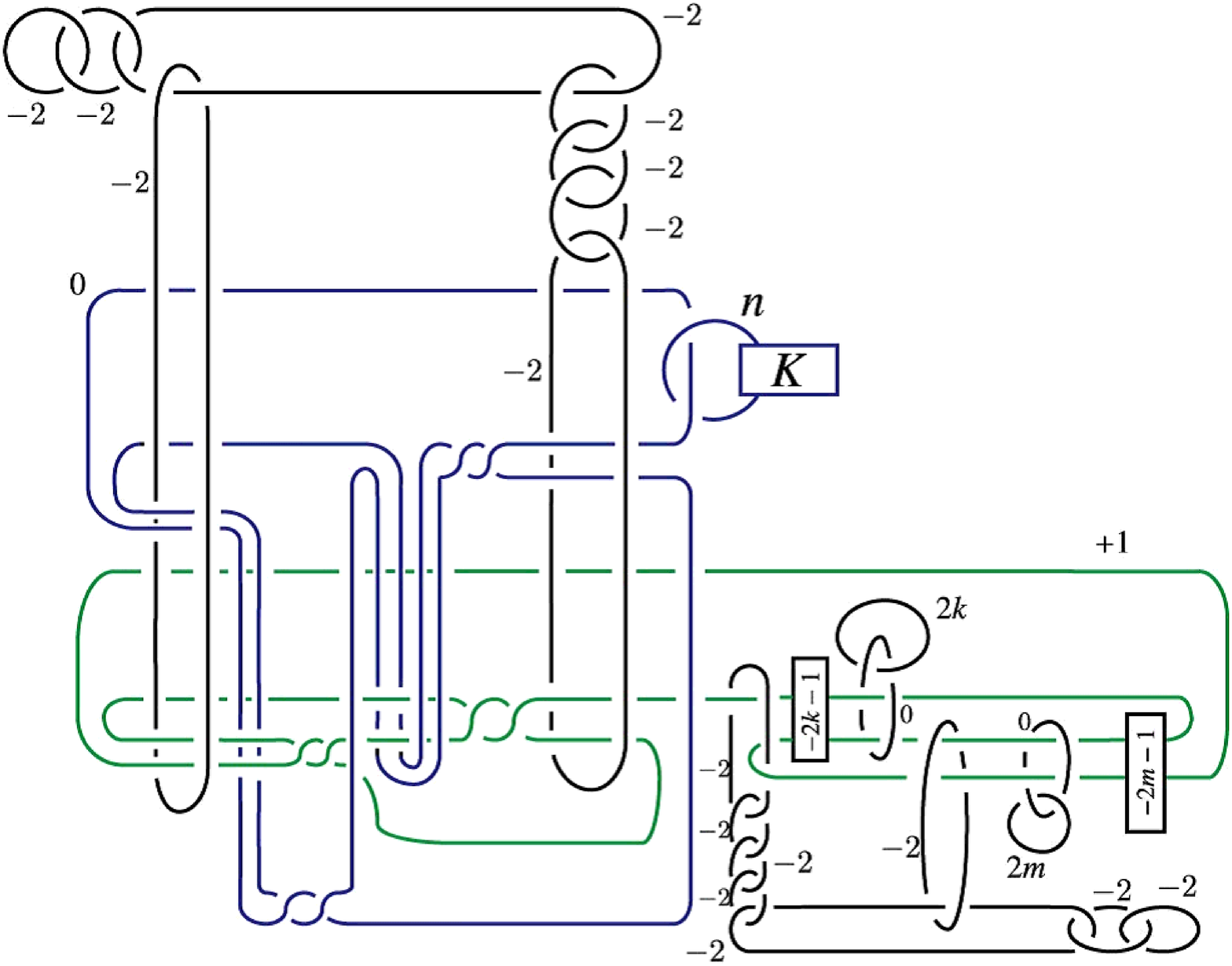}
   \caption{}
  \label{prop2-23}
 \end{minipage}
\end{figure}
Putting $k=m=-1$ in \figref{prop2-23}, we get \figref{YnK}.
\endproof

\begin{Remark}\label{rem4}
If $M_n (K)$ bounds a smooth contractible $4$-manifold $W_n (K)$ and $n$ is even, we have a smooth homotopy $K3$ surface $Z_n(K) \cup_{\partial} (-W_{n}(K))$.
\end{Remark}

Morgan and Szab\'o proved the following adjunction inequality.
\begin{Theorem}[see {\cite[Corollary~1.2]{MS}}]\label{thmMS}
Let $X$ be a smooth closed homotopy $K3$ surface and $g(x)$ be the minimal genus of smoothly embedded closed oriented surfaces representing $x$, where $x\in H_2(X; \Z)$. For every $x\in H_2(X; \Z), x \neq 0, x \cdot x \geq 0$, we have 
\begin{equation*}
2g(x)-2 \geq x \cdot x \, .
\end{equation*}
\end{Theorem}

The author is informed from Mikio Furuta \cite{F} that the following stronger version can be proved essentially in the same way as \thmref{thmMS}.

\begin{Furuta}[\cite{F}]
Let $X$ be a smooth closed spin rational homology $K3$ surface and $g(x)$ be the minimal genus of smoothly embedded closed oriented surfaces representing $x$, where $x\in H_2(X; \Z)$. For every $x\in H_2(X; \Z), x \neq 0, x \cdot x \geq 0$, we have 
\begin{equation*}
2g(x)-2 \geq x \cdot x \,.
\end{equation*}
\end{Furuta}

We will prove \thmref{thmt}.
\proof[Proof of \thmref{thmt}]
Let $N_{K, n}$ be a $4$-dimensional handlebody constructed by attaching a $2$-handle to $D^4$ along the $n$-framed knot $K$ of \figref{ZnK}.
Let $x$ be a generator of $H_2(N_{K, n}; \Z)$.
This homology class $x$ is represented by smooth closed oriented surface $\Sigma$ with genus $g^n_s(K)$ and $x \cdot x =n$.
The handlebody $N_{K, n}$ is a submanifold of $Z_n(K)$.

Suppose that $n>2g^n_s(K)-2$, $M_n(K)$ $(=\partial (Z_n(K)))$ bounds a smooth spin rational $4$-ball $W_n(K)$ and $n$ is even, then we have a smooth spin rational homology $K3$ surface $Z_n(K) \cup_{\partial} (-W_{n}(K))$ in which $x$ persists and is represented by $\Sigma$ with genus $g^n_s(K)$ and $x \cdot x =n$.
We apply the adjunction inequality to $x$ in $Z_n(K) \cup_{\partial} (-W_{n}(K))$. Then we get
\begin{equation*}
2g^n_s(K) -2 \geq x \cdot x = n.
\end{equation*}
This contradicts the assumption $n>2g^n_s(K)-2$.
Therefore if $n>2g^n_s(K)-2$ and $n$ is even, $M_n(K)$ does not bound any smooth spin rational $4$-ball.
By \remref{rem2}, we know that if $n$ is odd, $M_n(K)$ does not bound any smooth spin rational $4$-ball.
Therefore we conclude that if $n>2g^n_s(K)-2$, $M_n(K)$ does not bound any smooth spin rational $4$-ball.
\endproof

We will show \corref{cort}.
\proof[Proof of \corref{cort}]
By \remref{rem3}, if $D_-(LHT \sharp K, n)$ is a slice knot, $M_n(K)$ bounds a contractible $4$-manifold. By \thmref{thmt}, if $n>2g^n_s(K)-2$, $M_n(K)$ does not bound any smooth spin rational $4$-ball. Therefore if $n>2g^n_s(K)-2$, $D_-(LHT \sharp K, n)$ is not a slice knot.
\endproof
%
%
%
%
%
\section{Proof of \thmref{thm:sec}}\label{sec:3}
We show that if $(a_1, a_2, a_3, a_4)$ is any permutation of $(-1, -1, -2k-1,-2m-1)$, the knot represented by \figref{thmt-2} is not a slice knot, where $m\in\Z$ and $k\in\Z_{\geq0}$.
\proof[Proof of \thmref{thm:sec}]
Let $Y$ be the $4$-dimensional handlebody represented by \figref{thmt2-2} with intersection form $2E_8\oplus3\displaystyle\bigl(\begin{smallmatrix} 0 & 1 \\ 1 & 0 \end{smallmatrix} \bigr)\oplus1$.
By the proof of \propref{prop2}, $X_{-6}(U)$ and $Y$ have the same boundary. 

If $k\geq0$ and the $+1$-framed knot in \figref{thmt2-2} is a slice knot, then we can blow it down and get a smooth simply connected $4$-dimensional handlebody $Z$ represented by \figref{thmt2-3} with intersection form $2E_8\oplus3\displaystyle\bigl(\begin{smallmatrix} 0 & 1 \\ 1 & 0 \end{smallmatrix} \bigr)$. Because $M_{-6}(U)$ bounds a contractible $4$-manifold $W_{-6}(U)$ (see \remref{rem3}), we have a homotopy $K3$ surface $Z \cup_{\partial} (-W_{-6}(U))$. Let $x$ be the element of $H_2(Z \cup_{\partial} (-W_{-6}(U)); \Z)$ which is generated by attaching a $2$-handle along the $2k$-framed unknot in \figref{thmt2-3}. The homology class $x$ is represented by a smooth $S^2$ and $x \cdot x = 2k$. By \thmref{thmMS}, we have the following inequality:
 \begin{equation*}
-2 \geq x \cdot x = 2k.
\end{equation*}
This contradicts the assumption $k\geq0$. Therefore if $k\geq0$, the $+1$-framed knot in \figref{thmt2-2} is not a slice knot.
By the handle calculus from \figref{thmt2-4} to \figref{thmt2-5}, we can show that if $k\geq0$, the $+1$-framed knot in \figref{thmt2-5} is not a slice knot essentially in the same way as above.
Similarly we can prove that if $k\geq0$ and $(a_1, a_2, a_3, a_4)$ is any permutation of $(-1, -1, -2k-1,-2m-1)$, the knot represented by \figref{thmt-2} is not a slice knot.
\endproof
\begin{figure}[H]
 \begin{minipage}{0.4\hsize}
  \centering
   \includegraphics[height=21mm]{W-6U.eps}
  \caption{$W_{-6}(U)$}
 \end{minipage}%
\begin{minipage}{0.2\hsize}
  \centering
   \includegraphics[height=15mm]{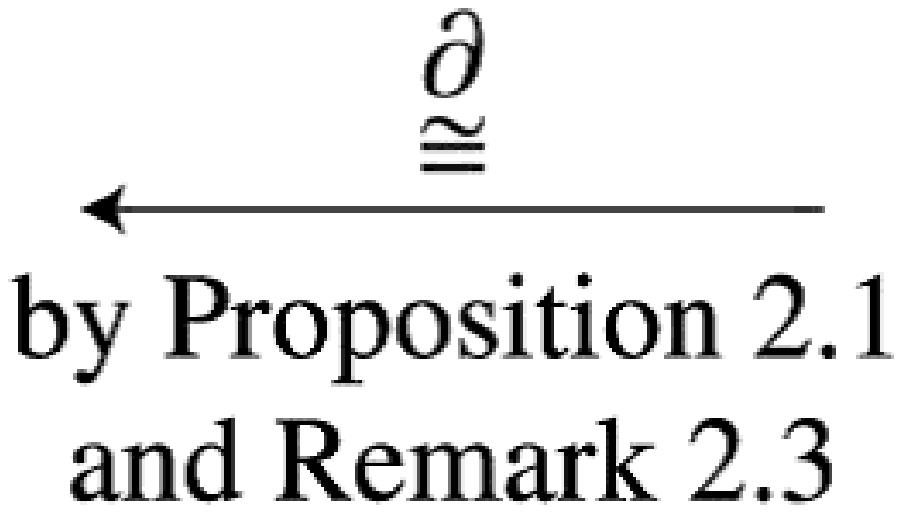}
 \end{minipage}%
 \begin{minipage}{0.4\hsize}
  \centering
   \includegraphics[height=20mm]{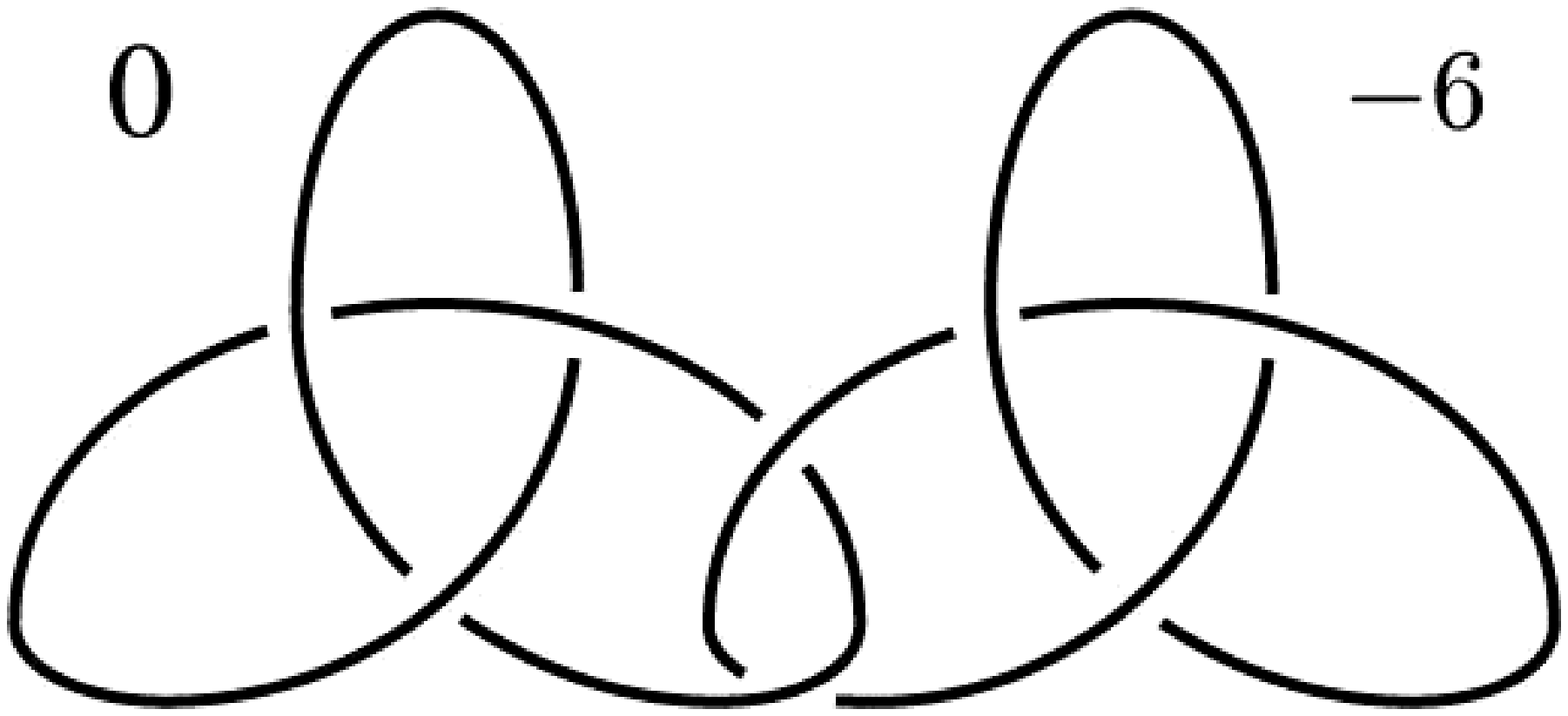}
  \caption{$X_{-6}(U)$}
 \end{minipage}
\end{figure}

\begin{figure}[H]
\begin{minipage}{0.4\hsize}
  \centering
   \includegraphics[height=20mm]{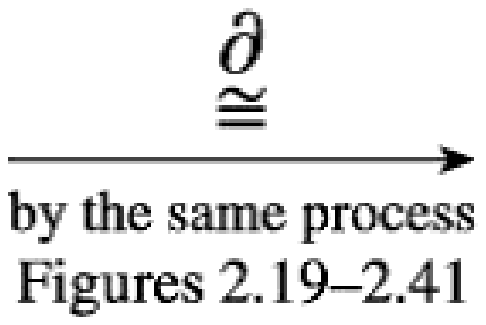}
 \end{minipage}%
 \begin{minipage}{0.6\hsize}
  \centering
   \includegraphics[height=72mm]{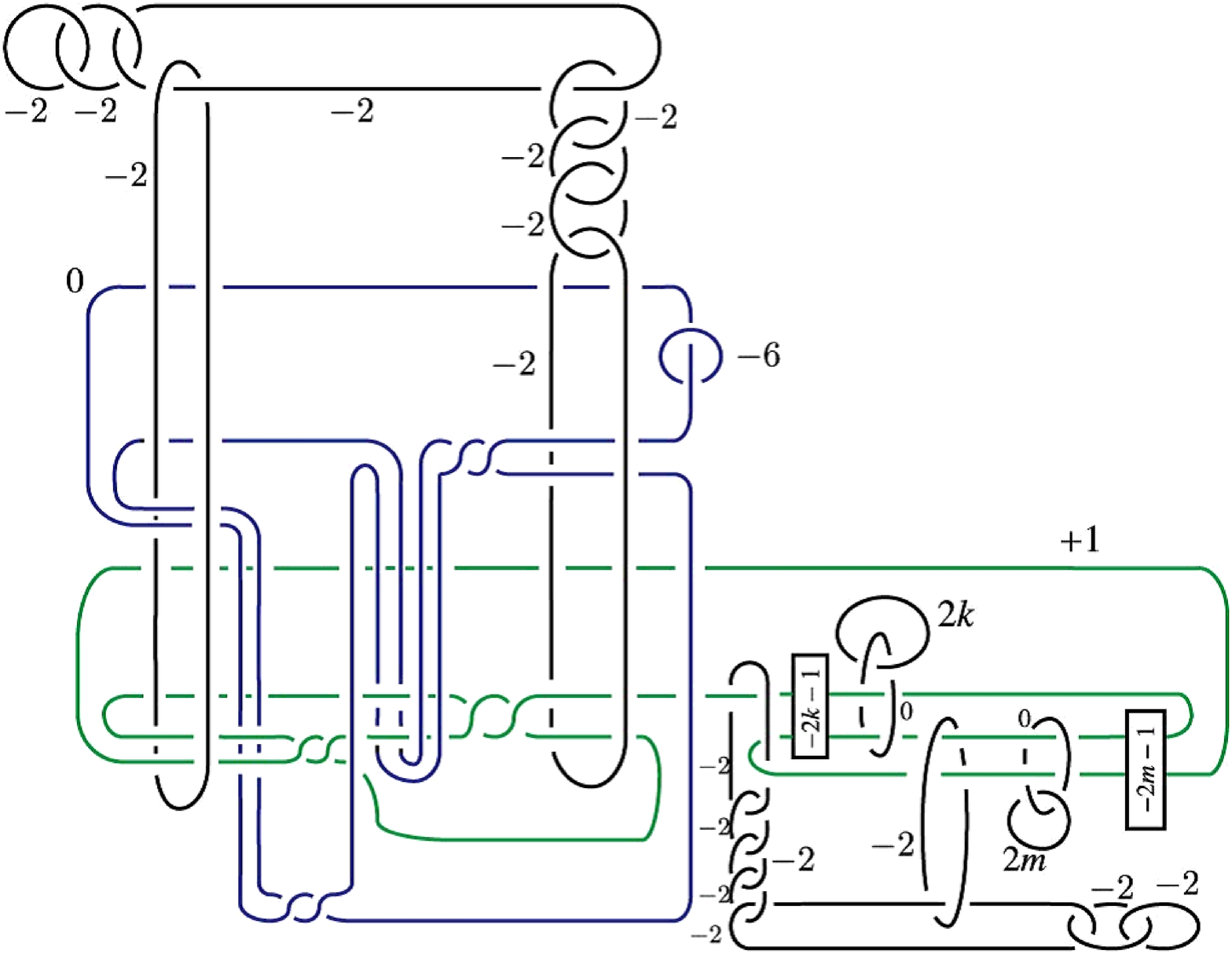}
  \caption{Y}
  \label{thmt2-2}
 \end{minipage}
\end{figure}

\begin{figure}[H]
\begin{minipage}{0.3\hsize}
  \begin{center}
   \includegraphics[height=25mm]{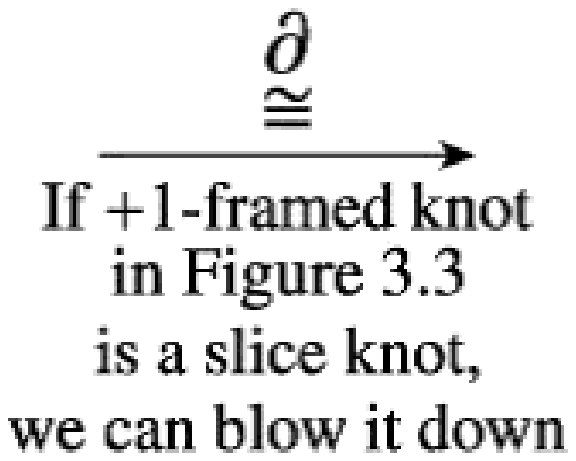}
  \end{center}
 \end{minipage}%
 \begin{minipage}{0.7\hsize}
  \begin{center}
   \includegraphics[height=72mm]{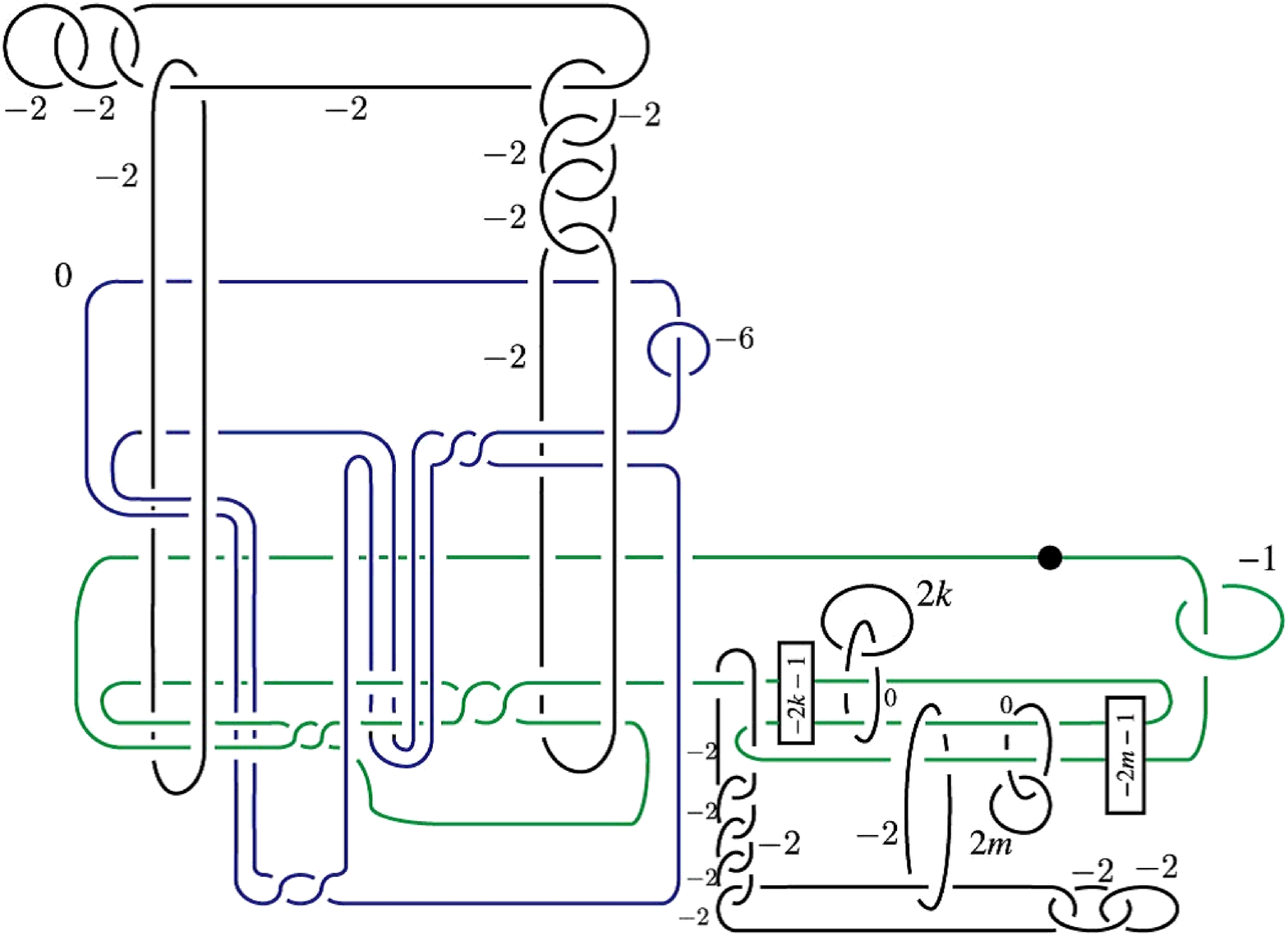}
  \end{center}
  \caption{Z}
  \label{thmt2-3}
 \end{minipage}
\end{figure}

\begin{figure}[H]
 \begin{minipage}{0.37\hsize}
  \begin{center}
   \includegraphics[height=22mm]{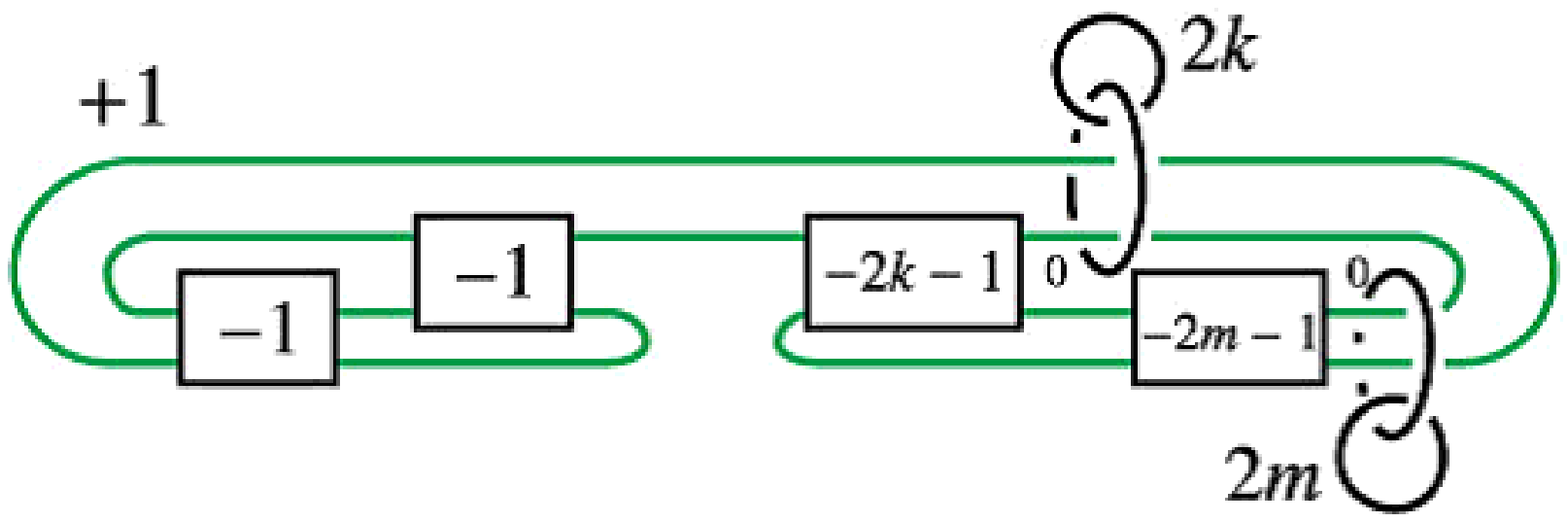}
  \end{center}
  \caption{}
  \label{thmt2-4}
 \end{minipage}%
\begin{minipage}{0.2\hsize}
  \begin{center}
   \includegraphics[height=9mm]{bdydiff.eps}
  \end{center}
 \end{minipage}%
 \begin{minipage}{0.36\hsize}
  \begin{center}
   \includegraphics[height=22mm]{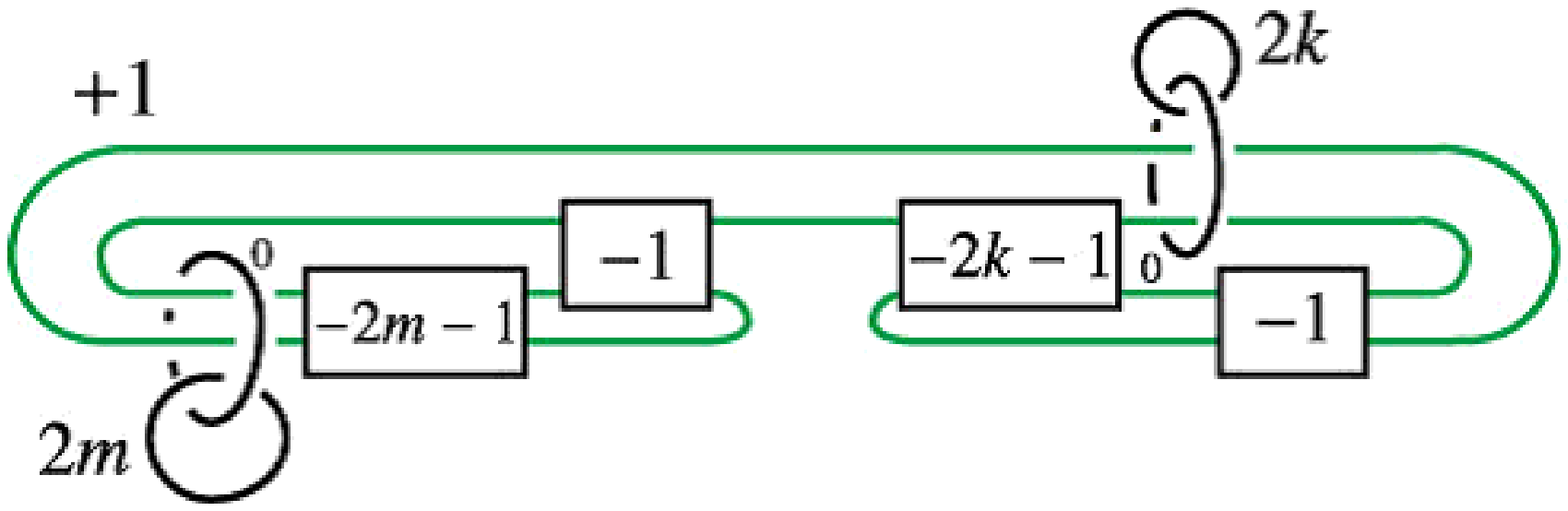}
  \end{center}
  \caption{}
  \label{thmt2-5}
 \end{minipage}
\end{figure}
%
%
%
%
\section{Some homotopy $K3$ surfaces}
We will exhibit some homotopy $K3$ surfaces obtained by the method of this paper.

\begin{figure}[H]
\begin{minipage}{0.45\hsize}
  \centering
   \includegraphics[height=20mm]{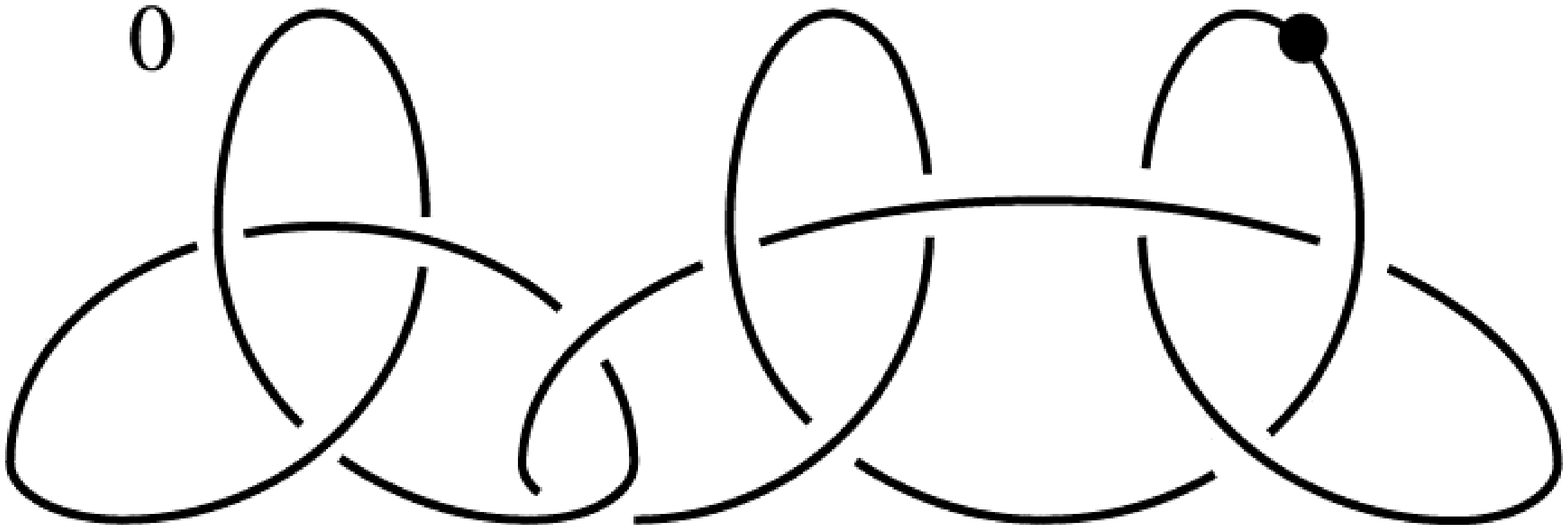}
   \caption{$W'_0 (RHT)$}
   \label{W'0RHT}
 \end{minipage}%
\begin{minipage}{0.15\hsize}
  \centering
   \includegraphics[height=10mm]{bdydiff.eps}
 \end{minipage}%
 \begin{minipage}{0.4\hsize}
  \centering
   \includegraphics[height=20mm]{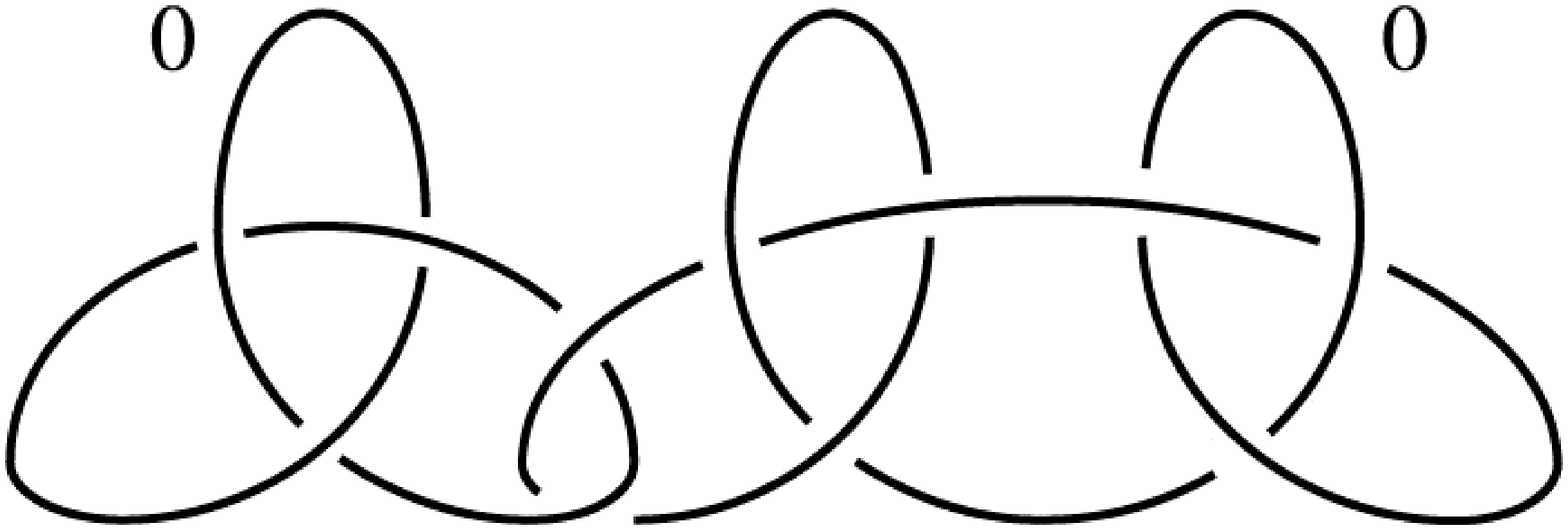}
  \caption{$X_0(RHT)$}
  \label{X0RHT}
 \end{minipage}
\end{figure}

(a). Since $RHT \sharp LHT$ is a slice knot, we have a smooth $S^2$ with self intersection $0$ in $X_0(RHT)$. By performing surgery on the $S^2$, we have a contractible $4$-manifold $W'_0(RHT)$ represented by \figref{W'0RHT}. By \propref{prop2}, $W'_0(RHT)$ and $Z_0(RHT)$ have the same boundary. Therefore we have a homotopy $K3$ surface $Z_0(RHT) \cup_{\partial} (-W'_0(RHT))$.

\begin{figure}[H]
\begin{minipage}{0.35\hsize}
  \centering
   \includegraphics[height=18mm]{sproc.eps}
 \end{minipage}%
 \begin{minipage}{0.65\hsize}
  \centering
   \includegraphics[height=71mm]{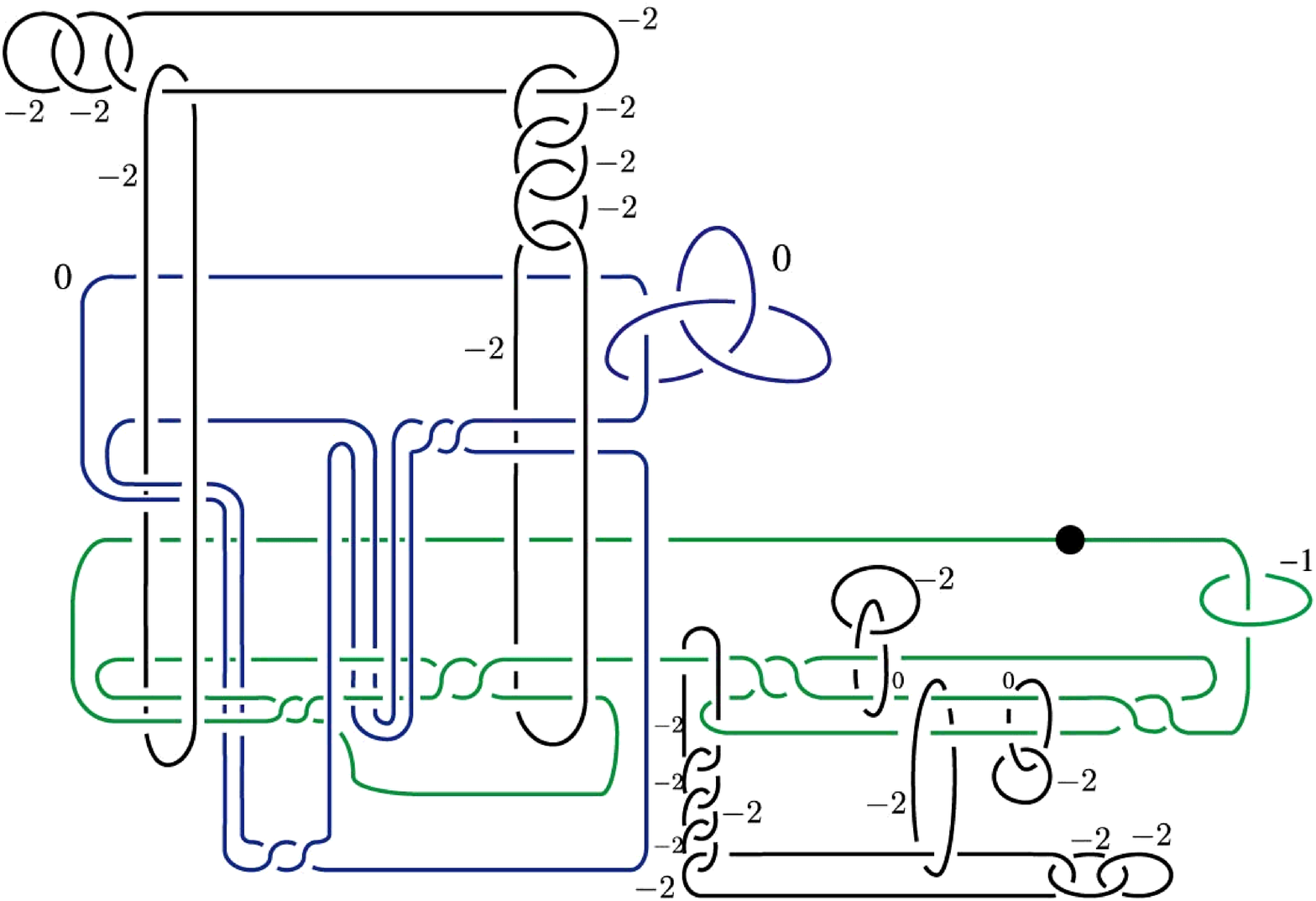}
  \caption{$Z_0(RHT)$}
  \label{Z0RHT}
 \end{minipage}
\end{figure}
\begin{figure}[H]
\begin{minipage}{0.51\hsize}
  \centering
   \includegraphics[height=30mm]{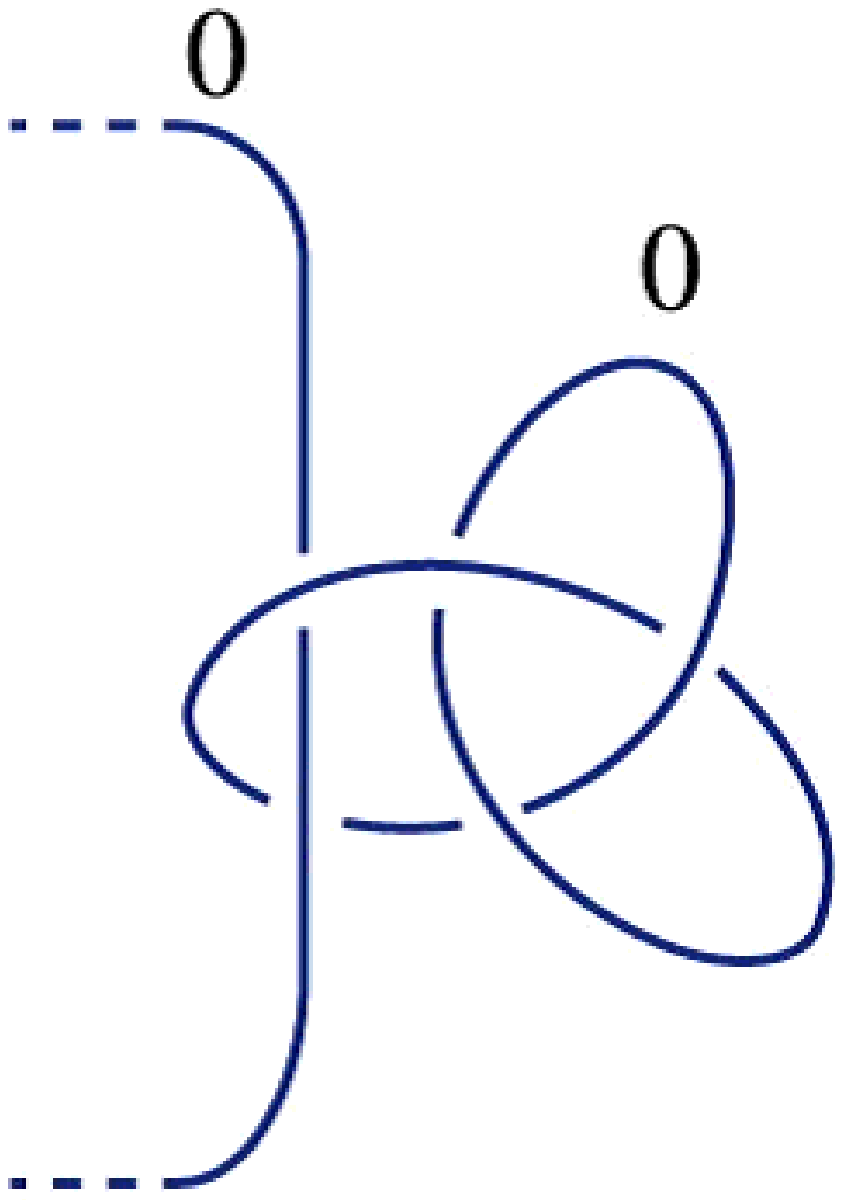}
   \caption{$0$-framed $RHT$ in $Z_0(RHT)$}
   \label{RHTinZ0RHT}
 \end{minipage}%
\begin{minipage}{0.07\hsize}
  \centering
   \includegraphics[height=10mm]{diff.eps}
 \end{minipage}%
 \begin{minipage}{0.42\hsize}
  \centering
   \includegraphics[height=30mm]{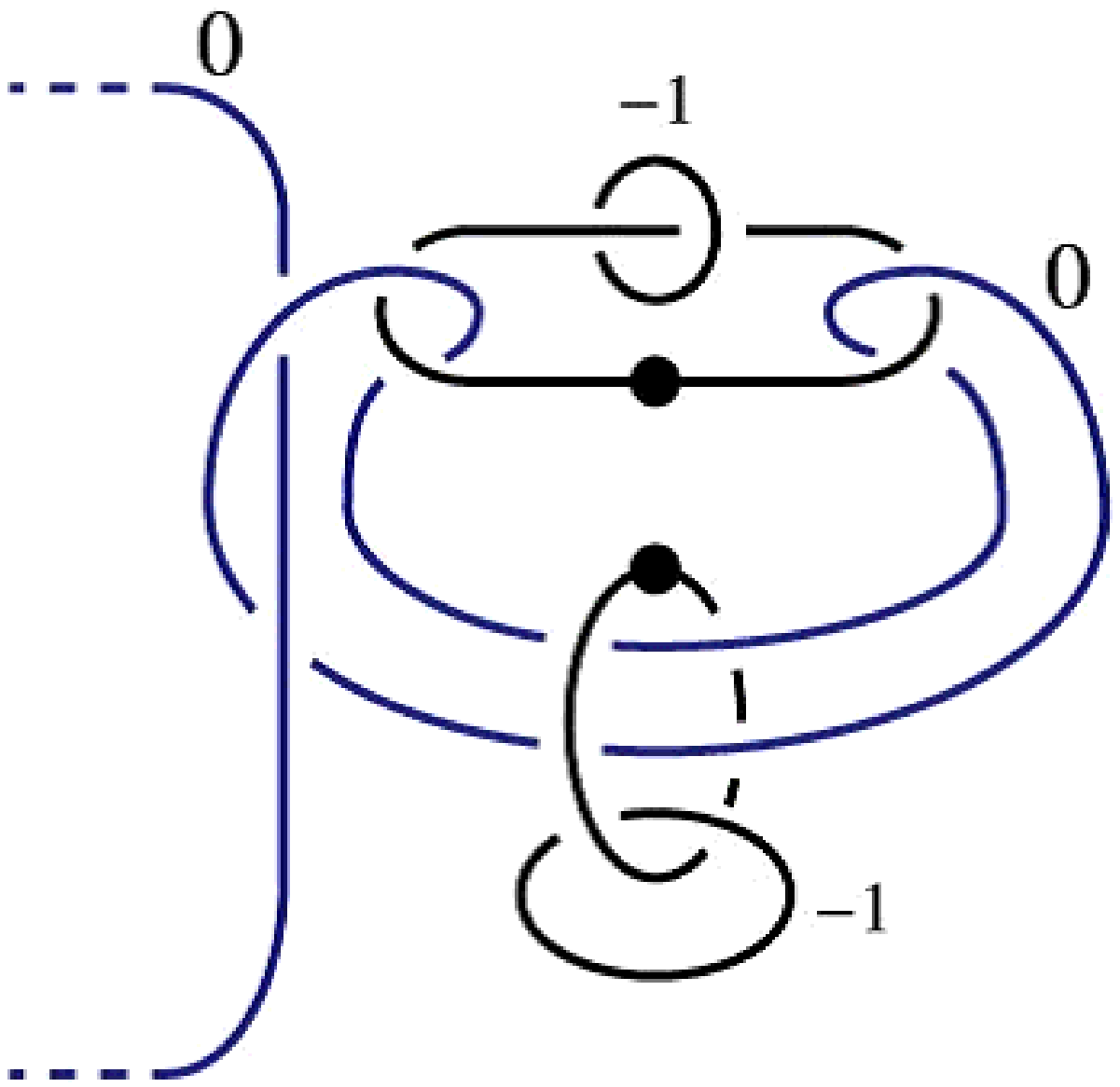}
  \caption{$T^2 \times D^2$ in $Z_0(RHT)$}
  \label{RHTinZ0RHT2}
 \end{minipage}
\end{figure}

Let $X$ be the homotopy $K3$ surface $Z_0(RHT) \cup_{\partial} (-W'_0(RHT))$. By  $g^0_s(RHT)=g_4(RHT)=1$, we have a smooth $T^2$ with trivial normal bundle represented by \figref{RHTinZ0RHT2} in $Z_0(RHT) ( \subset X)$. Then we can perform {\it knot surgery} which is introduced by Fintushel and Stern \cite{FS} on the $T^2 \times D^2$ in $X$. By the existence of the $2$-handle with $0$-framing which is linking to the $0$-framed $RHT$ in $Z_0(RHT)$, we have $\pi_1(X\setminus T^2)=1$. If we perform the knot surgery using a knot whose symmetric Alexander polynomial is $t^2-t+1$ for example, then the resulting $4$-manifold is an exotic homotopy $K3$ surface (see \cite{FS}).

(b). By \propref{prop1}, $X_0(RHT)$ and $V_0(RHT)$ have the same boundary. Since $RHT \sharp LHT$ is a slice knot, $D_-(RHT\sharp LHT, 0)$ is a slice knot. By \remref{rem3}, we have a contractible $4$-manifold represented by \figref{W0RHT}.
\begin{figure}[H]
\begin{minipage}{0.4\hsize}
  \centering
   \includegraphics[height=15mm]{X0RHT.eps}
   \caption{$X_0(RHT)$}
 \end{minipage}%
\begin{minipage}{0.2\hsize}
  \centering
   \includegraphics[height=8mm]{bdydiff.eps}
 \end{minipage}%
 \begin{minipage}{0.4\hsize}
  \centering
   \includegraphics[height=18mm]{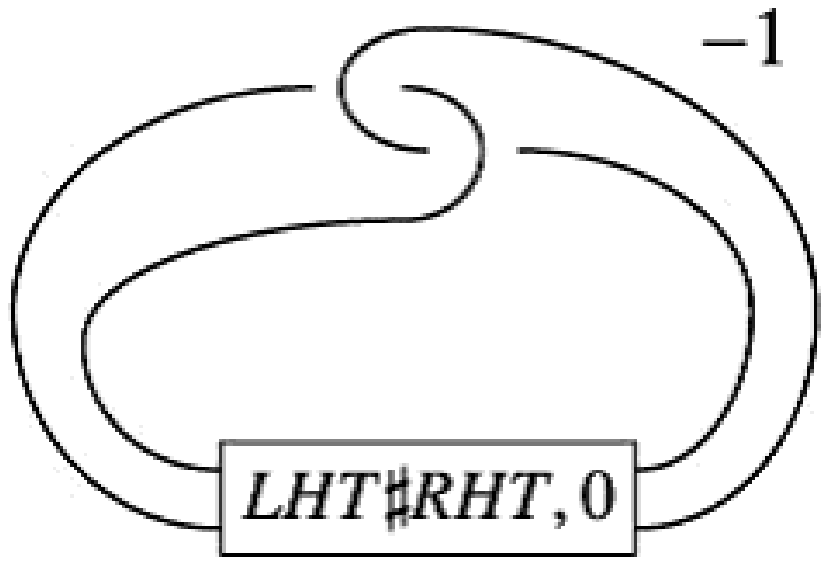}
  \caption{$V_0(RHT)$}
  \label{V0RHT}
 \end{minipage}
\end{figure}

\begin{figure}[H]
\begin{minipage}{0.4\hsize}
  \centering
   \includegraphics[height=8mm]{bdydiff.eps}
 \end{minipage}%
 \begin{minipage}{0.3\hsize}
  \centering
   \includegraphics[height=18mm]{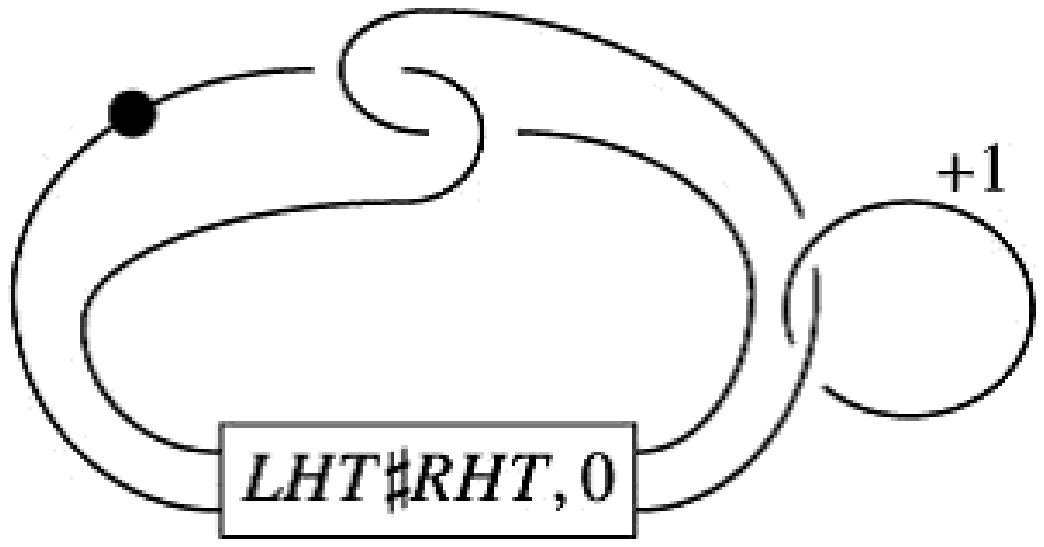}
  \caption{$W_0(RHT)$}
  \label{W0RHT}
 \end{minipage}
\end{figure}

Since $X_0(RHT)$ and $Z_0(RHT)$ have the same boundary, we have a homotopy $K3$ surface $Z_0(RHT) \cup_{\partial} (-W_0(RHT))$.
\begin{figure}[H]
\begin{minipage}{0.3\hsize}
  \centering
   \includegraphics[height=20mm]{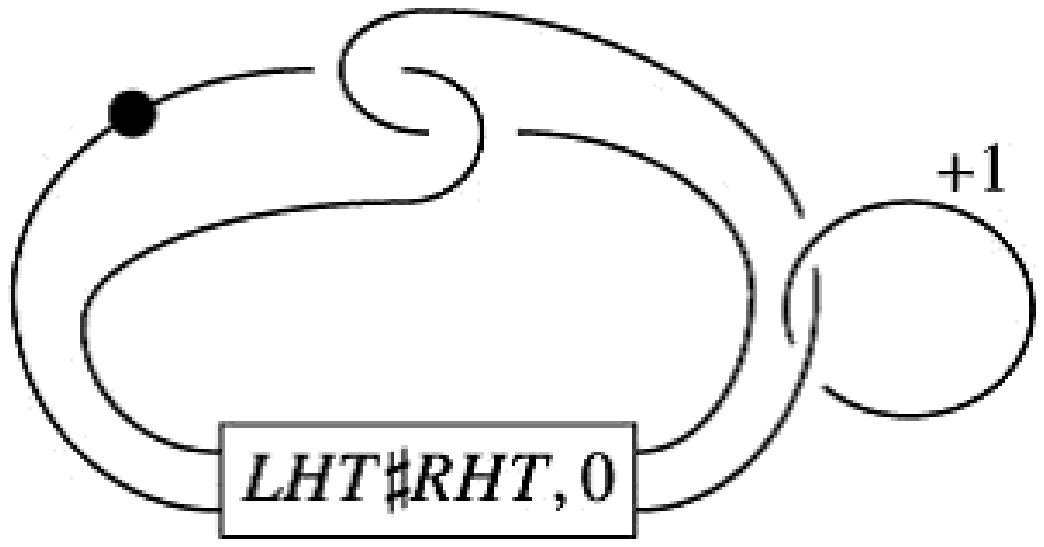}
   \caption{$W_0(RHT)$}
 \end{minipage}%
\begin{minipage}{0.08\hsize}
  \centering
   \includegraphics[height=9mm]{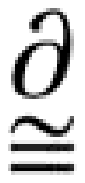}
 \end{minipage}%
 \begin{minipage}{0.6\hsize}
  \centering
   \includegraphics[height=60mm]{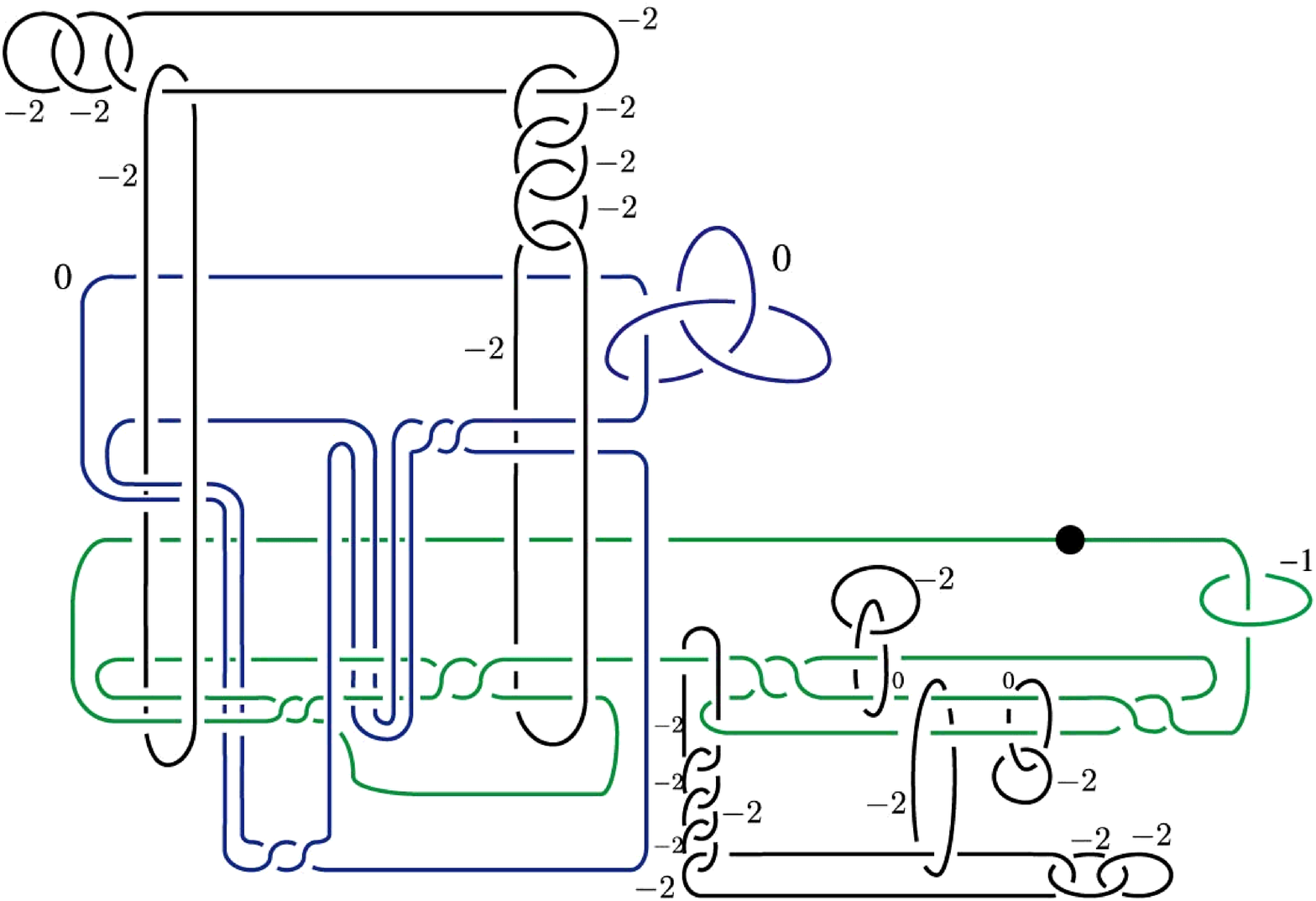}
  \caption{$Z_0(RHT)$}
 \end{minipage}
\end{figure}

\begin{figure}[H]
 \centering
  \includegraphics[height=65mm]{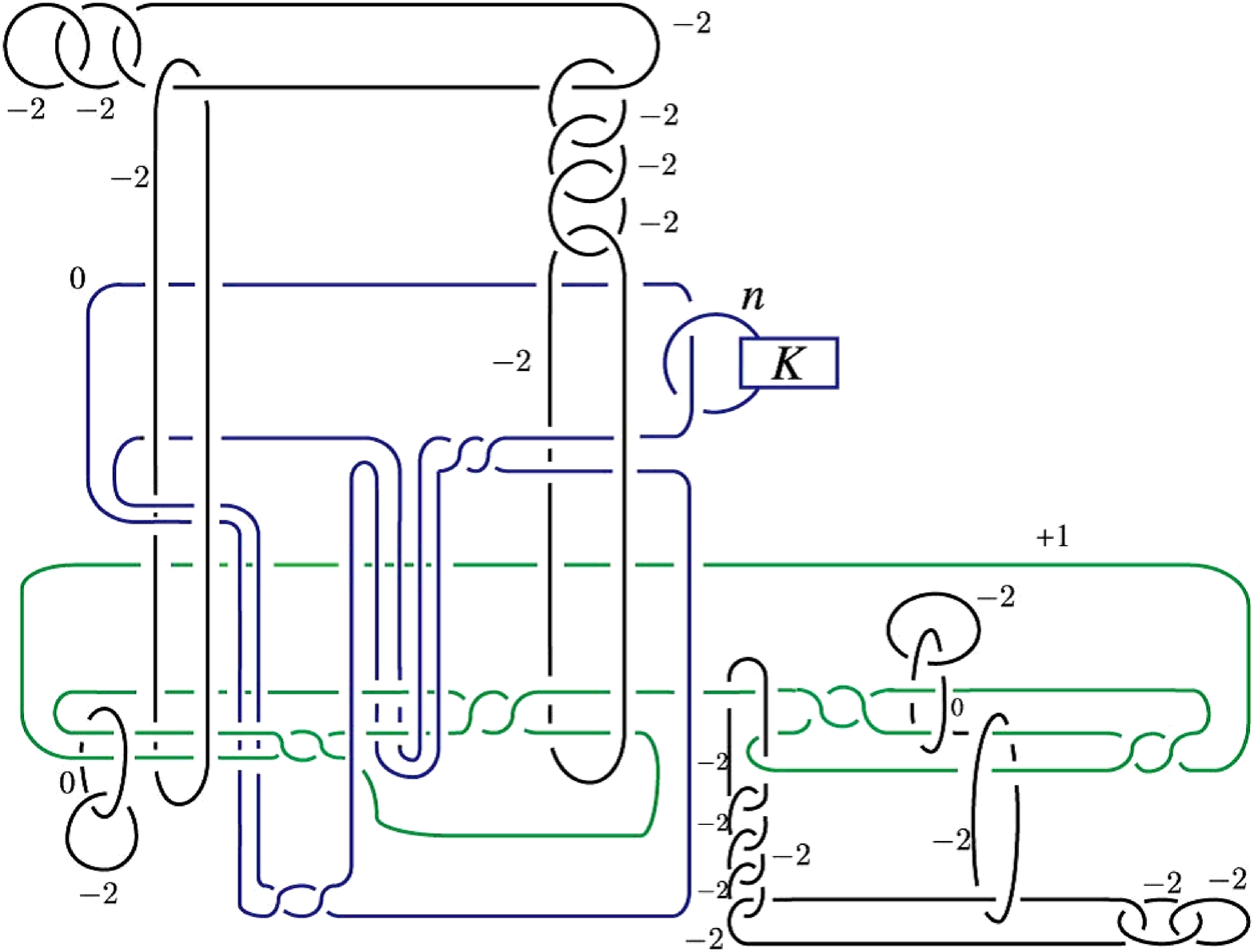}
 \caption{$Y'_n (K)$}
 \label{Y'nK}
\end{figure}

(c). Let $Y'_n(K)$ be the $4$-dimensional handlebody represented by \figref{Y'nK} with intersection form $2E_8\oplus2\displaystyle\bigl(\begin{smallmatrix} 0 & 1 \\ 1 & 0 \end{smallmatrix} \bigr)\oplus\displaystyle\bigl(\begin{smallmatrix} 0 & 1 \\ 1 & n \end{smallmatrix} \bigr)\oplus1$. Note that $Y_n(K)$ of \figref{YnK} and $Y'_n(K)$ have the same boundary.

\begin{figure}[H]
 \centering
  \includegraphics[height=65mm]{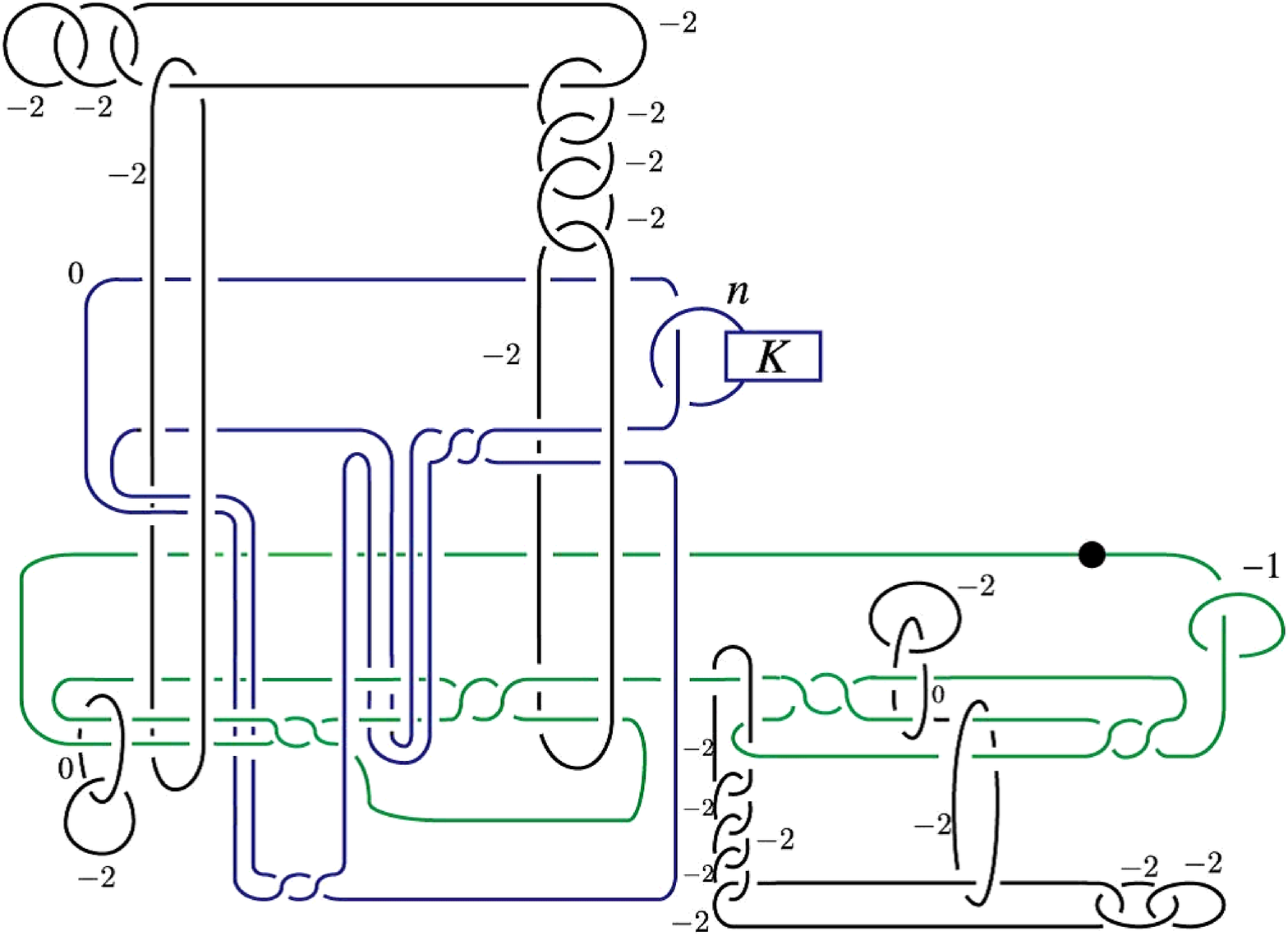}
 \caption{$Z'_n (K)$}
 \label{Z'nK}
\end{figure}

Since the $+1$-framed knot in \figref{Y'nK} is a slice knot (the connected sum of two figure eight knots), we can blow it down. Then we have a smooth simply connected $4$-dimensional handlebody $Z'_n(K)$ represented by \figref{Z'nK} with intersection form  $2E_8\oplus2\displaystyle\bigl(\begin{smallmatrix} 0 & 1 \\ 1 & 0 \end{smallmatrix} \bigr)\oplus\displaystyle\bigl(\begin{smallmatrix} 0 & 1 \\ 1 & n \end{smallmatrix} \bigr)$. Note that $Y'_n(K)$ and $Z'_n(K)$ have the same boundary.
Therefore we have homotopy $K3$ surfaces $Z'_0(RHT) \cup_{\partial} (-W'_0(RHT))$ and $Z'_0(RHT) \cup_{\partial} (-W_0(RHT))$.

\section{Appendix: On homology $3$-spheres which bound a contractible $4$-manifold}\label{sec:app}
\begin{figure}[H]
 \centering
  \includegraphics[height=28mm]{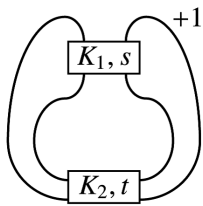}
  \caption{$V_{K_1, s}(K_2, t)$}
  \label{D1}
\end{figure}

We show a series of homology $3$-spheres which bound a contractible $4$-manifold.
Let $V_{K_1, s}(K_2, t)$ be the $4$-dimensional handlebody represented by \figref{D1} and $M_{K_1, s}(K_2, t)$ be $\partial(V_{K_1, s}(K_2, t))$, where $K_1$ and $K_2$ are knots and $s, t\in \Z$. If $K_1$ (or $K_2$) is an unknot $U$, a box in \figref{D1} represents the $s$ full twists (or $t$ full twists). Note that $M_{K_1, s}(K_2, t)$ is a homology $3$-sphere.

Let $T_{s, s+1}$ be a $(s, s+1)$-torus knot and $T_{1, 2}$ be an unknot $U$, where $s$ is a positive integer. We choose an orientation of $S^3$ such that $T_{2, 3}$ is $RHT$.
\begin{Theorem}\label{thm:appM}
The homology $3$-sphere $M_{T_{s, s+1}, s(s+1)}(K, n)$ bounds a contractible $4$-manifold, where $K$ is any knot and $n$ is any integer.
\begin{figure}[H]
 \centering
  \includegraphics[height=30mm]{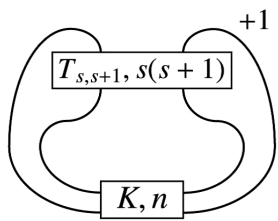}
 \caption{$V_{T_{s, s+1}, s(s+1)}(K, n)$}
 \label{VTKn}
\end{figure}
\end{Theorem}

\begin{example}
If $s$ is equal to $1$, $T_{1, 2}$ is an unknot $U$. The homology $3$-sphere $M_{U, 2}(K, n)$ bounds a contractible $4$-manifold, where $K$ is any knot and $n$ is any integer.
\end{example}
\begin{figure}[H]
 \centering
  \includegraphics[height=26mm]{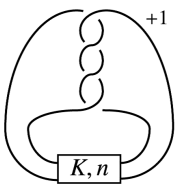}
  \caption{$V_{U, 2}(K, n)$}
  \label{VUKn}
\end{figure}

To prove \thmref{thm:appM}, we show the following \lemref{lem:app}.

\begin{Lemma}\label{lem:app}
The $4$-dimensional handlebodies represented by Figures~\ref{XrD1K3},~\ref{XtD2K2}~and~\ref{XsD3K1} have the same boundary, where $K_1, K_2$ and $K_3$ are knots and $s, t, r\in \Z$.
\begin{figure}[H]
 \begin{minipage}{0.4\hsize}
  \centering
   \includegraphics[height=35mm]{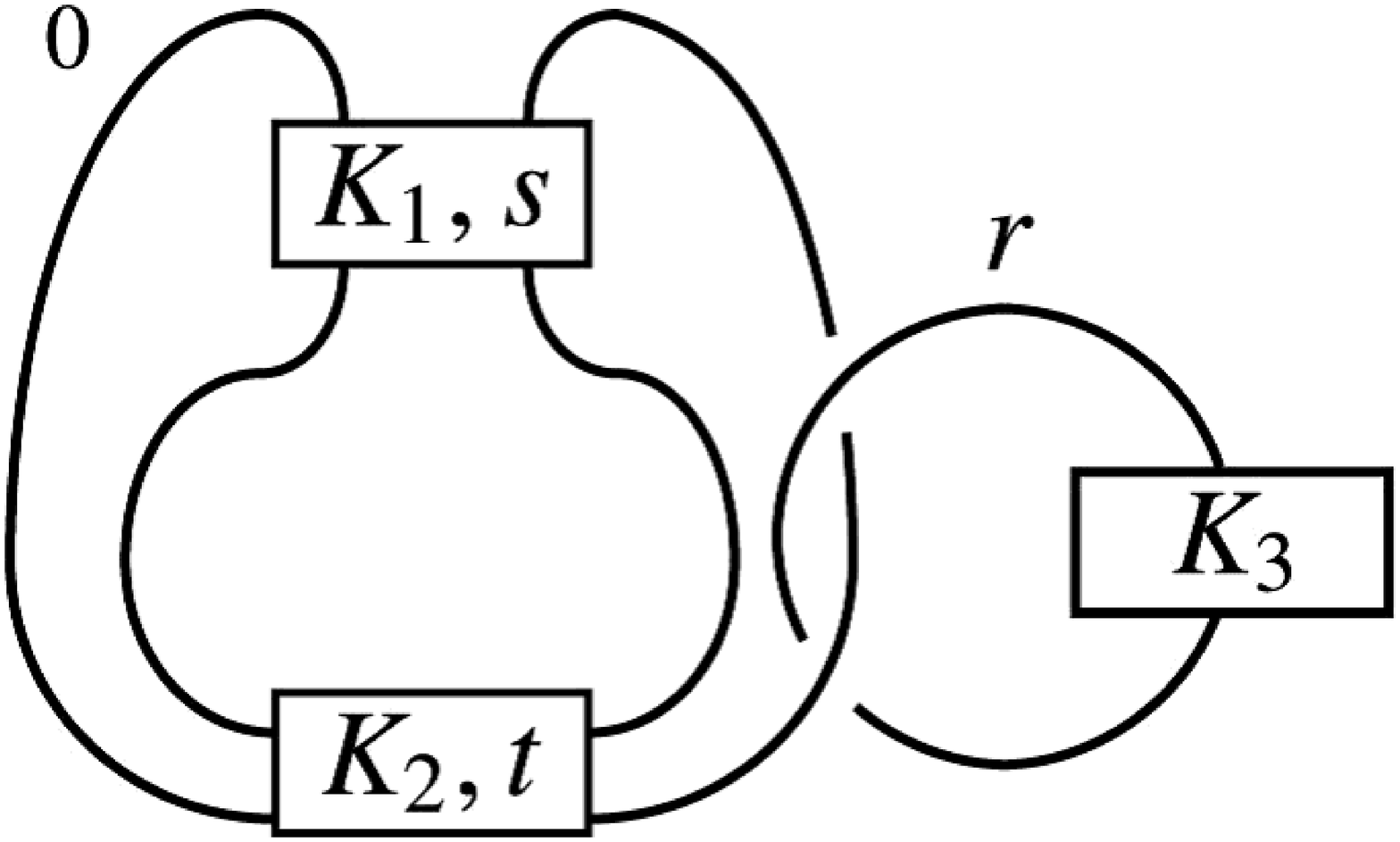}
 \caption{}
 \label{XrD1K3}
 \end{minipage}%
 \begin{minipage}{0.2\hsize}
  \centering
   \includegraphics[height=10mm]{bdydiff.eps}
 \end{minipage}%
 \begin{minipage}{0.4\hsize}
  \centering
   \includegraphics[height=35mm]{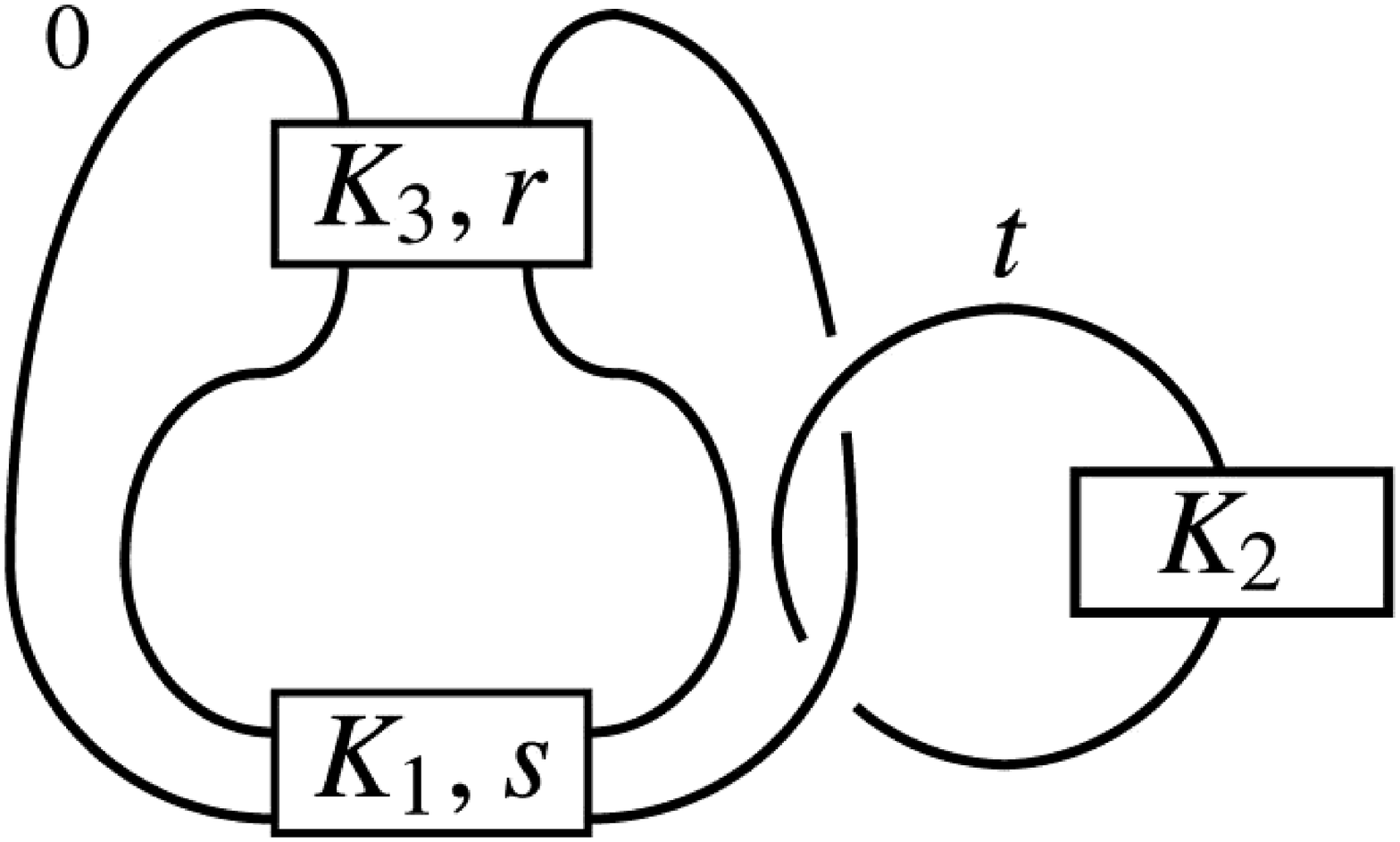}
  \caption{}
  \label{XtD2K2}
 \end{minipage}
\end{figure}
\begin{figure}[H]
\begin{minipage}{0.4\hsize}
  \centering
   \includegraphics[height=10mm]{bdydiff.eps}
 \end{minipage}%
 \begin{minipage}{0.3\hsize}
  \centering
   \includegraphics[height=35mm]{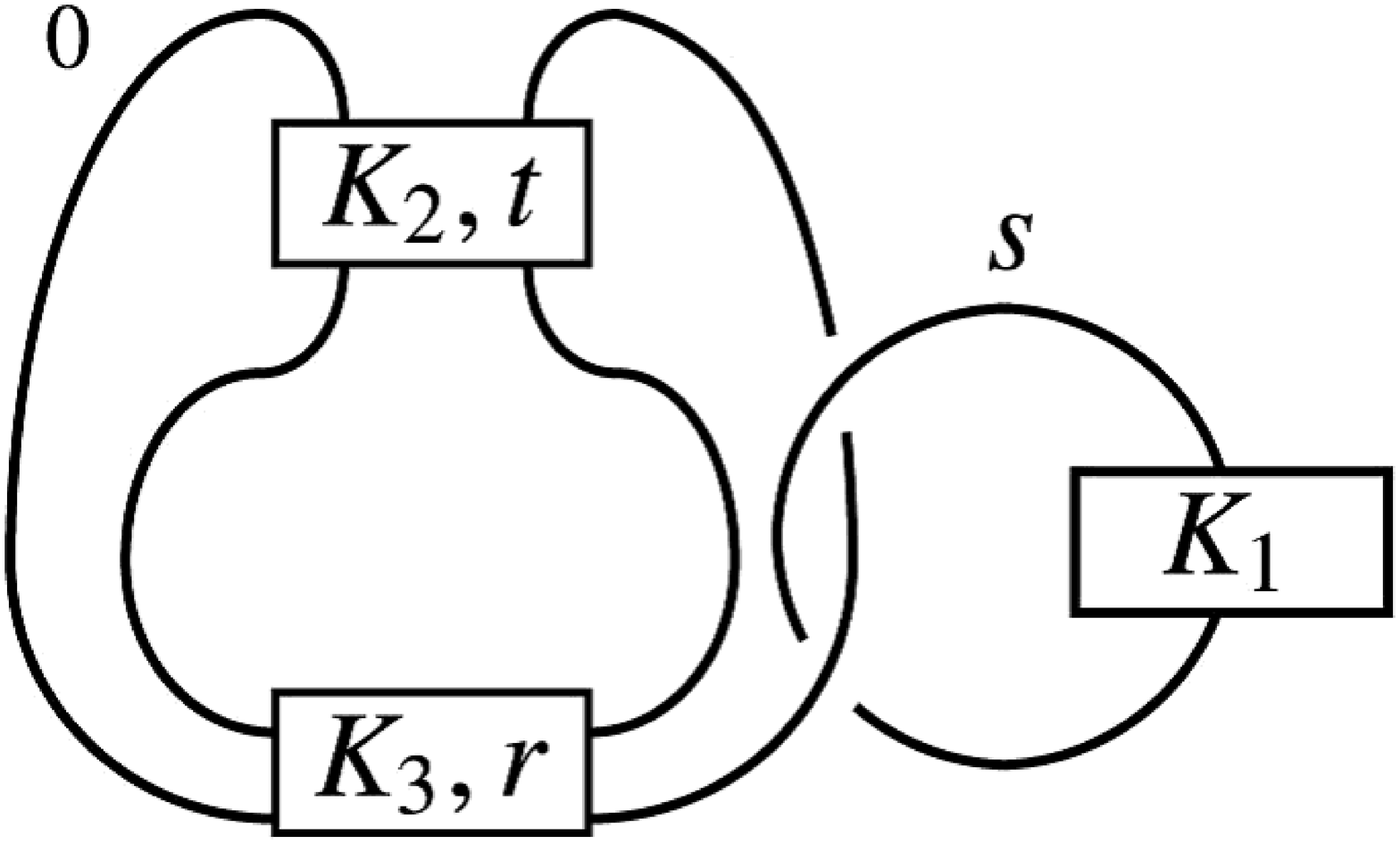}
  \caption{}
  \label{XsD3K1}
 \end{minipage}
\end{figure}
\end{Lemma}
\proof
We show \lemref{lem:app} by the following handle calculus:
\begin{figure}[H]
 \begin{minipage}{0.4\hsize}
  \centering
   \includegraphics[height=35mm]{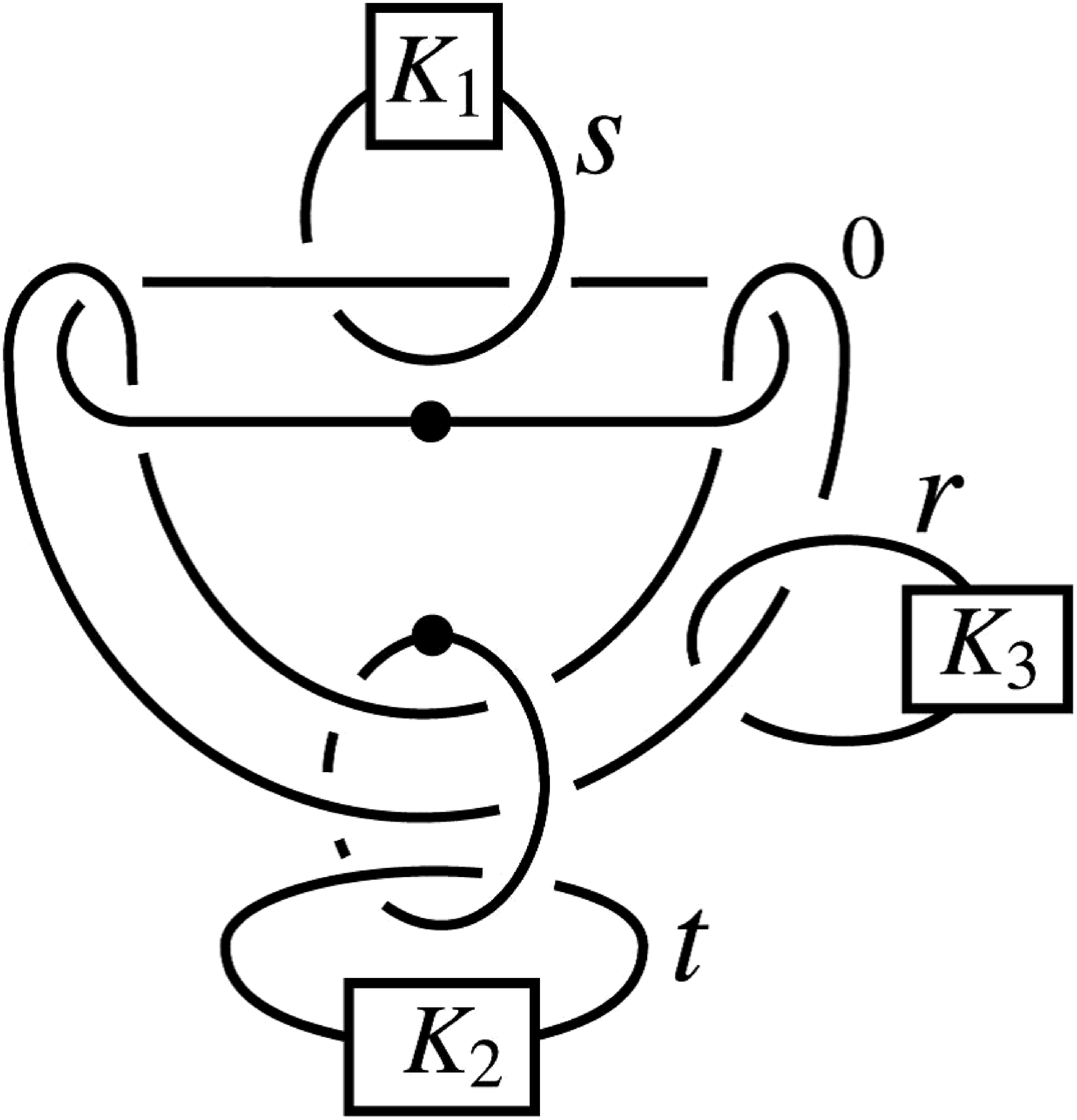}
  \caption{}
 \end{minipage}%
 \begin{minipage}{0.2\hsize}
  \centering
   \includegraphics[height=10mm]{diff.eps}
 \end{minipage}%
 \begin{minipage}{0.4\hsize}
  \centering
   \includegraphics[height=30mm]{XrD1K3.eps}
  \caption{}
 \end{minipage}
\end{figure}

\begin{figure}[H]
\begin{minipage}{0.1\hsize}
  \centering
   \includegraphics[height=20mm]{surgery4.eps}
 \end{minipage}%
 \begin{minipage}{0.3\hsize}
  \centering
   \includegraphics[height=35mm]{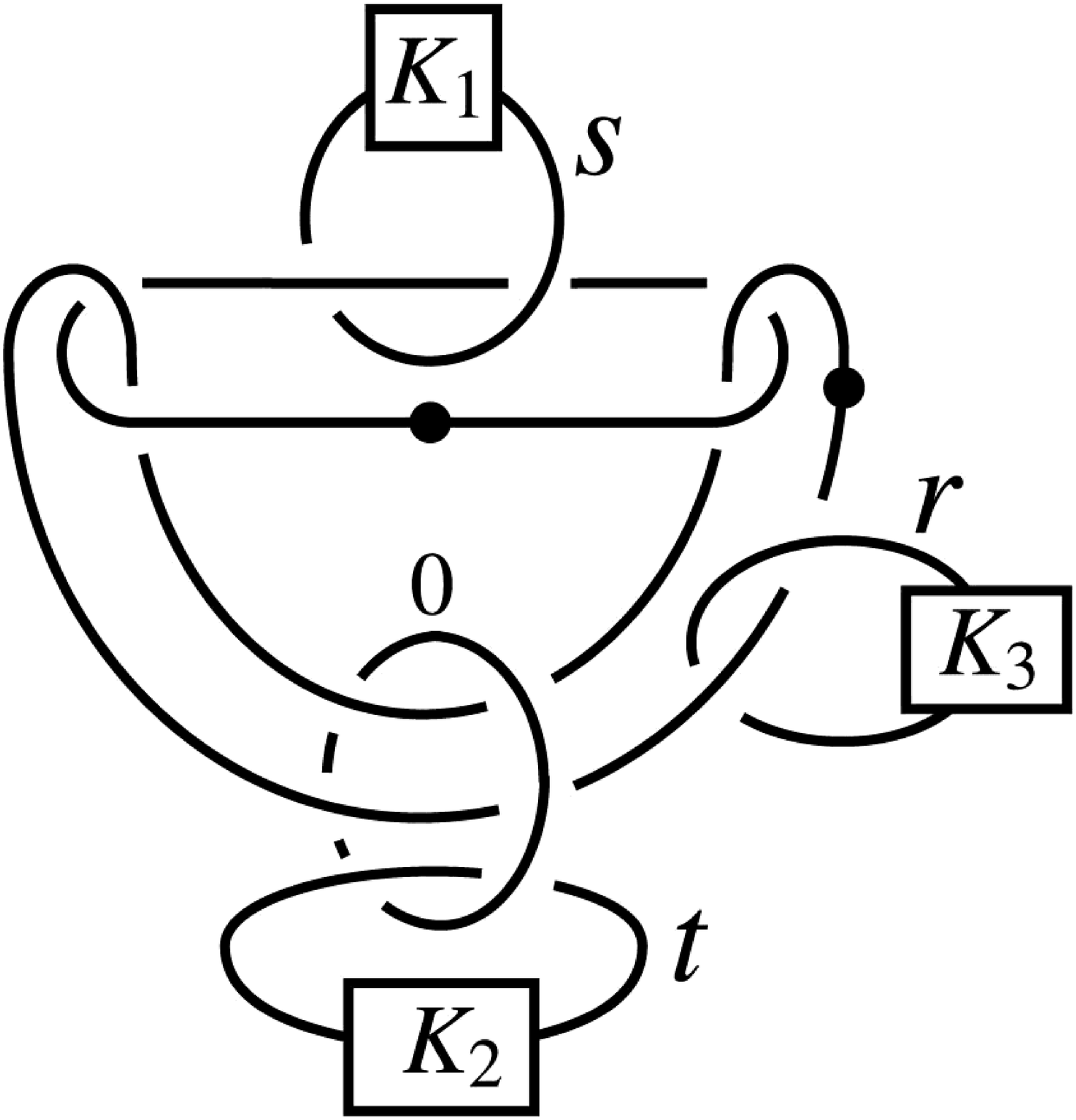}
   \caption{}
 \end{minipage}%
\begin{minipage}{0.2\hsize}
  \centering
   \includegraphics[height=10mm]{diff.eps}
 \end{minipage}%
 \begin{minipage}{0.4\hsize}
  \centering
   \includegraphics[height=30mm]{XtD2K2.eps}
  \caption{}
 \end{minipage}
\end{figure}

\begin{figure}[H]
\begin{minipage}{0.1\hsize}
  \centering
   \includegraphics[height=20mm]{surgery4.eps}
 \end{minipage}%
 \begin{minipage}{0.3\hsize}
  \centering
   \includegraphics[height=35mm]{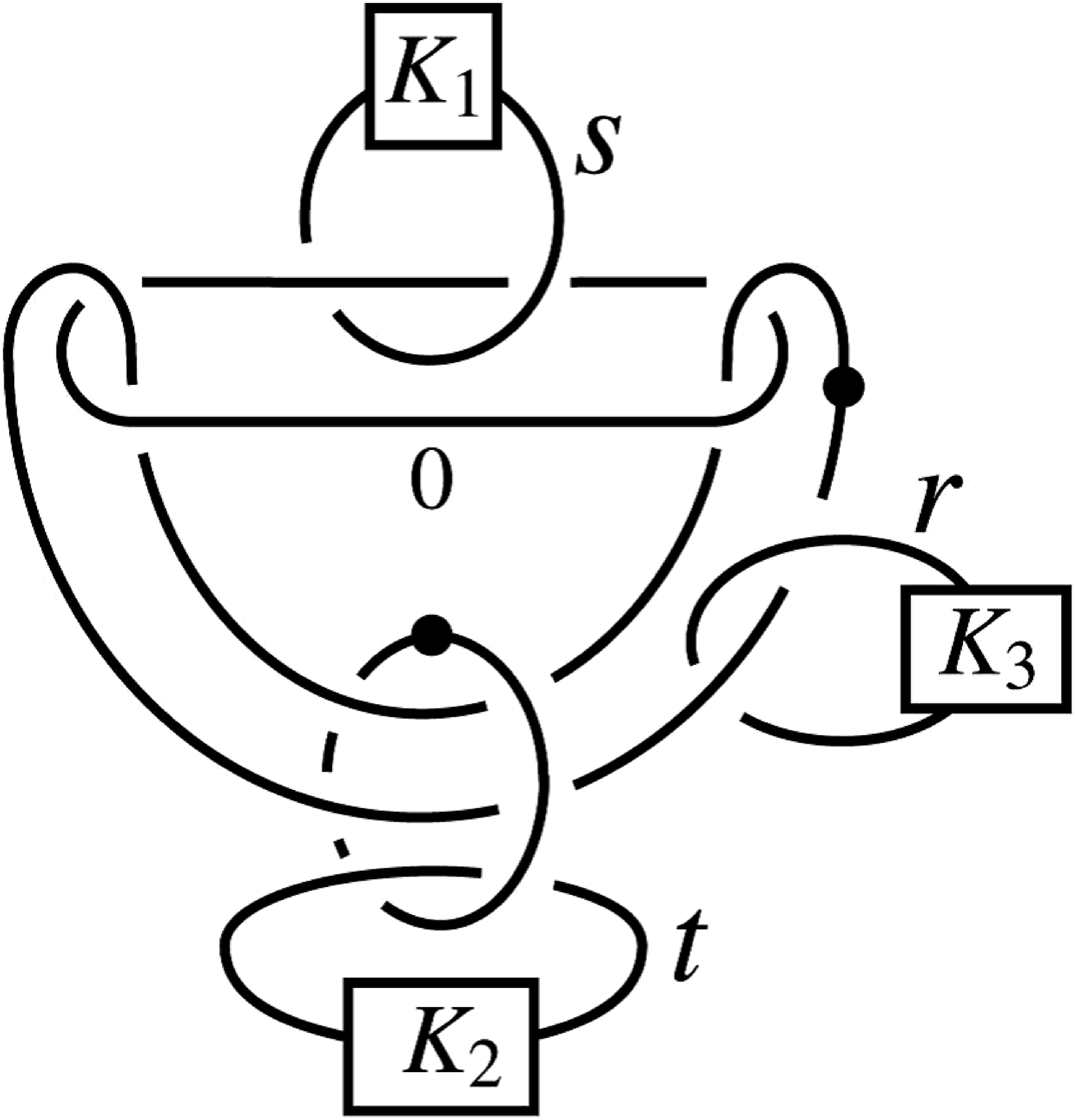}
   \caption{}
 \end{minipage}%
\begin{minipage}{0.2\hsize}
  \centering
   \includegraphics[height=10mm]{diff.eps}
 \end{minipage}%
 \begin{minipage}{0.4\hsize}
  \centering
   \includegraphics[height=30mm]{XsD3K1.eps}
  \caption{}
 \end{minipage}
\end{figure}
\endproof

We will prove \thmref{thm:appM}.

\proof[Proof of \thmref{thm:appM}]
We prove \thmref{thm:appM} by the following handle calculus:
\begin{figure}[H]
 \begin{minipage}{0.4\hsize}
  \centering
   \includegraphics[height=20mm]{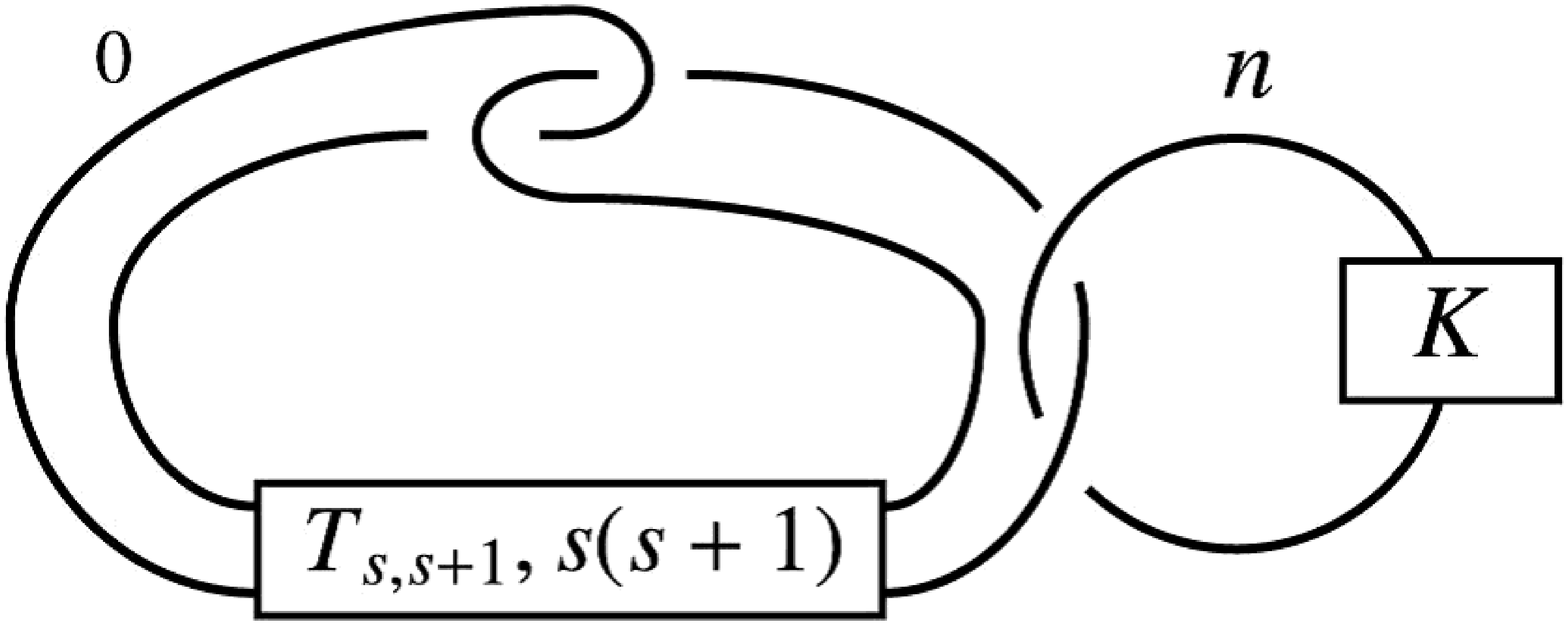}
  \caption{}
  \label{pVTKn1}
 \end{minipage}%
 \begin{minipage}{0.2\hsize}
  \centering
   \includegraphics[height=13mm]{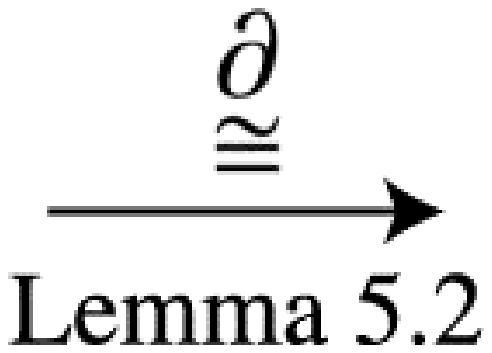}
 \end{minipage}%
 \begin{minipage}{0.4\hsize}
  \centering
   \includegraphics[height=20mm]{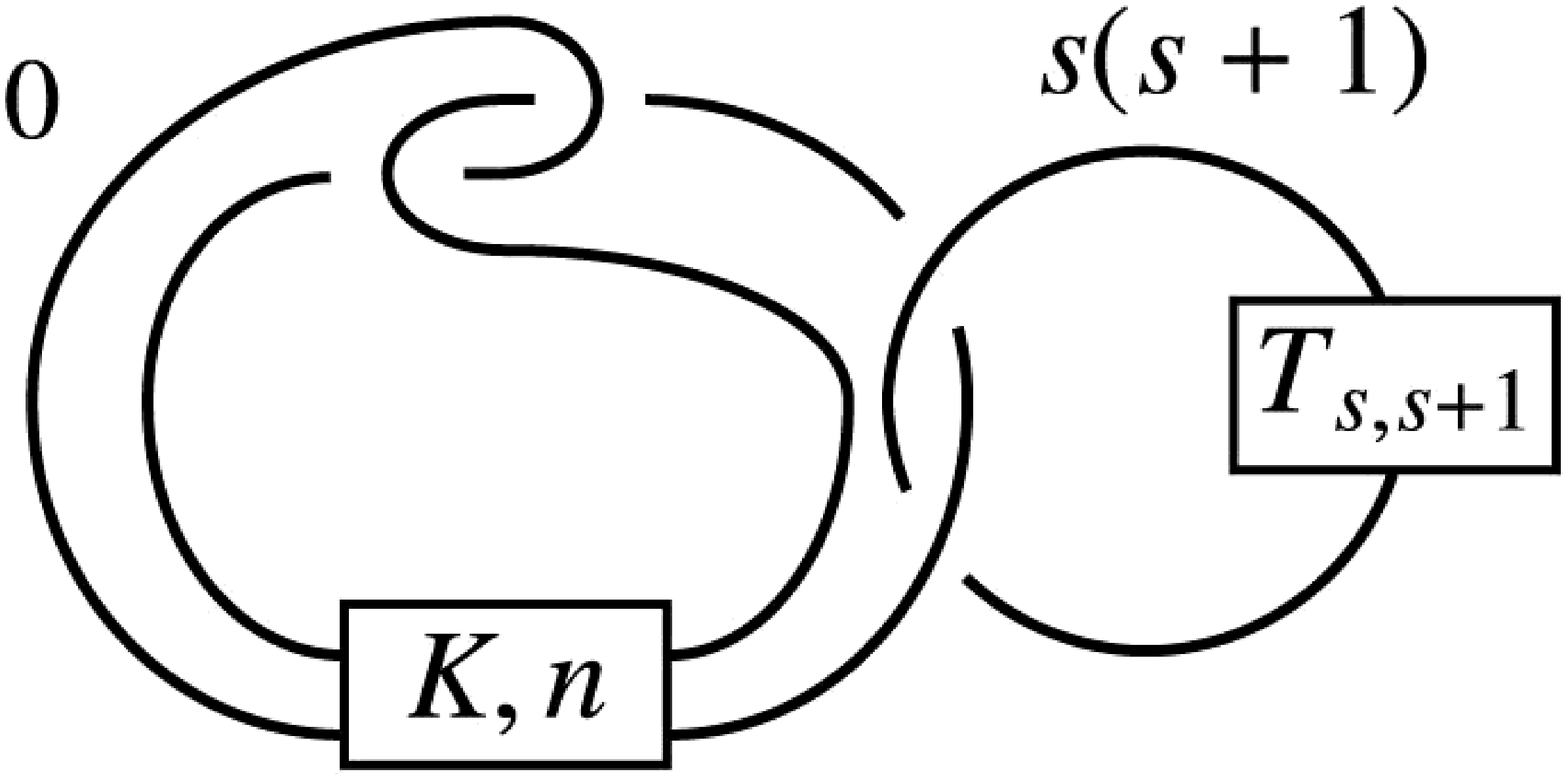}
  \caption{}
  \label{pVTKn2}
 \end{minipage}
\end{figure}

\begin{figure}[H]
\begin{minipage}{0.1\hsize}
  \centering
   \includegraphics[height=13mm]{bdydiff3.eps}
 \end{minipage}%
 \begin{minipage}{0.4\hsize}
  \centering
   \includegraphics[height=25mm]{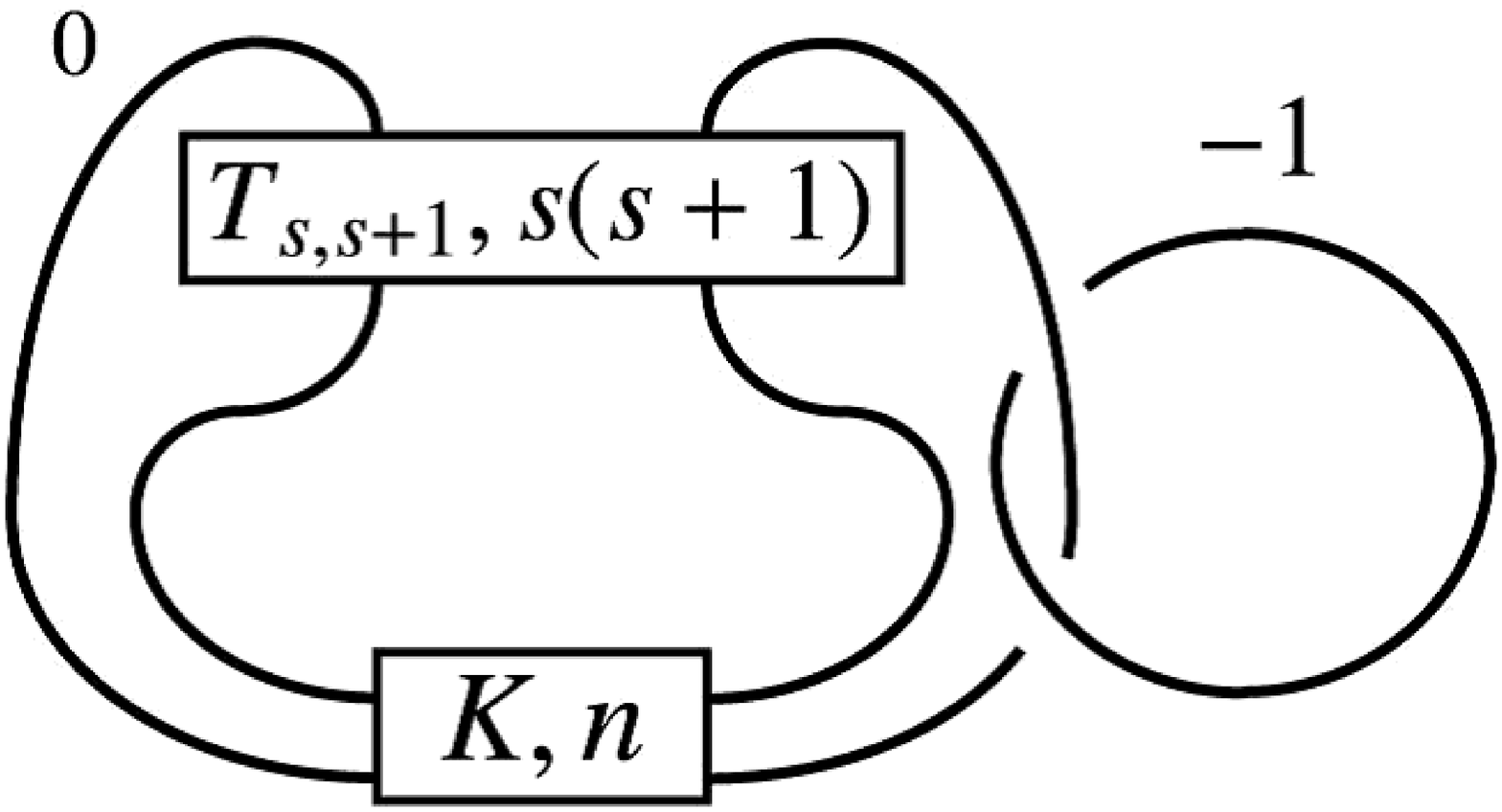}
   \caption{}
 \end{minipage}%
\begin{minipage}{0.1\hsize}
  \centering
   \includegraphics[height=14mm]{blowdown.eps}
 \end{minipage}%
 \begin{minipage}{0.4\hsize}
  \centering
   \includegraphics[height=25mm]{VTKn.eps}
  \caption{$V_{T_{s, s+1}, s(s+1)}(K, n)$}
 \end{minipage}
\end{figure}
\begin{figure}[H]
 \begin{minipage}{0.4\hsize}
  \centering
   \includegraphics[height=20mm]{proofVTKn1.eps}
  \caption{}
 \end{minipage}%
 \begin{minipage}{0.2\hsize}
  \centering
   \includegraphics[height=10mm]{bdydiff.eps}
 \end{minipage}%
 \begin{minipage}{0.4\hsize}
  \centering
   \includegraphics[height=20mm]{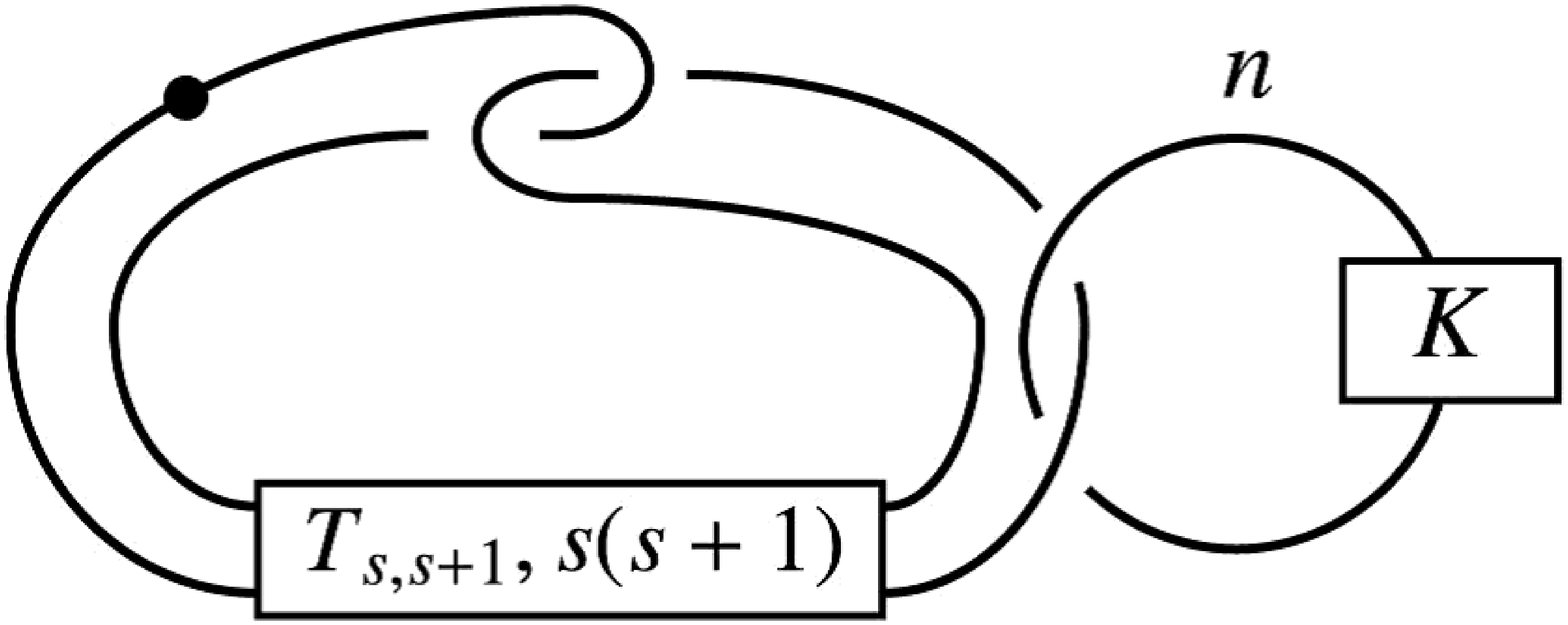}
  \caption{$W$}
  \label{pVTKn4}
 \end{minipage}
\end{figure}

Note that the left side knot in \figref{pVTKn1} is $D_+(T_{s, s+1}, s(s+1))$. If $s$ is equal to $1$, it is known that $D_+(U, 2)$ is a slice knot. If $s\geq2$, Litherland \cite{L} remarks that $D_+(T_{s, s+1}, s(s+1))$ is a slice knot. Therefore we have a smooth $S^2$ with self intersection $0$ in the $4$-dimensional handlebody represented by \figref{pVTKn1}. By performing surgery on the $S^2$, we have a contractible $4$-manifold $W$  represented by \figref{pVTKn4}.
\endproof

\begin{Remark}
By \thmref{thm:appM}, the $4$-manifold $V_{T_{s, s+1}, s(s+1)}(K, n) \cup_{\partial} (-W)$ is a homotppy $\CP^2$.
\end{Remark}

\begin{Remark}
In the case $s$ is equal to $1$, we could prove that this $4$-manifold is diffeomorphic to $\CP^2$.
\end{Remark}

\begin{Remark}
Let $t$ be a positive integer. If $K$ is $T_{t, t+1}$ and $n$ is equal to $t(t+1)$ in \figref{pVTKn1}, we have a smooth $S^2$ with self intersection $0$ in the $4$-dimensional handlebody represented by \figref{pVTKn2}. By performing surgery on the $S^2$, we have a contractible $4$-manifold represented by \figref{apprem}. By the proof of \thmref{thm:appM}, the $4$-dimensional handlebodies represented by Figures~\ref{apprem} and \ref{apprem2} have the same boundary. Therefore by gluing these two contractible $4$-manifolds, we get a homotopy $S^4$. We do not know whether or not this homotopy $S^4$ is diffeomorphic to $S^4$.
\begin{figure}[H]
 \begin{minipage}{0.4\hsize}
  \centering
   \includegraphics[height=20mm]{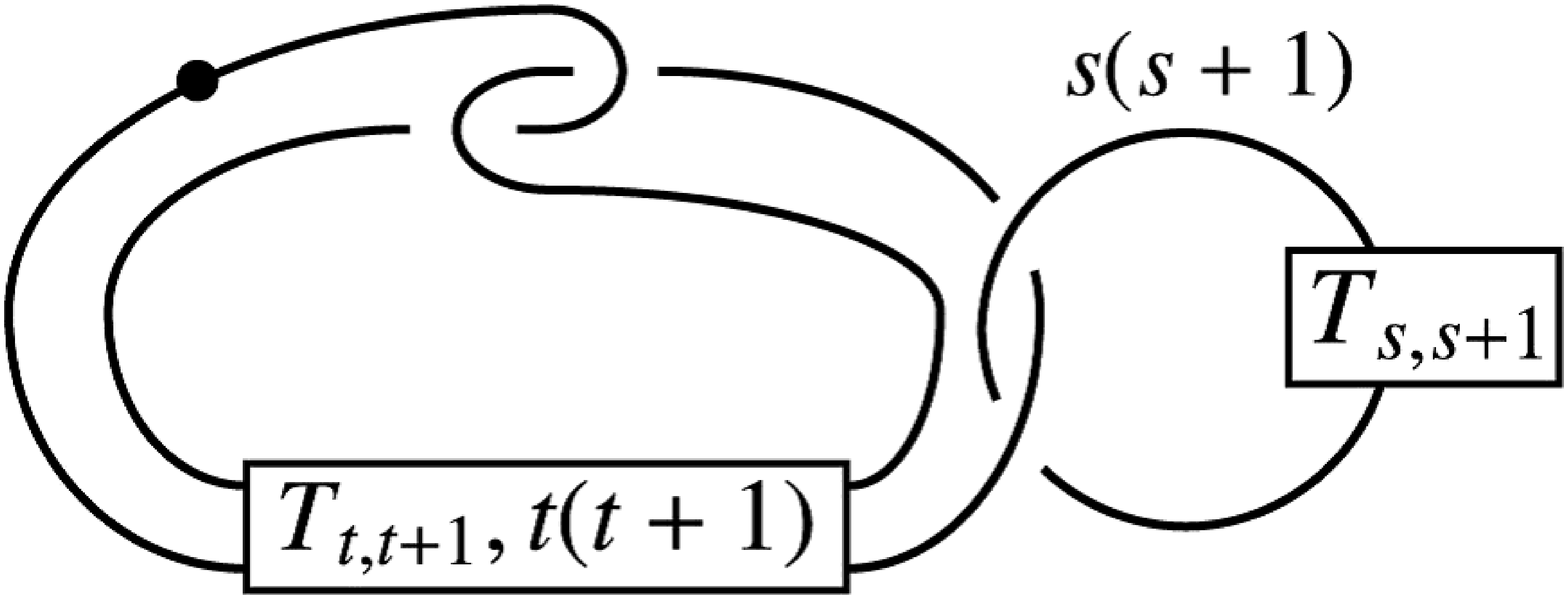}
  \caption{}
  \label{apprem}
 \end{minipage}%
 \begin{minipage}{0.2\hsize}
  \centering
   \includegraphics[height=10mm]{bdydiff.eps}
 \end{minipage}%
 \begin{minipage}{0.4\hsize}
  \centering
   \includegraphics[height=20mm]{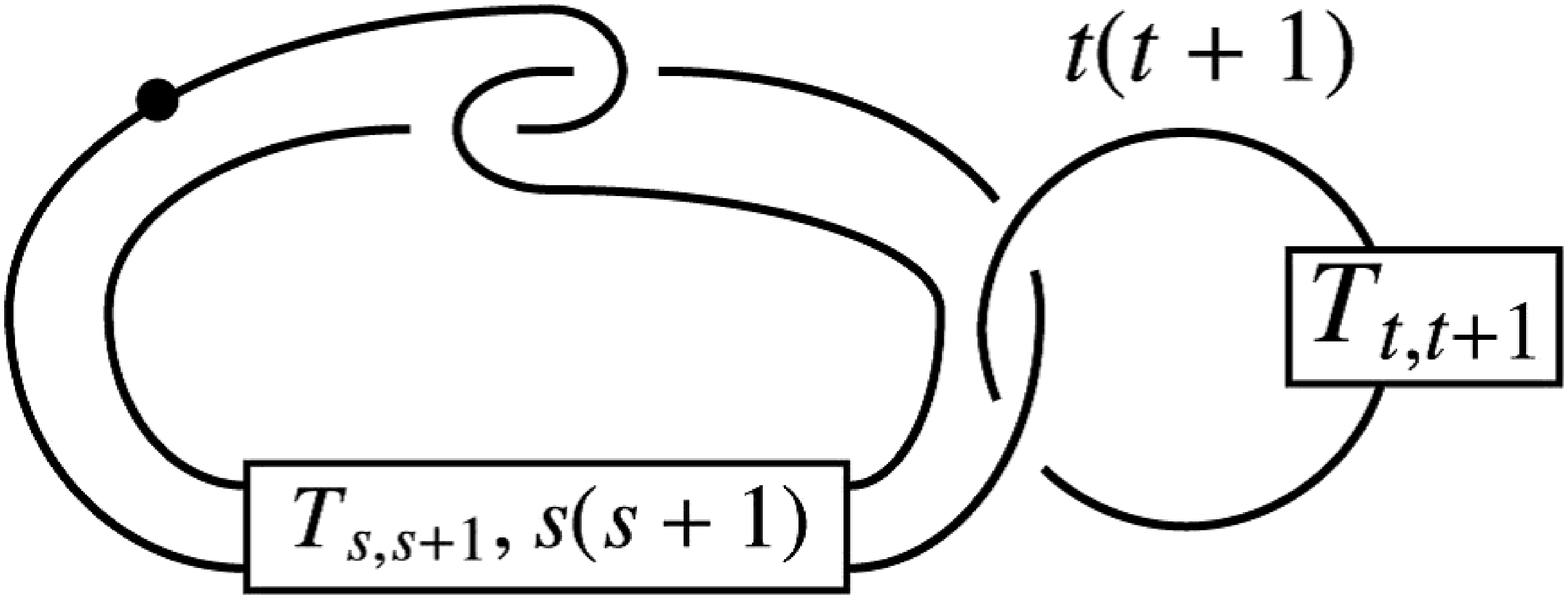}
  \caption{}
  \label{apprem2}
 \end{minipage}
\end{figure}
\end{Remark}

\end{document}